\documentclass[a4paper,11pt]{article}
\usepackage{amssymb}
\usepackage{amsfonts}
\usepackage{amsmath,tabu}
\usepackage{amsthm}
\usepackage{cite}
\usepackage{amsmath, upgreek}
\usepackage{amssymb}
\usepackage{amscd}
\usepackage{xspace}
\usepackage{verbatim}
\usepackage{xcolor}
\usepackage{mathrsfs}
\usepackage{color}

\DeclareMathOperator*{\esssup}{ess\,sup}

\renewcommand{\j}{{\mathbf{j}}}
\newcommand{\s}{s}

\newcommand{\C}{\mathbb{C}}
\newcommand{\R}{\mathbb{R}}

\newcommand{\N}{\mathbb{N}}

\newcommand{\cuad}{{\sqcap\kern-.68em\sqcup}}

\newcommand{\norm}[1]{\|#1\|}

\newcommand{\bel}[1]{\begin{equation}\label{#1}}
\newcommand{\ee}{\end{equation}}
\numberwithin{equation}{section}
\newtheorem{theorem}{Theorem}[section]
\newtheorem{definition}[theorem]{Definition}
\newtheorem{proposition}[theorem]{Proposition}

\newtheorem{lemma}[theorem]{Lemma}
\newtheorem{corollary}[theorem]{Corollary}
\newtheorem{remark}[theorem]{Remark}
\newcommand{\bremark}{\begin{remark} \em}
\newcommand{\eremark}{\end{remark} }

\newcommand{\nind}{\noindent}

\newcommand{\cC}{{\mathcal C}}
\newcommand{\cD}{{\mathcal D}}
\newcommand{\cE}{{\mathcal E}}
\newcommand{\cF}{{\mathcal F}}

\newcommand{\cH}{{\mathbb H}}

\newcommand{\cL}{{\mathcal L}}

\newcommand{\cP}{{\mathcal P}}
\newcommand{\cQ}{{\mathcal Q}}

\newcommand{\cS}{{\mathcal S}}

\newcommand{\dd}{{\,\mathrm{d}}} 
\newcommand{\bbX}{{\,\mathbb{X}}} 
\newcommand{\lnlap}{\cL_{\text{\tiny $\Delta \,$}}\!}
\newcommand{\bth}[1]{\def\name{Theorem}
\begin{sub}\label{t:#1}}
\newcommand{\blemma}[1]{\def\name{Lemma}
\begin{sub}\label{l:#1}}
\newcommand{\bcor}[1]{\def\name{Corollary}
\begin{sub}\label{c:#1}}
\newcommand{\bdef}[1]{\def\name{Definition}
\begin{sub}\label{d:#1}}
\newcommand{\bprop}[1]{\def\name{Proposition}
\begin{sub}\label{p:#1}}
\newcommand{\dsps}{\displaystyle}
\newcommand{\ti}{\times}

\newcommand{\BBR}{\eqref}

\newcommand{\BA}{
\begin{array}}
\newcommand{\EA}{
\end{array}}

\def\ga{\alpha}     \def\gb{\beta}       
       \def\gd{\delta}      \def\ge{\epsilon}
\def\gth{\theta}                         
\def\gf{\varphi}           
            \def\gl{\lambda}
\def\gm{\mu}                 \def\gp{\pi}
    \def\gr{\rho}        
\def\gs{\sigma}       \def\gt{\tau}
      \def\gw{\omega}
                \def\gz{\zeta}
\def\Gg{\Gamma}     \def\Gd{\Delta}

\def\Gw{\Omega}              

\def\CS{{\mathcal S}}      
      \def\CP{{\mathcal P}}
      \def\CC{{\mathcal C}}
\def\CD{{\mathcal D}}   \def\CE{{\mathcal E}}   \def\CF{{\mathcal F}}
      
      \def\CL{{\mathcal L}}

   \def\BBH {\mathbb H}    
       
   \def\BBN {\mathbb N}    
   \def\BBR {\mathbb R}

\DeclareMathOperator{\dist}{\rm dist}

\usepackage{cite}
\begin{document}

\begin{center}{\bf  \large   The Cauchy problem associated to the logarithmic Laplacian\\[2mm]
with an application to the fundamental solution 
 }\bigskip

\medskip%\medskip
 
{\small
 Huyuan Chen\footnote{chenhuyuan@yeah.net}

\medskip
Department of Mathematics, Jiangxi Normal University,\\
Nanchang, Jiangxi 330022, PR China \\[6pt]

 The University of Sydney, School of Mathematics and Statistics, NSW 2006, Australia\\[14pt]
 
 Laurent V\'eron\footnote{veronl@univ-tours.fr}\medskip

Institut Denis Poisson, CNRS-UMR 7013\\ Universit\'e de Tours,  37200 Tours, France \\[14pt]

  }%\medskip

 %\\[6mm] 

 \begin{abstract}
Let  $\lnlap$ be the logarithmic Laplacian operator  with Fourier symbol $2\ln |\zeta|$, we study the expression of  the diffusion kernel  which is associated to the equation
$$\partial_tu+ \lnlap u=0 \ \  {\rm in}\ \,  (0,\tfrac N2) \times \R^N,\quad\quad  u(0,\cdot)=0\ \ {\rm in}\  \, \R^N\setminus \{0\}.$$
We apply our results to give a classification of the solutions of
$$\left\{ \arraycolsep=1pt
\begin{array}{lll}
\displaystyle \partial_tu+ \lnlap u=0\quad \  &{\rm in}\ \   (0,T)\times \R^N,\\[2.5mm]
 \phantom{ \lnlap \ \,   }
\displaystyle   u(0,\cdot)=f\quad \ &{\rm{in}}\  \   \R^N
\end{array}
\right.
$$
and obtain an expression of  the fundamental solution of the associated stationary equation in $\R^N$,
and of the fundamental solution $u$ in a bounded domain, i.e.
$\lnlap u=k\delta_0$ in the sense of distributions in $\Omega$, such that $u=0$ in  $\R^N\setminus\Omega.$

 \end{abstract}

\end{center}
 \tableofcontents \vspace{1mm}
 \vspace{1mm}

\vspace{2mm}

%%%%%%%%%%%%%%%%%%%%%%%%%%%%%%%%%%%%%%%%%%%%%%%%%%%%%%%%%%%%%%%%%%%%%%%%%%%%%%%%%%%%%%%%%%%%%%%%%%%%%%%%%%%%%%%%%%%%%%%%%%
%%%%%%%%%%%%%%%%%%%%%%%%%%%%%%%%%%%%%%%%%%%%%%%%%%%%%%%%%%%%%%%%%%%%%%%%%%%%%%%%%%%%%%%%%%%%%%%%%%%%%%%%%%%%%%%%%%%%%%%%%%

\setcounter{equation}{0}
\section{Introduction }

%\subsection{Cauchy problem}

The logarithmic Laplacian $\lnlap$ in $\R^N$ ($N\geq 1$)  is defined by the expression
$$
 \lnlap  u(x) = c_{_N} \int_{\R^N  } \frac{ u(x)1_{B_1(x)}(y)-u(y)}{|x-y|^{N} } \dd y + \rho_N u(x),
   $$
where
$$
c_{_N}:= \pi^{-  N/2}  \Gamma(\tfrac{N}2) = \frac{2}{\omega_{_{N}}}, \qquad \rho_N:=2 \ln 2 + \psi(\tfrac{N}{2}) +\psi(1),
$$
$\omega_{_{N}}:=|\mathbb{S}^{N}|=\displaystyle\int_{\mathbb{S}^{N}}dS$, $\mathbb{S}^{N}$ is the unit sphere in $\R^N$,  $\Gamma$ is the Gamma function,   $\psi = \frac{\Gamma'}{\Gamma}$ is the associated Digamma function. 
In this article our first aim is to study  the Cauchy problem 
\begin{equation}\label{eq 1.1}
\left\{ \arraycolsep=1pt
\begin{array}{lll}
\displaystyle \partial_tu+ \lnlap u=0\quad \  &{\rm in}\ \   (0,T)\times \R^N,\\[2.5mm]
 \phantom{ \lnlap \ \,   }
\displaystyle   u(0,\cdot)=f\quad \ &{\rm{in}}\  \ \,  \R^N,
\end{array}
\right.
\end{equation}
where $f$ is a measurable function in $\R^N$ and $T\in(0,\frac{N}{2})$.\smallskip

%%%%
The historical model for diffusion equation is the heat equation 
\begin{equation}\label{eq 1.1 lap} 
u_{t}(t, x)-\Delta u(t,x)=0, \quad(x, t) \in (0, T)  \times \mathbb{R}^{N}, 
\end{equation}
for which  Widder proved in \cite{W} the uniqueness among nonnegative classical solutions and provided  the representation 
$$u(t,x)=\frac{1}{(4 \pi t)^{\frac{N}{2}}} \int_{\mathbb{R}^{N}}  
e^{-\tfrac{|x-y|^{2}}{4 t}} u(0,y)\dd y.$$
By a classical solution $u$ we intend a function $u$ belonging to $\mathcal{C}\left([0, T) \times \mathbb{R}^{N} \right)$, such that $u_{t}, u_{x_{i} x_{j}} \in \mathcal{C}\left((0, T) \times \mathbb{R}^{N} \right)$ satisfying (\ref{eq 1.1 lap})
everywhere (see \cite{H,W1} and the references therein).

In recent years, there has been a renewed and increasing interest in the study of diffusion equations involving linear and nonlinear integro-differential operators, and this growing interest is justified both by important progresses made in
understanding nonlocal phenomena from a PDE or a probabilistic point of view, see e.g.
\cite{CR,CS0,CS1,CF,C,CV,CFQ,RS1,RS,Si,musina-nazarov,L,T} and the references therein, and by their wide range of applications.
Among nonlocal differential order operators the simplest examples are the fractional powers of the Laplacian, which exhibit  many phenomenological  properties.
Recall that for $\s\in(0,1)$, the fractional Laplacian
$(-\Delta)^\s$ can be written as a singular integral operator defined in the principle value sense
$$% \begin{equation}\label{fl 1}
 (-\Delta)^\s  u(x)=c_{N,\s} \lim_{\epsilon\to0^+} \int_{\R^N\setminus B_\epsilon(x) }\frac{u(x)-
u(y)}{|x-y|^{N+2\s}}  \dd y ,
$$%\end{equation}
where $c_{N,\s}=2^{2\s}\pi^{-\frac N2}\s\frac{\Gamma(\frac{N+2\s}2)}{\Gamma(1-\s)}>0$ is a normalized constant such that  for a function $u \in C^\infty_c(\R^N)$,   
\begin{equation*}
  \label{eq:Fourier-representation}
\mathcal{F}((-\Delta)^\s u)(\xi):=(2\pi)^{-\frac{N}{2}}\int_{\R^N}e^{-ix\cdot \xi}((-\Delta)^\s u)(x) \dd x= |\xi|^{2\s}\widehat u (\xi)\quad\ \text{for all $\xi \in \R^N$}.
\end{equation*}
Here and in the sequel both $\mathcal{F}$ and $\widehat \cdot$ denote the Fourier transform. 
A representation formula  for the {\it $s$-fractional diffusion
equation} 
 \begin{eqnarray}\label{eq 1.1 frac} 
u_{t}+(-\Delta)^{s} u=0 \quad  {\rm in }\ \,  (0, T) \times\mathbb{R}^N, 
 \end{eqnarray} 
was recently proved in \cite{BPS}. It is based upon the $s$-fractional diffusion kernel
 \begin{eqnarray}\label{eq-2}
\cP_{s}(t,x)= t^{-\frac{N}{2s}}  {\bf P}\left(  t^{-\frac{1}{2s}}x\right),
 \end{eqnarray}
where
$$%\begin{eqnarray}\label{eq-3}
{\bf P}(x):=\int_{\mathbb{R}^N} \mathrm{e}^{i x \cdot \xi-|\xi|^{2s}} \mathrm{d} \xi. 
$$%\end{eqnarray}
%Generally,   given a polynomial $L$ with real coefficients,   the heat kernel of the differential operator $L(D)$ could be expressed by
%$$  \cP_{_L}(t,x)=\int_{\mathbb{R}^N} \mathrm{e}^{i x \cdot \xi} L(i \xi)
%\mathrm{e}^{-t|\xi|^{\alpha}} \mathrm{d} \xi.$$
%Since $\mathrm{e}^{-t|\xi| \alpha^{2s}}$ is a tempered distribution 
Note that  $\cP_{s} \in \mathcal{C}^{\infty}\left((0, \infty)\times \mathbb{R}^N \right)$ (see e.g. \cite{BJE,DGV}). 
Under suitable conditions on $f$, a  solution of (\ref{eq 1.1 frac})  with initial data $u(0,\cdot) = f$ is expressed by
the formula
$$
u_{s,f}(t, x)=\int_{\mathbb{R}^N} \cP_{s}(t,x-y) f(y) \dd y.
$$
%%%%

It is well-known that the following two 
limits hold when $\s$ tends $0$ or $1$:
\begin{equation*}
\lim_{\s\to1^-}(-\Delta)^\s u(x)=-\Delta u(x)
\quad \text{and}\quad \lim_{\s\to0^+}  (-\Delta)^\s  u(x) = u(x)\quad\ \text{for $u\in \mathcal{C}^2_c(\R^N)$,}
\end{equation*}
see e.g. \cite{DPV}.  Furthermore  the following surprising expansion at $s=0$ is  proved in \cite{CT}
$$
(-\Delta)^s  u(x) = u(x) + s \lnlap u (x) + o(s) \quad \text{as }\,s\to 0^+\quad\text{ for all  }u \in \cC^2_c(\R^N) \text{ and }
x\in \R^N$$
 where the formal operator
$$%\begin{equation}\label{deriv}
\lnlap:= \frac{d}{d\s}\Big|_{\s=0} (-\Delta)^\s
$$%\end{equation}
 is the {\em logarithmic Laplacian}; more precisely (see e.g. \cite{CT}),
\begin{enumerate}
\item[(i)] for $1 < p \le \infty$, there holds $\lnlap  u \in L^p(\R^N)$ and $\frac{(-\Delta)^\s u- u}{\s} \to \lnlap  u$ in $L^p(\R^N)$ as $\s \to 0^+$;
\item[(ii)] $\mathcal{F}(\lnlap u)(\xi) = (2 \ln |\xi|)\,\widehat u (\xi)$ % = (\ln |\xi|^2)\, \widehat u (\xi)$
 \, for a.e. $\xi \in \R^N$.
\end{enumerate}
Recently, the following related topics on the logarithmic Laplacian has been investigated: the eigenvalues estimates \cite{CV1,LW}, log-Sobolev inequality \cite{FKT}, semilinear problems \cite{JSW,CDK}.
The domain of definition of $ \lnlap$ is the space  $Dom_{\CL_{\Gd}}(\BBR^N):=\mathscr{L}\left(\mathbb{R}^N\right)\cap C_D(\R^N)$, where 
$$%\begin{eqnarray}\label{eq-5}
\mathscr{L}\left(\mathbb{R}^N\right)=\left\{u: \mathbb{R}^N \rightarrow \mathbb{R}\ \, 
  {\rm measurable: }\ \, \int_{\mathbb{R}^N} \frac{|u(x)|}{1+|x|^{N}} \dd x<+\infty\right\}
$$
and $C_D(O)$ is the space of uniformly Dini continuous functions  $u$ in $O$ such that for any $x\in O$
$$\int_0^1\sup_{y\in O,\, |y-x|\leq r} |u(x)-u(y)|\frac{dr}r <+\infty .$$
If $u \in \mathscr{L} \left(\mathbb{R}^N\right)$ then
$\lnlap u$ can be defined as a distribution with the duality product $\left\langle \lnlap u,\varphi\right\rangle$
for every $\varphi\in\cS(\R^N)$ where $\cS(\R^N)$ denotes the Schwartz space of $C^\infty$ fast decaying functions.
%in \cite{CV,LW} estimates of the eigenvalues of the logarithmic Laplacian have been obtained  and some related topics c ould see \cite{JSW,FKT}. 
\smallskip

Applying the inverse Fourier transform to the Fourier expression of the diffusion equation associated to the kernel defined in (ii), we obtain that  the formal diffusion kernel of the logarithmic Laplacian has the following expression
 \begin{eqnarray}\label{eq-3-log}
\cP_{\ln}(t,x):=\int_{\mathbb{R}^N} |\xi|^{-2t}\mathrm{e}^{i x \cdot \xi} \mathrm{d} \xi .
 \end{eqnarray}
From now on, this formula defines the logarithmic diffusion kernel.

\begin{proposition}\label{pr 1.1}
Let the logarithmic diffusion kernel $\cP_{\ln}$ be defined by (\ref{eq-3-log}), then
\begin{equation}\label{hk 1}
\cP_{\ln}(t,x)=\cP_0(t)   |x|^{2t-N} \quad {\rm for} \ \ (t,x)\in \Big(0,\tfrac N2\Big)\times \left(\R^N\setminus\{0\}\right),
\end{equation}
 where 
\begin{equation}\label{e 2.1}
\cP_0(t)=\pi^{-\frac{N}{2}}\, 4^{-t} \frac{\Gamma(\frac{N-2t}{2})}{\Gamma(t)}.
\end{equation}
 Furthermore, we have that \smallskip

\noindent (i) \ $\cP_{\ln}(t,\cdot)\not\in L^1(\R^N) $, but $\cP_{\ln}(t,\cdot)\in L^{\frac N{N-2t},\infty}(\R^N) $
for every $t\in(0,\frac{N}{2})$, where  $L^{\frac N{N-2t},\infty}(\R^N)$ is the Marcinkiewicz space with exponent $\tfrac N{N-2t}$;\smallskip

\noindent (ii) $\cP_{\ln}(t,\cdot)$ blows up uniformly in any compact set of $\R^N$  as $t\to (\frac{N}{2})^-.$ More precisely there holds for $x\not=0$
$$ \displaystyle\lim_{t\to (\frac N2 )^-}   (N-2t)\cP_{\ln}(t,x) = 2^{2-N} \omega_{_{N}}.$$
\smallskip

\noindent (iii)  For any $x\not=0$,
$$\lim_{t\to 0^+}\cP_{\ln}(t,x)=0  $$ 
and $\displaystyle\lim_{t\to0^+}\cP_{\ln}(t,\cdot) =\delta_0$ in the sense of distributions, i.e.  
 $$\lim_{t\to0^+}\int_{\R^N}\cP_{\ln}(t,x)\varphi(x)\dd x=\varphi(0)\quad{\rm for}\ \ \varphi\in \cC_c(\R^N). $$
\end{proposition}
 
In the next table we emphasise the striking differences between $\cP_{\ln}$ and $\cP_s$ for $s\in (0,1)$, 
 \vspace{0.1cm}

\begin{center}
\renewcommand{\arraystretch}{1.6}
    \begin{tabular}{| l | l | l | l | l |} \hline
      & Lifespan & Time asymptotics  & Behaviour at $x=0$& Decay  as $|x|\to\infty$  \\         \hline
    $\cP_{s}$ & $(0,+\infty) $& decay rate $t^{-\frac{N}{2s}} $ at infinty & smooth& $ |x|^{-2s-N}$  
\\ \hline
$\cP_{\ln}$ & $(0,\, \frac{N}{2})$& blow-up rate  $\frac1{N-2t}$ & singular $|x|^{2t-N}$ & $ |x|^{2t-N}$ \\ \hline
\end{tabular}
\renewcommand{\arraystretch}{1}
\end{center}
  \vspace{0.1cm}
 Furthermore, $\cP_{\ln}(t,\cdot)$ is the fundamental solution of the fractional Laplacian $(-\Delta)^t$ for $t\in (0,1)$.
 In the sequel we call the expression below the {\it  logarithmic diffusion equation}
   \begin{equation}\label{L1}
\partial_tu+ \lnlap u=0,
 \end{equation}
without presuming in what sense it holds.

 \begin{definition}\label{strong sol}
 
We say that a function $u$ defined in  $\cC\Big([0, T) \times  \mathbb{R}^N \Big)$ is a strong solution of the logarithmic diffusion equation (\ref{L1}) with initial data $f\in C( \mathbb{R}^N)$, if $u(t,.)\in Dom_{\CL_\Gd}(\BBR^N)$ for all $t>0$
and the following conditions hold: \smallskip

\nind (i) $\partial_tu \in \cC\Big((0, T) \times \mathbb{R}^N \Big)$;\smallskip 

\nind(ii) $u \in \cC\Big([0, T) \times  \mathbb{R}^N \Big)$;  \smallskip

\nind(iii) the equation is satisfied pointwisely,  i.e. for every $(t,x) \in (0, T)\times \mathbb{R}^N$  
\begin{equation}\label{e 1-v}
\partial_t u(t,x)+\lnlap u(t,x) =0,
\end{equation}
 and $u(0,x)=f(x)$  for every $x \in \mathbb{R}^N$.
 \end{definition}

As a consequence of the equation, in this framework the function $(t,x)\mapsto \lnlap u(t,x)$ is continuous in $(0, T) \times \mathbb{R}^N$. \medskip

We  prove the following analogue of Widder's representation for the logarithmic Cauchy problem.

\begin{theorem}\label{thm1}
Let  $0<T\leq \frac{N}{2}$ and $f$ be a nonnegative function in  $\cC(\R^N) \cap L^1(\R^N,(1+|x|)^{2T-N} dx)$. 
Then   problem (\ref{eq 1.1})  admits a unique strong solution $u=u_f$, which  is positive and has the form 
\begin{equation}\label{e 1-0}
u_f(t,x)=\int_{\mathbb{R}^N} \cP_{\ln}(t,x-y) f(y) \dd y=\CP_0(t) I_{2t}\ast f(x),\quad\forall\, (t,x)\in [0,T) \times \R^N,
\end{equation}
where $I_{2t}(x)=|x|^{2t-N}$ is the Riesz kernel of order $2t$. Moreover
\begin{equation}\label{e1-1}
\lim_{t\to(\frac{N}{2})^-} (N-2t) u_f(t,x)=\frac{2}{(4\gp)^{\frac N2}} \|f\|_{L^1}\quad  {\rm for\ any\ }\, x\in \R^N. 
\end{equation}
Furthermore, the following implications hold:\smallskip

\nind (i) If $f$ satisfies for some $M>0$ and $\tau<-N$   
 $$ 0\lvertneqq  f(x)\leq M(1+|x|)^{\tau}\quad {\it for}\ \  x\in\R^N,$$
 then there exists $c=c(N,f)>1$ such that for any $(t,x)\in (0,T)\times\R^N$
   \begin{equation}\label{e 1.0-0}
\frac{\cP_0(t)}{ct}    (1+|x|)^{2t-N}\leq  u_f(t,x)   \leq c \frac{\cP_0(t)}{t}   (1+|x|)^{2t-N},
 \end{equation}
 where $\cP_0$ is defined in (\ref{e 2.1}). 
 
\nind (ii) If $f$ satisfies for some $M>1$ and $\tau\in(-N, -2T]$
 $$\frac1M(1+|x|)^{\tau}\leq f(x)\leq M(1+|x|)^{\tau}\quad {\it for}\ \  x\in\R^N,$$
 then  there exists $c=c(N,M,\gt)>1$   such that for any $(t,x)\in (0,T)\times\R^N$
   \begin{equation}\label{e 1.0-1}
\frac1c \frac{\cP_0(t)}{2t}  (1+|x|)^{ 2t+\tau  } \leq  u_f(t,x)\leq c \frac{\cP_0(t)}{2t}  (1+|x|)^{  2t+ \tau}.
 \end{equation}
 
\nind (iii)  If it is assumed that, for some $M>1$, 
 $$\frac1M(1+|x|)^{-N}\leq f(x)\leq M(1+|x|)^{-N}\quad {\it for}\ \  x\in\R^N,$$
  then  there exists $c=n(N,M)>1$    such that for any $(t,x)\in (0,T)\times\R^N$
   \begin{equation}\label{e 1.0-2}
\frac1c \frac{\cP_0(t)}{2t}  (1+|x|)^{ 2t-N  } \ln(e+|x|)\leq  u_f(t,x)\leq c \frac{\cP_0(t)}{2t}  (1+|x|)^{ 2t-N  } \ln(e+|x|).
 \end{equation}
\end{theorem}

Involving a signed initial data, we have the same Widder's   representation by the linearity of the operators.  

 %The counter part of this result is the following characterization of the initial trace of any nonnegative solution of (\ref{e 1-v}).
 \begin{corollary}\label{thmTr}
Let  $0<T\leq \frac{N}{2}$ and $f\in L^1(\R^N,(1+|x|)^{2T-N} dx)\cap \cC(\R^N)$.
Then   problem (\ref{eq 1.1})  admits a unique strong solution $u=u_f$ expressed by (\ref{e 1-0}). Moreover,  (\ref{e1-1}) holds. 
%\begin{equation}\label{e 1-0}
%u_f(t,x)=\int_{\mathbb{R}^N} \cP_{\ln}(t,x-y) f(y) \dd y=\CP_0(t) I_{2t}\ast f(x),\quad\forall\, (t,x)\in [0,T) \times \R^N.
%\end{equation}
%where $I_{2t}(x)=|x|^{2t-N}$ is the Riesz kernel of order $2t$. 
%\begin{equation}\label{e1-1}
%\lim_{t\to(\frac{N}{2})^-} (N-2t)u_f(t,x)=\frac{2}{(4\gp)^{\frac N2}\Gamma(\frac N2)}\int_{\R^N}f(y)dy\quad {\it for\ any}\ \, x\in \R^N. 
%\end{equation}
\end{corollary}

Note that the bounds of the strong solution  in Theorem \ref{thm1} show that the function $u_f(t,\cdot)$ does not belong to $L^1(\R^N)$. Hence we propose a definition of weak solutions to (\ref{L1}) valid for more general initial data $f$. 

 \begin{definition}\label{weak sol}
  We say that the function $u$ is a weak solution of the logarithmic diffusion equation (\ref{e 1-v}) with initial data $f\in L_{  {\rm loc }}^{1} (\mathbb{R}^N )$
 if the following conditions hold:\smallskip

\nind (i)  $u \in L^{1}\left(\left(0, T^{\prime}\right], \, L^1 \big(\mathbb{R}^N, (1+|x|)^{2(T-T')-N}dx\big)\right)$
for every $T^{\prime}<T$;  

\nind (ii)  $u \in \cC\left((0, T),\,  L_{  {\rm loc }}^{1} (\mathbb{R}^N )\right)$;\smallskip

\nind(iii)  For $\varphi \in \bbX_0$ and any $T^{\prime}\in(0,T)$ there holds
\begin{equation}\label{def w1}
 \int_{0}^{T^{\prime}} \int_{\mathbb{R}^N}\left[- \partial_t \varphi  +  \lnlap
\varphi \right] u  \dd t\dd x  =\int_{\mathbb{R}^N} \varphi(0, \cdot) f\dd x-\int_{\mathbb{R}^N} u\left(T^{\prime}, x \right)
\varphi\left(T^{\prime},x \right) \dd x,
\end{equation}
where
$\bbX_0=\left\{\gz\in \cC^1([0,T)\times\mathbb{R}^N)\!:\gz(t,\cdot)\in \cC^1_c(\mathbb{R}^N) {\rm\ for\ all\ }t\in [0,T)\right\}$.
 \end{definition}

When $f=\mu$ is a  Radon measure, identity (\ref{def w1}) is replaced by
  \begin{equation}\label{def w1-d_0}
\int_{0}^{T^{\prime}} \int_{\mathbb{R}^N}\left[- \partial_t \varphi  +  \lnlap
\varphi \right] u  \dd t\dd x  =\int_{\mathbb{R}^N} \varphi(0, \cdot)d\gm-\int_{\mathbb{R}^N} u\left(T^{\prime}, x \right)
\varphi\left(T^{\prime},x \right) \dd x.
\end{equation}

Our first existence and uniqueness result which shows the role of the Widder's type representation is as follows.
\begin{theorem}\label{th 1.2}
Assume  $0<T\leq \frac{N}{2}$ and $  f\in  L^1(\R^N,(1+|x|)^{2T-N} dx)$, 
 then problem (\ref{eq 1.1}) admits a unique weak solution   $u_f$ expressed by
 \begin{eqnarray}\label{express 0}
u_f(t,x)=\int_{\mathbb{R}^N} \cP_{\ln}(t,x-y) f(y) \dd y\quad{\rm for \ all}\ t\in (0,T),\ {\rm a.e. }\;\ x\in\R^N.
\end{eqnarray}
\end{theorem}

The result is also valid if $f$ is replaced by a Radon measure $\gm$ in $\BBR^N$ satisfying 

\bel{R1}
\int_{\BBR^N}(1+|x|)^{2T-N} d|\gm|(x)<\infty.
\ee
Of particular importance is the case where $\gm$ is a Dirac mass. 
\begin{corollary}\label{th 1.3}
Let  $0<T\leq \frac{N}{2}$, then $\cP_{\ln}$ is the unique weak solution of  
\begin{equation}\label{eq 1.1-d_0}
\left\{ \arraycolsep=1pt
\begin{array}{lll}
\displaystyle \partial_tu+ \lnlap u=0\quad \  &{\rm in}\ \   (0,T)\times \R^N,\\[2.5mm]
 \phantom{ \lnlap \ \,   }
\displaystyle   u(0,\cdot)=\delta_0\quad \ &{\rm in}\  \   \R^N.
\end{array}
\right.
\end{equation}
\end{corollary}

\nind Notice that,

\nind (a)  from the uniqueness, the classical solution obtained in Theorem \ref{thm1} is 
a weak solution of  (\ref{eq 1.1});  \smallskip

\nind (b)   for a  nonnegative function $f\in  L^1(\R^N)$,  the weak solution  $u_{s,f}$ of the fractional diffusion equation of (\ref{eq 1.1 frac}) exists for   
 $(0, +\infty)\times \R^N$, while the solution  of the nonlocal diffusion equation (\ref{eq 1.1})
 exists only in  $(0, \frac N2)\times \R^N$, and it blows up locally in $\R^N$
as $t\to (\frac N2)^-$ in accordance with the expression (\ref{e1-1}). \medskip

%\subsection{Fundamental solutions}

The final aim of this study is to use the diffusion kernel $\cP_{\ln}$ to give an expression of the {\it fundamental solution} of the logarithmic Laplacian,  that is the solution of
\begin{equation}\label{eq 8.1}
\lnlap u=k\delta_0\quad\   {\rm in}\ \  \cD'(\R^N).
 \end{equation}
 
\begin{theorem}\label{th 8.1}
 Let $N\geq 3$, then the operator $\lnlap$ has a radially symmetric fundamental  solution $\Phi_{\ln} $ in $\R^N$ expressed by
 \bel{fund}
 \Phi_{\ln}(x)=\int_0^1\cP_0(t)|x|^{2t-N}dt+\int_0^1\cP_0(t)\int_{\BBR^N}|x-y|^{2t-N}\Phi_1(y)dydt\quad {\rm for\ }\, x\in \BBR^N\setminus\{0\},
 \ee
 where 
  \bel{fund1}
\cP_0(t)=\gp^{-\frac N4}\frac{\Gamma (\frac{N-2t}{2})}{4^t\Gamma (t)}
 \ee
 and $\Phi_1$ is the fundamental of the Helmoltz equation $-\Gd u-u=\delta_0$ in $\cD'(\BBR^N)$ expressed via Hankel functions.

Furthermore there exist two positive constants $c=c_1(N)$ and $c_1=c_2(N)$ such that  
\begin{equation}\label{7xy}
\left|\Phi_{\ln} (x)-c_0\frac{ |x|^{-N}}{\ln^2 |x|}\right| \leq c_1\frac{ |x|^{-N}}{\ln^3 |x|}\qquad  {\rm for\ all\ }\ x\in B_{\frac1e}(0)\setminus\{0\}
\end{equation}
 and
\begin{equation}\label{best inf}
\left|\Phi_{\ln}(x)\right|\leq c_2\frac{|x|^{\frac {3-N}2}}{\ln|x|}\qquad {\rm for\ }\,|x|\geq 2,
\end{equation}
where 
$$c_0=\tfrac14\pi^{-\frac{N}{2}}\Gamma\left(\tfrac N2\right). $$

 \end{theorem}
 
 Our idea is to apply  the parabolic operator $\partial_t+\lnlap$
 to connect the fundamental solutions of the logarithmic Laplacian with   the one of the Helmholtz operator $-\Gd-I$. This requires the lifespan of $\cP_{\ln}$ to contain  the interval $(0,1]$. 
When $N\geq3$, it is well defined since $\cP_0(t)$ for $t=1<\frac N2$, but it fails  when $N=2$,
since in that case $\cP_0(1)=+\infty$. {\it However we conjecture that the fundamental solution $\Phi_{\ln}$ exists for $N=1,2$ and it decays 
like $|x|^{\frac{1-N}{2}}$ as $|x|\to+\infty$. }

In a bounded Lipschitz domain $\Omega$ containing $0$, we consider the problem below, the solution of which is the fundamental solution in 
$\Gw$,
\begin{equation}\label{eq 7.2}
\left\{ \arraycolsep=1pt
\begin{array}{lll}
  \lnlap u=\delta_0 \quad \  &{\rm in}\ \   \cD'(\Omega),\\[2mm]
 \phantom{ \lnlap     }
\displaystyle   u  = 0 \quad \ &{\rm{in}}\  \    \R^N\setminus\Omega,
\end{array}
\right.
\end{equation} 
where $\cD'(\Omega)$ denotes the space of distributions in $\Gw$ that is the topological dual of the space $\cD(\Omega)$ of $C^\infty$ functions with compact support endowed with the 
inductive limit topology (see \cite{Sch}).  We denote by $\sigma(\lnlap,\Gw)$ the spectrum of $\lnlap$ in $\Omega$ with zero  condition in $\Gw^c$, that is the set of $\gl\in\BBR$ such that
there exists a nontrivial solution to 
\begin{equation}\label{eq 7.1} 
\left\{ \arraycolsep=1pt
\begin{array}{lll}
  \lnlap u=\lambda  u  \quad \  &{\rm in}\ \   \Omega,\\[2mm]
 \phantom{ \lnlap     }
\displaystyle   u  = 0 \quad \ &{\rm{in}}\  \    \R^N\setminus\Omega.
\end{array}
\right.
\end{equation}
It is proved in \cite{CT} that $\sigma(\lnlap,\Gw)=\{\lambda_k^{^L}(\Omega)\}_{k\geq 1}$ where  $k\mapsto \lambda_k^{^L}(\Omega)$  is
  increasing and  
 $$\lim_{k\to+\infty}\lambda_k^{^L}(\Omega)\to+\infty. $$
Furthermore $ker{(  \lnlap+\lambda_1^{^L}(\Omega))}$ is 1-dimensional and generated by a positive function $\varphi_1$.
First we characterise the fundamental solution of $\lnlap$ in $\Omega$.

 \begin{theorem}\label{th 7.1}
 Let $N\geq 1$,  $\Omega$ be a bounded Lipschitz domain satisfying a uniform exterior sphere condition,  containing $B_{1}$ and such that $0\notin \sigma(\lnlap,\Gw)$. 
 Then
problem (\ref{eq 7.2}) possesses a  fundamental solution
$\Phi^\Omega_{\ln}$ and  there exists $c_3 >0$ such that  
\bel{E1}\left|\Phi^\Omega_{\ln}(x)-c_0\frac{|x|^{-N}}{\ln^2 |x|}\right|  \leq c_3\frac{ |x|^{-N}}{\ln^3 |x|}\quad {\it for}\
x\in B_{\frac1e}(0)\setminus\{0\}.\ee
Furthermore, there holds near $\partial\Gw$
\begin{equation}\label{7xx}
|\Phi^\Omega_{\ln}(x)|=O\left(\frac{1}{|\ln\gr(x)|^\gt}\right)\quad\text{for all }\gt\in \left(0,\frac 12\right),
\end{equation}
where $\gr(x)=\dist(x,\partial\Gw)$.
 \end{theorem}
 
We also gives some estimates of $\Phi_{\ln}$ is some Orlicz classes and we prove, at least if $N\geq 4$, that under mild assumptions on  $f$ the function $u=\Phi_{\ln}\ast f$ is Dini continuous. It is therefore a strong solution of 
\begin{equation}\label{AAA}
\CL_\Gd u=f\qquad\text{in }\, \BBR^N.
\end{equation}
 \begin{theorem}\label{cr 5.1} Let $N\geq 4$, $\gth>1$ and $f\in L^1(\BBR^N)\cap L^{\infty,\gth}(\BBR^N)$, where
 \begin{equation}\label{AAB}
L^{\infty,\gth}(\BBR^N):=\left\{h\in L^\infty_{loc}(\BBR^N)\!:\, \esssup_{y\in\BBR^N}(1+|y|)^\gth |f(y)|<\infty\right\}. 
\end{equation}
Then for any $R>0$ there exists $c_4=c_4(N,R,\gth)$ such that 
 \begin{equation}\label{AAC}
|\Phi_{\ln}\ast f(x)-\Phi_{\ln}\ast f(x')|\leq \frac{c_4}{1+\ln^2(|x-x'|)}\left(\norm f_{L^1}+\norm f_{L^{\infty,\gth}}\right)\quad  \text{ for  all }x,x'\in B_R,
\end{equation}
where
 \begin{equation}\label{AAD}
\norm f_{L^{\infty,\gth}}=\esssup_{y\in\BBR^N}(1+|y|)^\gth |f(y)|.
\end{equation}
\end{theorem}

\medskip
 
In the previous theorem the assumption $N\geq 4$ is made necessary by our method based upon the key estimates obtained in Theorem \ref {P3} (see  the remark after its proof). {\it We  conjecture that the above result is still valid when $N\geq 3$, up to some modifications of the estimates. }

 \smallskip

The paper is organized as follows. Section 2 is devoted to the study of the
diffusion kernel $\cP_{\ln}$ of $\partial_t+\lnlap$ and to the existence and the representation 
of classical solutions of (\ref{eq 1.1}).   Section 3 is devoted to the Cauchy problem associated to the logarithmic Laplacian, first in  the weak sense and then in the classical sense, based upon the proof of the uniqueness of weak solutions.  In section 4,  we construct the fundamental solution of the  logarithmic Laplacian using the diffusion kernel,  representation formula for solutions of the Cauchy problem and precise estimates of the kernel of Helmholtz operator in $\BBR^N$. In Section 5 we use the estimates of 
$\Phi_{\ln}$ to prove some imbedding properties in Orlicz spaces,  to obtain point-wise estimate of $\nabla \Phi_{ln}$ in $\BBR^N\setminus\{0\}$ and to prove that the solution of the Dirichlet problem in $\BBR^N$ is Dini continuous.

 \section{Basic properties}
 In the course of this paper $c$ will denote a generic positive constant, depending on different parameters (which can be specified in some cases), but not on the variable, the value of which can vary from one occurence to another. In some cases, in order to avoid misunderstanding, we introduce constants $c'$ or $c''$. 
  \subsection{The logarithmic diffusion kernel}

  \noindent{\it Proof of Proposition \ref{pr 1.1}. }
 We recall (see e.g. \cite[Chapter V, p. 117]{S}) that
$$
 \mathcal{F}(|\cdot|^\tau) =
\sigma(\tau)|\cdot|^{ -N-\tau}  \quad {\rm in}\ \ \mathcal{S}'(\R^N)
$$
with
$$\sigma(\tau):=2^{\tau+N}\pi^{\tfrac{N}{2}}\frac{\Gamma(\frac{\tau+N}2)}{\Gamma(-\frac\tau2)}.$$ Consequently, if $\tau \not = -N$, we have that
$$
  \mathcal{F}^{-1}\bigl(|\cdot|^{-N-\tau}\bigr) =  \frac{1}{\sigma(\tau)} |\cdot|^{\tau}   \quad {\rm in}\ \ \mathcal{S}'(\R^N)
$$
and
$$
 \mathcal{F}^{-1}(1) =  \delta_0 \quad {\rm in}\ \ \mathcal{S}'(\R^N).
$$
From the properties of the Fourier transform there holds
\bel{X1}\partial_t\widehat{\cP_{\ln}}(t,\xi)+2\ln|\xi| \,\widehat{\cP_{\ln}}(t,\xi)=0,\ee  
 which implies
 $$\widehat{\cP_{\ln}}(t,\xi)=e^{-2t\ln|\xi|  }=|\xi|^{-2t} \quad{\rm for}\ \ (t,\xi)\in(0,\tfrac N2)\ti\R^N\setminus\{0\}.$$
By the inverse Fourier transform, 
 $$\cP_{\ln}(t,x)=  \cP_0(t) |x|^{2t-N} \quad {\rm for} \ \ (t,x)\in (0,\tfrac N2)\times \R^N, $$
 where, we recall it, 
  $$%\begin{equation}\label{e 2.1}
 \cP_0(t)=\pi^{-\frac{N}{2}}\, 4^{-t} \frac{\Gamma(\frac{N-2t}{2})}{\Gamma(t)}.
 $$%\end{equation}
Clearly $\cP_0$ is a positive and smooth function in $(0,\frac{N}{2})$ and it vanishes at $t=0$. 
The expression of $\cP_{\ln}$ shows that
 $\cP_{\ln}(t,\cdot)$ does not belong to  $L^1(\R^N)$ for any $t\in(0,\frac N2)$, this is part (i). 
 
 Since $x\Gg(x)\to 1$ as $x\downarrow 0$ 
 we obtain 
$$   \cP_0(t)\cdot \tfrac{N-2t}2\to (4\pi)^{-\tfrac{N}{2}}\  \Gamma^{-1}(\tfrac{N}{2})  \quad{\rm as}\ \  t\to (\tfrac N2)^-. $$
This implies that  for any $x\not=0$
$$ \cP_{\ln}(t,x)\cdot \tfrac{N-2t}2 \to (4\pi)^{-\tfrac{N}{2}}\  \Gamma^{-1}(\tfrac{N}{2}) \quad{\rm as}\ \  t\to (\tfrac N2)^-.$$
Now $(ii)$ follows by the fact that 
$$\lim_{t\to (\tfrac N2)^-}\Gamma(\tfrac{N-2t}{2})(\tfrac{N}2- t)=1\quad{\rm and}\quad \lim_{t\to (\tfrac N2)^-} |x|^{2t-N}=1\quad {\rm for}\ x\not=0. $$
Because $\dsps\lim_{\tau\to 0}\Gamma(\tau)=\infty$, we have for any $x\not=0$,
$$\lim_{t\to 0^+}\cP_{\ln}(t,x)=0.  $$

 For the proof of part $(iii)$, if $\varphi\in \cC_c (\R^N)$ has support in $B_{\gs_0}$, for any $\ge>0$ there exists $\gd>0$ such that 
 $$|\varphi(x)-\varphi(0)|\leq \ge\quad \text{ if }\ \,|x|\leq\gd.$$
 Then 
$$\BA{lll}\dsps
 \int_{\R^N}\cP_{\ln}(t,x)\varphi(x) \dd x =  \pi^{-\tfrac{N}{2}} 4^{-t} \frac{\Gamma(\frac{N-2t}{2})}{\Gamma(1-t)} t\left(\int_{B_{\gs_0}}  |x|^{2t-N} \dd x\right)\varphi(0)
 \\[5mm]\phantom{\dsps
 \int_{\R^N}\cP_{\ln}(t,x)\left(\varphi(x)-\varphi(0)\right) \dd x}
 \dsps+ \pi^{-\tfrac{N}{2}} 4^{-t} \frac{\Gamma(\frac{N-2t}{2})}{\Gamma(1-t)} t\int_{B_{\gs_0}}  |x|^{2t-N}\left(\varphi(x)-\varphi(0)\right) \dd x.
 \EA$$
We have that
$$\pi^{-\tfrac{N}{2}} 4^{-t} \frac{\Gamma(\frac{N-2t}{2})}{\Gamma(1-t)}\to \pi^{-\frac{N}{2}}\Gamma(\tfrac N2) \quad{\rm as}\ \ t\to0^+,$$ 
$$t\int_{B_{\gs_0}}  |x|^{2t-N} \dd x=\frac{\gs_0^{2t}\gw_{_{N}}}{2}\to \frac{\gw_{_{N}}}{2}\quad{\rm as}\ \ t\to0^+
$$
 and
$$\BA{lll}
\dsps  t\left|\int_{B_{\gs_0}}  |x|^{2t-N}\left(\varphi(x)-\varphi(0)\right) \dd x\right|\leq t\int_{B_{\gd}}  |x|^{2t-N}\left|\varphi(x)-\varphi(0)\right| \dd x\\[4mm]\phantom{-------------------}\dsps
+t\int_{B_{\gs_0}\setminus B_{\gd}}  |x|^{2t-N}\left|\varphi(x)-\varphi(0)\right| \dd x\\[0mm]
\phantom{\dsps  t\left|\int_{B_{\gs_0}}  |x|^{2t-N}\left(\varphi(x)-\varphi(0)\right) \dd x\right|}
\leq I+II.
\EA $$
Then
 $$ 
 I\leq \ge\gw_{_{N}}\gd^{2t} \to \ge\gw_{_{N}}\quad\text{as }t\to 0$$
 and
 $$II\leq \norm\varphi_{L^\infty}\gw_{_{N}}(\gs_0^{2t}-\gd^{2t})\to 0\quad\text{as }t\to 0.
 $$
Since $\ge$ is arbitrary and $$\pi^{-\tfrac{N}{2}}\Gamma(\tfrac N2)   \omega_{_{N}}=2, $$
we obtain 
 $$\lim_{t\to0^+}\int_{\R^N}\cP_{\ln}(t,x)\varphi(x)\dd x=\varphi(0). $$
\hfill$\Box$
%%%%%%%%%%%%%%%%%%%%%%%%%%%%%%%%%%%%%%%%%%%%%%%%%%%%%%%%%%%%%%%%%%%%%%%%%%%%%%%%%%%%%%%%%%%%%%%%%%%%%%%%%%%%%%%%%%%%%%%%%%%%%%%%%%%%
 
 \begin{lemma}\label{lm a-1} Assume $\gt<0$ and $(t,x)\in (0,-\frac\gt 2)\ti\BBR^N$, then\smallskip
 
\nind (i) If $-N<\tau< 0$ there exists $c=c(N,\gt)>1$ such that 
$$\frac{1}{ct} \min\left\{|x| ^{\tau+2t},t+|x|^{2t}\right\}\leq  \int_{\R^N} |x-y|^{2t-N} (1+|y|)^{\tau} \dd y\leq \frac{c}{t} |x| ^{2t+\gt}.$$

\nind (ii) If $\tau=-N$  there exists $c=c(N,\gt)>1$ such that 
$$\frac{1}{ct} \min\left\{|x| ^{2t-N}\ln (e+|x|),t+|x|^{2t}\right\}\leq  \int_{\R^N} |x-y|^{2t-N} (1+|y|)^{-N} \dd y\leq \frac{c}{ t} |x| ^{2t-N}\ln (e+|x|).$$

\nind (iii) If $\tau<-N$ there exists $c=c(N,\gt)>1$ such that 
$$\frac{1}{ct} \min\left\{|x| ^{2t-N},t+|x|^{2t}\right\}\leq  \int_{\R^N} |x-y|^{2t-N} (1+|y|)^{\tau} \dd y\leq \frac{c}{ t}\max\left\{|x| ^{2t+\gt},|x|^{2t-N}\right\}.$$
 \end{lemma}
 %%%PROOF%%%%%%%%%%%%%%%%%%%%%%%%%%%%%%%%%%%%
\proof  Assume $x\neq 0$. The integral over $\BBR^N$ is the sum of the integrals over $B_{\frac{|x|}2}(0)$, over $B_{\frac{|x|}2}(x)$ and over $\BBR^N\setminus\left(B_{\frac{|x|}2}(0)\cup B_{\frac{|x|}2}(x)\right)$.
Clearly
\bel{Ex1}\BA{lll}\dsps
 \int_{B_{\frac{|x|}2}(x)} |x-y|^{2t-N} (1+|y|)^{\tau} \dd y    \leq     \frac{(2+ |x|)^{\tau}}{2^{\tau}} \int_{B_{\frac{|x|}2}(x)}|x-y|^{2t-N} \dd y
=  \frac{\gw_{_{N}}}{ 2^{2t+\gt+1}}\frac{|x|^{2t}(2+|x|)^{\tau}}{t}
\EA
\ee
and
\bel{Ex2}\BA{lll}\dsps
 \int_{B_{\frac{|x|}2}(x)} |x-y|^{2t-N} (1+|y|)^{\tau} \dd y    \geq  \frac{(2+ 3|x|)^{\tau}}{2^{\tau}}\int_{B_{\frac{|x|}2}(x)}|x-y|^{2t-N} \dd y
\geq  \frac{\gw_{_{N}}}{ 2^{2t+\gt+1}}\frac{|x|^{2t}(2+3|x|)^{\gt}}{t}.
\EA
\ee
Concerning the integral over $B_{\frac{|x|}2}(0)$, one has
$$\BA{lll}\dsps\frac{|3x|^{2t-N}\gw_{_{N}}}{2^{2t-N}}\int_0^{\frac{|x|}2}(1+r)^\gt r^{N-1}dr \leq \int_{B_{\frac{|x|}2}(0)} |x-y|^{2t-N} (1+|y|)^{\tau} \dd y\\[4mm]\phantom{----------------------}\dsps\leq \frac{|x|^{2t-N}\gw_{_{N}}}{2^{2t-N}}\int_0^{\frac{|x|}2}(1+r)^\gt r^{N-1}dr.
\EA$$
If $|x|\to\infty$, there holds 
$$\int_0^{\frac{|x|}2}(1+r)^\gt r^{N-1}dr=\left\{\BA{lll}\dsps \frac{|x|^{N+\gt}}{2^{N+\gt}(N+\gt)}(1+o(1))\quad&\text{if }\,&-N<\gt<0,\\[3mm]
\ln|x|(1+o(1))&\text{if }\;&\gt=-N,\\[2mm]
c^*(1+o(1))&\text{if }\;&\gt<-N
\EA\right.$$
for some explicit value $c^*$, and if $x\to 0$, 
$$\int_0^{\frac{|x|}2}(1+r)^\gt r^{N-1}dr=\frac{|x|^N}{2^NN}(1+o(1)).
$$
This implies
\bel{Ex3}\int_{B_{\frac{|x|}2}(0)} |x-y|^{2t-N} (1+|y|)^{\tau} \dd y\leq c\left\{\BA{lll}\inf\{|x|^{2t+\gt},|x|^{2t}\}\quad&\text{if }\,&-N<\gt<0,\\[2mm]
\inf\{|x|^{2t-N}\ln(|x|+e) |,|x|^{2t}\}\quad&\text{if }\,&\gt=-N,\\[2mm]
\inf\{|x|^{2t-N} |,|x|^{2t}\}\quad&\text{if }\,&\gt<-N
\EA\right.
\ee
and
\bel{Ex4}\int_{B_{\frac{|x|}2}(0)} |x-y|^{2t-N} (1+|y|)^{\tau} \dd y\geq \frac 1c\left\{\BA{lll}\inf\{|x|^{2t+\gt},|x|^{2t}\}\quad&\text{if }\,-N<\gt<0,\\[2mm]
\inf\{|x|^{2t-N}\ln(|x|+e) ,|x|^{2t}\}\quad&\text{if }\,-N=\gt,\\[2mm]
\inf\{|x|^{2t-N} ,|x|^{2t}\}\quad&\text{if }\,-N>\gt,
\EA\right.
\ee
where the constant $c$ depends only on $N$. \smallskip

  Finally, we consider the integral over $\R^N\setminus \Big( B_{\frac{|x|}2}(0)\cup B_{\frac{|x|}2}(x)\Big)$,
$$\BA{lll}\dsps
 \int_{\R^N\setminus \big( B_{\frac{|x|}2}(0)\cup B_{\frac{|x|}2}(x)\big)} |x-y|^{2t-N} (1+|y|)^{\tau} \dd y   \\[4mm]
 \phantom{--------------}\dsps=|x|^{2t+\gt}
  \int_{\R^N\setminus \big( B_{\frac{1}2}(0)\cup B_{\frac{1}2}({\bf e}_x)\big)} |{\bf e}_x-z|^{2t-N} (|x|^{-1}+|z|)^{\tau} \dd z,
\EA
$$
 where ${\bf e}_x=\frac x{|x|}$. This integral is invariant by rotation, thus we can assume that ${\bf e}_N=(0,0,...,1)\in\R^N$. The function
 $$|x|\mapsto  H(|x|):=\int_{\R^N\setminus \big( B_{\frac{1}2}(0)\cup B_{\frac{1}2}({\bf e}_N)\big)} |{\bf e}_N-z|^{2t-N} (|x|^{-1}+|z|)^{\tau} \dd z
 $$
is increasing, therefore
$$\dsps 0=\lim_{|y|\to 0}H(|y|)\leq H(|x|)\leq \int_{\R^N\setminus \big( B_{\frac{1}2}(0)\cup B_{\frac{1}2}({\bf e}_N)\big)} |{\bf e}_N-z|^{2t-N} |z|^{\tau} \dd z:=H(\infty)<\infty.
$$
This yields
\bel{Ex5}
\int_{\R^N\setminus \big( B_{\frac{|x|}2}(0)\cup B_{\frac{|x|}2}(x)\big)} |x-y|^{2t-N} (1+|y|)^{\tau} \dd y\leq |x|^{2t+\gt}H(\infty).
\ee
For the lower bound, $B_{\frac{1}2}(0)\cup B_{\frac{1}2}({\bf e}_N)\subset B_2(0)$, therefore 
$$\BA{lll}\dsps
\int_{\R^N\setminus \big( B_{\frac{1}2}(0)\cup B_{\frac{1}2}({\bf e}_N)\big)} |{\bf e}_N-z|^{2t-N} (|x|^{-1}+|z|)^{\tau} \dd z
 \geq  \int_{\R^N\setminus B_2(0)}  |{\bf e}_N-z|^{2t-N} (|x|^{-1}+|z|)^{\tau} \dd z\\[4mm]
 \phantom{\dsps
\int_{\R^N\setminus \Big( B_{\frac{1}2}(0)\cup B_{\frac{1}2}({\bf e}_N)\Big)} |{\bf e}_N-z|^{2t-N} (|x|^{-1}+|z|)^{\tau} \dd z}\dsps
\geq 2^{2t-N}\int_{\R^N\setminus B_2(0)} |z|^{2t-N}(|x|^{-1}+|z|)^{\tau} \dd z\\[4mm]
 \phantom{\dsps
\dsps
\int_{\R^N\setminus \Big( B_{\frac{1}2}(0)\cup B_{\frac{1}2}({\bf e}_N)\Big)} |{\bf e}_N-z|^{2t-N} (|x|^{-1}+|z|)^{\tau} \dd z}\dsps
\geq2^{2t-N}\gw_{_{N}} \int_2^\infty r^{2t-1} (|x|^{-1}+r)^{\tau} \dd r.
\EA$$
But the following relation holds
$$\int_2^\infty r^{2t-1} (|x|^{-1}+r)^{\tau} \dd r=|x|^{-2t-\gt}\int_{2|x|}^\infty s^{2t-1} (s+1)^{\tau} \dd s\geq c\min\left\{1,|x|^{-\gt-2t}\right\}.
$$
This implies
\bel{Ex6}
\int_{\R^N\setminus \big( B_{\frac{|x|}2}(0)\cup B_{\frac{|x|}2}(x)\big)} |x-y|^{2t-N} (1+|y|)^{\tau} \dd y\geq c\min\left\{|x|^{2t+\gt},1\right\}.
\ee
Since $0<t<-\frac\gt 2$ the proof follows by combining (\ref{Ex1})--(\ref{Ex6}).
\hfill$\Box$\medskip

\nind \begin{remark}\label{X>1}When $|x|\geq 1$, the estimates in Lemma \ref{lm a-1} with possibly a new constant $c=c(N,\gt)>1$ can be expressed under the following form:
\bel{Ex7}\BA{lllll}
 {\rm (i)} 
\quad&\dsps\frac{1}{ct}|x| ^{2t+\gt}\leq  \int_{\R^N} |x-y|^{2t-N} (1+|y|)^{\tau} \dd y\leq \frac{c}{t} |x| ^{2t+\gt}\quad\qquad\text{if } \tau\in(-N, 0);\\[4mm]
 {\rm (ii)} \dsps
&\dsps\frac{1}{ct} |x| ^{2t-N}\ln (e+|x|)\leq  \int_{\R^N} |x-y|^{2t-N} (1+|y|)^{-N} \dd y
%\\[4mm]
\dsps\leq \frac{c}{ t} |x| ^{2t-N}\ln (e+|x|)
\quad&\text{if }\tau= -N;\qquad\quad  \\[4mm]
  {\rm (iii)}\dsps
&\dsps\frac{1}{ct} |x| ^{2t-N}\leq  \int_{\R^N} |x-y|^{2t-N} (1+|y|)^{\tau} \dd y\leq \frac{c}{ t}|x| ^{2t-N}\qquad \text{if }\tau< -N.
\EA
\ee
\end{remark}
%%%%%%%%%%%%%%%%%%%%%%%%%%%%%%%%%%%%%%%%%%%%%%%%%%%%%%%%%%%%%%%%%%%%%%%%%%%%%%%%%%%%%%%%%%%%%%%%%%%%%%%%%%%%%%%%%%%%%%%%%%%%%%%%%%%
 
 \begin{proposition}\label{pr 2.1}
 Let  $T\in(0,\frac N2)$ and  $f\in \cC^1(\R^N)$ such that $f,\, |\nabla f|\in   L^1(\R^N).$ 
Then the function $u_f$ defined by
\bel{Ex6*}u_f(t,x)=\int_{\R^N} \cP_{\ln}(t,x-y) f(y)\dd y, \ee
  is a strong solution of (\ref{eq 1.1}). 
 
\nind (i) Since  the mapping: $f\mapsto   \cP_{\ln}(t,\cdot)\ast f$ is increasing and $\cP_0(t)>0$, if assume furthermore that  $f\neq 0$ is nonnegative and 
 $$  f(x)\leq M(1+|x|)^{\tau}\quad {\it for}\ \  x\in\R^N$$
 for some $M>0$ and $\tau<-N$, then  there exists $c>2$    such that for any $(t,x)\in (0,T)\times\R^N$,
   \begin{equation}\label{e1 2.2}
\frac{1}{ct}  \cP_0(t)  (1+|x|)^{2t-N}\leq  u_f(t,x)   \leq  \frac{c}{t} \cP_0(t)  (1+|x|)^{2t-N},
 \end{equation}
 where $\cP_0$ is defined in (\ref{e 2.1}). 
 
\nind (ii) Moreover, if it is assumed that
 $$\frac1M(1+|x|)^{\tau}\leq f(x)\leq M(1+|x|)^{\tau}\quad {\it for}\ \  x\in\R^N$$
 for some $\tau\in(-N, -2T]$
and some $M>1$,  then  there exists $c>1$   such that for any $(t,x)\in (0,T)\times\R^N$,
   \begin{equation}\label{e1 2.2-tau}
  \frac{1}{ct}  \cP_0(t) (1+|x|)^{ 2t+\tau  } \leq  u_f(t,x)\leq  \frac{c}{t}  \cP_0(t)   (1+|x|)^{  2t+ \tau}.
 \end{equation}
 
\nind (iii) At last,  if assume also that, for some $M>1$, 
 $$\frac1M(1+|x|)^{-N}\leq f(x)\leq M(1+|x|)^{-N}\quad {\it for}\ \  x\in\R^N,$$
  then  there exists $c>1$    such that for any $(t,x)\in (0,T)\times\R^N$,
   \begin{equation}\label{e1 2.2-N}
\frac{1}{ct}  \cP_0(t)  (1+|x|)^{ 2t-N  } \ln(e+|x|)\leq  u_f(t,x)\leq \frac{c}{t}  \cP_0(t)  (1+|x|)^{ 2t-N  } \ln(e+|x|).
 \end{equation}
 \end{proposition}
 %%%%%%%%PROOF%%%%%%%%%%%%%%%%%%%%%%%%%%%%%%%%%%%%%%%%%%%%%%%%%
\noindent \proof First we note that if 
\begin{equation}\label{conv}
u_f(t,x)=\int_{\R^N} \cP_{\ln}(t,z) f(x-z)\dd z =\cP_{0}(t)I_{2t}\ast f(x),
 \end{equation}
 then classicaly
$$\widehat{u_f}(t,\xi)=\cP_{0}(t)\widehat{I_{2t}\ast f}(\xi)=\cP_{0}(t)4^t\gp^{\frac N2}\frac{\Gg(t)}{\Gg(\frac N2-t)}|\xi|^{-2t}\hat f(\xi)=|\xi|^{-2t}\hat f(\xi),
$$
because of the identity (\ref{e 2.1}). Therefore 
\begin{equation}\label{Four}\partial_t\widehat{u_f}(t,\xi) =-2\ln|\xi|\hat f(\xi)\ \ \Longrightarrow \ \ \cF\left(\partial_tu_f+ \lnlap u_f\right)=0.
 \end{equation}
Thus the equation is verified in the space $\CS'(\BBR^N)$ (see \cite{Sch} for the duality between $\CS(\BBR^N)$ and $\CS'(\BBR^N)$ and the Fourier transform therein).
\smallskip

  (i)  Since $u_f$ is defined by (\ref{Ex6*}) 
and $\nabla f\in L^1(\BBR^N)\cap \cC(\BBR^N)$, for all integers $1\leq i\leq N$  the expressions
$$\partial_{x_i}u_f(t,x)=\int_{\R^N} \cP_{\ln}(t,z) \partial_{x_i}f(x-z)\dd z ,$$
are well-defined and are continuous functions of $(x,t)$.  
As a consequence,  $u_f(t,\cdot)\in \cC^1(\R^N)$.

 For the time-derivative, we note that  
$$\partial_t \cP_{\ln}(t,x)= \cP'_0(t) |x|^{2t-N}+2\cP_0(t)|x|^{2t-N} \ln|x|  ,   $$
then
\begin{eqnarray*}
\partial_t u_f(t,x)&=&  \int_{\R^N} \partial_t \cP_{\ln}(t,z)  f(x-z)\dd z 
\\&=&  \cP'_0(t) \int_{\R^N}|z|^{2t-N}f(x-z)\dd z+ 2\cP_0(t)\int_{\R^N}|z|^{2t-N}(\ln|z|) f(x-z)\dd z,
\end{eqnarray*}
which is continuous at any point  $(t,x) \in (0,T)\times\R^N$.  Next we prove that
\begin{equation}\label{Four+}\lim_{t\to0^+} u_f(t,x)=f(x)\quad\text{for any }x\in\BBR^N,  \end{equation}
by adapting the proof of Proposition \ref{pr 1.1}-(iii). For $x\in \R^N$, $R>|x|+1$ we write
$$
u_f(t,x)=\int_{B_R(x)}\cP_{\ln}(t,z)  f(x-z)\dd z+\int_{B^c_R(x)}\cP_{\ln}(t,z)  f(x-z)\dd z.
$$
Since $f$ is uniformly continuous on $B_R(x)$, we have again
$$\lim_{t\to0^+}\int_{B_R(x)}\cP_{\ln}(t,z)  f(x-z)\dd z=f(x).
$$
By H\"older and Gagliardo-Nirenberg inequalities
$$\BA{lll}\dsps\left|\int_{B^c_R(x)}\cP_{\ln}(t,z)  f(x-z)\dd z\right|\leq \cP_{0}(t)\left(\int_{B^c_R(x)}|z|^{(2t-N)N}\dd z\right)^{\frac 1N}
\left(\int_{B^c_R(x)}|f(x-z)|^{\frac{N}{N-1}}\dd z\right)^{\frac{N-1}{N}}\\[4mm]
\phantom{\dsps\left|\int_{B^c_R(x)}\cP_{\ln}(t,z)  f(x-z)\dd z\right|}
\dsps\leq c\cP_{0}(t)\left(\int_{B^c_1}|z|^{(2t-N)N}\dd z\right)^{\frac 1N}\int_{\BBR^N}|\nabla f|\dd z\\[4mm]
\phantom{\dsps\left|\int_{B^c_R(x)}\cP_{\ln}(t,z)  f(x-z)\dd z\right|}
\dsps\leq c'\cP_{0}(t)\int_{\BBR^N}|\nabla f|\dd z.
\EA$$
Since $\cP_{0}(t)\to 0$ as $t\to 0$, (\ref{Four+}) follows.  This implies that
$$u_f\in \cC\Big([0,\tfrac2N)\times \R^N \Big)\cap \cC^1\Big((0,\tfrac2N)\times \R^N \Big).$$
Since (\ref{Four}) holds, it follows that $u_f$ satisfies (\ref{L1}); as a consequence it is a strong solution of (\ref{eq 1.1}).  

  Let $R>1$ such that   ${\rm supp}(f)\subset \bar B_R$, then if $|x|>2R$, 
$$\frac12|x| \leq |x-y|\leq  3|x| \quad{\rm for}\ \  y\in B_R,$$
  and if $f$ is a nonnegative non identically zero function, 
\begin{eqnarray*}
 \frac{|x|^{2t-N}}{2^{N-2t}}\| f\|_{L^1(\R^N)}<  \int_{\R^N}|x-y|^{2t-N} f(y) \dd y  =    \int_{B_R} |x-y|^{2t-N}  f(y)  \dd y 
 \leq      \frac{|x|^{2t-N}}{3^{2t-N}} \| f\|_{L^1(\R^N)},
\end{eqnarray*}
and if $|x|\leq 2R$, we have that 
\begin{eqnarray*}
 \int_{\R^N}|x-y|^{2t-N} f(y) \dd y   \leq    \norm f_{L^\infty(\R^N)}\int_{B_{3R}}|z|^{2t-N}  \dd z =\frac{cR^{2t}}{t}\norm f_{L^\infty(\R^N)}
\end{eqnarray*}
for some $c=c(N)>0$.   

  If we assume that $0\leq f(x)\leq M(1+|x|)^\gt$ and $\gt<-N$, then the right-hand side part of (\ref{e1 2.2}) follows directly from 
Lemma\ref{lm a-1}-(iii). For the left-hand side, we can suppose that $\norm{f{\bf 1}_{B_R}}_{L^1}>0$ for some $R>0$ and we have for 
$|x|\geq 2R$, 
$$\frac{|x|^{2t-N}}{2^{N-2t}}\norm{f{\bf 1}_{B_R}}_{L^1(\BBR^N)}\leq 
\int_{B_R}|x-y|^{2t-N}f(y)dy.
$$
This implies (\ref{e1 2.2}) up to changing the constant $c$.
\smallskip

  (ii)  If we assume  $\frac1M(1+|x|)^{\tau} \leq f(x)\leq M(1+|x|)^{\tau}$,  then (\ref{e1 2.2-tau}) follows from Lemma \ref{lm a-1}-(i).

\smallskip

  (iii) If we assume  $\frac1M(1+|x|)^{-N} \leq f(x)\leq M(1+|x|)^{-N}$, then (\ref{e1 2.2-N}) follows from Lemma \ref{lm a-1}-(ii).  
\hfill$\Box$\medskip 

%%%%%%%%%%%%%%%%%%%%%%%%%%%%%%%%%%%%%%%%%%%%%%%%%%%%%%%%%%%%%%%%%%%%%%%%%%%%%%%%%%%%%%%%%%%%%%%%%%%%%%%%%%%%%%%%%%%%%%%%%%%%%%%%%%%%%%%%%%%%%%
  
  The next results will be used later on

   \begin{lemma}\label{re 2.1} There holds
        \begin{equation}\label{G0}
     \lim_{t\to 0}t\Gg(t)=1,
   \end{equation}
     \begin{equation}\label{G1}
     \lim_{t\to 0}t^2\Gg'(t)=-1
   \end{equation} 
   and
     \begin{equation}\label{G2}
     \lim_{t\to 0}t^3\Gg''(t)=2.
     \end{equation} 
In particular
\begin{equation}\label{G3} \BA{lll}
\cP_0'(0)=
  \pi^{-\tfrac{N}{2}}\Gamma(\tfrac N2)\,\text{ and }\;
\cP_0''(0)=2\pi^{-\tfrac{N}{2}}\left(\Gamma(\tfrac N2)-\Gamma'(\tfrac N2)\right).
\EA
\end{equation} 
   \end{lemma}
   \proof 
Since we have 
$$\Gg(t)=\int_0^\infty z^{t-1}e^{-z}dz=\left[\frac{z^t}{t}e^{-z}\right]_{z=0}^\infty+\frac{1}{t}\int_0^\infty
z^te^{-z}dz,
$$
we obtain (\ref{G0}). For $\Gg'$ there holds
$$\BA{lll}
\dsps\Gg'(t)=\int_0^\infty z^{t-1}\ln ze^{-z}dz=\left[\left(\frac{z^t}{t}\ln z-\frac{z^t}{t^2}\right)e^{-z}\right]_{z=0}^\infty+\int_0^\infty
\left(\frac{z^t}{t}\ln z-\frac{z^t}{t^2}\right)e^{-z}dz\\[4mm]
\phantom{\Gg'(t)}\dsps=\frac{1}{t^2}\int_0^\infty
\left(tz^t\ln z-z^t\right)e^{-z}dz.
\EA$$
Then (\ref{G1}) follows. For $\Gg''$ we also have
$$\BA{lll}
\dsps\Gg''(t)=\int_0^\infty z^{t-1}\ln^2 ze^{-z}dz=\left[\left(\frac{z^t}{t}\ln^2 z-2\frac{z^t}{t^2}\ln z+\frac{2z^t}{t^3}\right)e^{-z}\right]_{z=0}^\infty
\\[4mm]
\phantom{----------------\Gg'(t)}\dsps
+\int_0^\infty
\left(\frac{z^t}{t}\ln^2 z-2\frac{z^t}{t^2}\ln z+\frac{2z^t}{t^3}\right)e^{-z}dz\\[4mm]
\phantom{\Gg'(t)}\dsps=\frac{1}{t^3}\int_0^\infty
\left(t^2z^t\ln^2 z-2tz^t\ln z+2z^t\right)e^{-z}dz.
\EA$$
Hence (\ref{G2}) follows by integration by parts. The first formula in (\ref{G3}) follows directly from (\ref{G1}). The second is a little more involved and left to the reader.
 \hfill$\Box$\medskip 
 
  Given $t_0\in(0,\frac{N}{2})$, 
  we consider the backward problem 
  \begin{equation}\label{eq 2.1}
\left\{ \arraycolsep=1pt
\begin{array}{lll}
\displaystyle \partial_tu- \lnlap u=0\quad \  &{\rm in}\ \   [0,t_0)\times \R^N,\\[2mm]
 \phantom{ \lnlap \     }
\displaystyle   u(t_0,\cdot)=g\quad \ &{\rm{in}}\  \   \R^N.
\end{array}
\right.
\end{equation} 
 The existence of strong solutions of (\ref{eq 2.1}) is established below.
 \begin{corollary}\label{cr 2.1} 
Let $t_0\in(0,\frac{N}{2})$ and for $g\in \cC^1_c(\R^N)$, we denote
 \begin{equation}\label{e dual}
 v_g(t,x)=\int_{\R^N} \cP_{\ln}(t_0-t,x-y) g(y)\dd y, 
 \end{equation}
 then $v_g$ is a strong solution of  (\ref{eq 2.1}),
 $$v_g\in \cC\Big([0,t_0]\times \R^N \Big)\cap \cC^1\Big((0,t_0]\times \R^N \Big), $$
 and we have, for some $c>0$ dependent of $\| f\|_{L^1(\R^N)} $ and $\| f\|_{L^\infty(\R^N)}$
 \begin{equation}\label{e 2.2}
 |v_g(t,x)|\leq c \frac{\cP_0(t_0-t)}{ t_0-t }  (1+|x|)^{2(t_0-t)-N}\quad {\it for}\ (t,x)\in (0,t_0)\times\R^N.
 \end{equation}
 
\end{corollary}
 \noindent\proof 
 Note that the backward diffusion kernel $\cQ_{\ln}$ is 
$$
 \cQ_{\ln}(t,x)=\cP_{\ln}(t_0-t,x) \quad {\rm in}\ \ [0,t_0)\times \R^N.
$$ 
Then the proof is the same as the one of Proposition \ref{pr 2.1}.  \hfill$\Box$
 %%%%%%%%%%%%%%%%%%%%%%%%%%%%%%%%%%%%%%%%%%%%%%%%%%%%%%%%%%%%%%%%%%%%%%%%%%%%%%%%%%%%%%%%%%%%%%%%%%%%%%%%%%%%%%%%%%%%%%%%%%%%%%%%%%%%%%%%%%%%%%%%%%%%%%%%%%%%%%%%%%%%%%%%%%%%%%%%%%%%%%%%%%%%%%%%%%%%%%
 \subsection{Maximum principles}
 We first recall some maximum principles for the logarithmic Laplacian valid in bounded domains \cite[Corollary 1.9]{CT}:
 \begin{lemma}
\label{sec:main-result-max-principle}
Let $\Omega \subset \R^N$ be a bounded Lipschitz domain. Then $\lnlap$ satisfies the maximum principle in $\Omega$ if one of the following conditions is satisfied:\smallskip

\nind (i) $h_\Omega(x) + \rho_N \geq 0$ on $\Omega$, where
\begin{equation}
  \label{eq:def-h-Omega}
h_\Omega(x)=  c_{_N} \Bigl( \int_{B_1(x)\setminus \Omega} \frac{1}{|x-y|^N}dy - \int_{\Omega \setminus B_1(x)} \frac{1}{|x-y|^N}dy\Bigr);
\end{equation}
\nind (ii) $|\Omega| \le 2^N \exp \Bigl(\frac{N}{2} \Bigl(\psi(\frac{N}{2})-\gamma\Bigr)\Bigr)|B_1(0)|$.
\end{lemma}\medskip

The role of $h_\Omega(x) + \rho_N$ is clear since the  following expression,
 \begin{equation}
 \lnlap  u(x) = c_{_N} \int_{\Omega} \frac{ u(x) -u(y)}{|x-y|^{N} } dy - c_{_N} \int_{\R^N \setminus \Omega} \frac{u(y)}{|x-y|^{N}}\,dy + \Big(h_\Omega(x)+\rho_N\Big) u(x)\label{representation-regional}
\end{equation}
 is an   alternative definition of the logarithmic Laplacian.\smallskip
 
For the associated diffusion operator,  we have the following maximum principle:
 
 \begin{corollary}\label{cr new}
 Let $\Omega \subset \R^N$ be a bounded Lipschitz domain. Assume $v \in
\cC([0,T)\times \bar{\Omega} )$ such that $\partial_tv \in \cC\Big((0, T) \times \bar{\Omega} \Big)$ satisfies pointwise 
\begin{equation}\label{eq 2.2}
\left\{ \arraycolsep=1pt
\begin{array}{lll}
\displaystyle \partial_t v+ \lnlap v\geq 0\quad \  &{\rm in}\ \   (0,T)\times \Omega,\\[2mm]
 \phantom{ \lnlap ---   }
\displaystyle   v  \geq 0 \quad \ &{\rm in}\  \   \Big((0,T)\times (\R^N\setminus\Omega)\Big)\cup \Big(\{0\}\times \R^N \Big).
\end{array}
\right.
\end{equation}
 If $h_\Omega + \rho_N \geq0$ on $\Omega$, then
   $$  v\geq 0 \quad{\it in}\ \, [0, T)\times \Omega.$$
\end{corollary}
\noindent\proof We proceed by contradiction, assuming that
for   any $T^{\prime} \in(0, T)$, we have that  
\begin{equation}  \label{e 22-0}
v\left(t_0, x_{0}\right):=\min _{ \left[0, T^{\prime}\right]\times \bar{\Omega}} v(t,x)< 0.
\end{equation} 
Clearly  $\left(t_0, x_{0}\right)$   lies in $(0, T)\times  \Omega$.  Therefore   $v_{t}\left(x_{0}, t_{0}\right) \leq 0$ and thus 
\begin{equation}  \label{e 22-00}
\lnlap v(t_0,x_0) \geq 0.
\end{equation} 
This yields a contradiction by using the expression (\ref{representation-regional}), since
$$\BA{lll}\dsps
  \lnlap v(t_0, x_0) = c_{_N}\int_{\Omega} \frac{ v(t_0, x_0) -v(t_0, y)}{|x_0-y|^{N} } dy  
- c_{_N} \int_{\R^N \setminus \Omega} \frac{v(t_0, y)}{|x_0-y|^{N}}\,dy + \Big(h_\Omega(x_0)+\rho_N\Big)v(t_0, x_0)<  0,
\EA
$$
from  (\ref{e 22-0}) and $ \Big(h_\Omega(x_0)+\rho_N\Big) v(t_0, x_0)\geq 0$ from the assumption on     $h_\Omega+\rho_N$.
\hfill$\Box$\medskip

In the following result we get rid of the  restriction $h_\Omega + \rho_N \geq0$ provided the boundary value of $v$ is bounded from below.
 \begin{lemma}\label{lm comp}
 Let $\Omega \subset \R^N$ be a bounded Lipschitz domain and $v \in
\cC([0,T)\times \bar{\Omega} )$ such that $\partial_tv \in \cC\Big((0, T) \times \bar{\Omega} \Big)$ which satisfies pointwise the following inequalities
\begin{equation}\label{eq 2.2-stronger}
\left\{ \arraycolsep=1pt
\begin{array}{lll}
\displaystyle \partial_t v+ \lnlap v\geq 0\quad \  &{\rm in}\ \   (0,T)\times \Omega,\\[2mm]
 \phantom{ \lnlap ---   }
\displaystyle   v  \gvertneqq 0 \quad \ &{\rm in}\  \    (0,T)\times (\R^N\setminus\Omega), \\[2mm]
 \phantom{ \lnlap -    }
  v(0,\cdot) \geq c \quad \ &{\rm in}\  \  \Omega
\end{array}
\right.
\end{equation}
for some $c>0$. 
Then
   $$  v> 0 \quad{\it in}\ \, [0, T)\times \Omega.$$
\end{lemma}
\noindent\proof %For $\epsilon>0$, let $\Omega_\epsilon=\{x\in\Omega:{\rm dist}(x,\partial\Omega)<\epsilon\}$.  We need show that for any $\epsilon>0$
%$$  v> 0 \quad{\rm in}\ \, [0, T)\times \Omega_\epsilon.$$
We proceed by contradiction. Let $t_0>0$ be the smallest $t>0$ such that $v(t_0,x_0)=0$ for some $x_0\in \Omega$.
Then  $v(t , x)>0$ for all  $(t,x)\in(0,t_0)\times \Omega$.
Thus $\partial_tv(t_0, x_{0})\leq 0$  and 
$$
 \lnlap v(t_0, x_{0})= -c_{_N}\int_{\R^N} \frac{   v(t_0, y)}{|x_0-y|^{N} } dy \geq 0,
$$
where the equality holds only if   $u(t_0,\cdot) \equiv  0$ a.e. in $\R^N$. Since $v$ is continuous in $([0,T)\times \bar{\Omega} )$ nonnegative and not identically $0$ on  the boundary, we are led to a contradiction.
\hfill$\Box$\medskip

%%%%%%%%%%%%%%%%%%%%%%%%%%%%%%%%%%%%%%%%%%%%%%%%%%%%%%%%%%%%%%%%%%%%%%%%%%%%%%%%%%%%%%%%%%%%%%%%%%%%%%%%%%%%%%%%%%%%%%%%%%%%%%%%%%%%%%%%%%%%%%
 %%%%%%%%%%%%%%%%%%%%%%%%%%%%%%%%%%%%%%%%%%%%%%%%%%%%%%%%%%%%%%%%%%%%%%%%%%%%%%%%%%%%%%%%%%%%%%%%%%%%%%%%%%%%%%%%%%%%%%%%%%%%%%%%%%%%%%%%%%%%%%
%%%%%%%%%%%%%%%%%%%%%%%%%%%%%%%%%%%%%%%%%%%%%%%%%%%%%%%%%%%%%%%%%%%%%%%%%%%%%%%%%%%%%%%%%%%%%%%%%%%%%%%%%%%%%%%%%%%%%%%%%%%%%%%%%%%%%%%%%%%%%%

 \section{ The Cauchy problem}
 
 \subsection{Uniqueness of weak Solutions}

Our uniqueness  result for weak solutions is the following.
%%%%%%%%%%%%%%%%%%%%%%%%%%%%%%%%%%%%%%%%%%%%%%%%%%%%%%%%%%%%%%%%%%%%%%%%%%%%%%%%%%%%%%%%%%%%%%%%%%%%%%%%%%%%%%%%%%%%%%%%%%%%%%%%%%%%%%%%%%%%%%
%%%%%%%%%%%%%%%%%%%%%%%%%%%%%%%%%%%%%%%%%%%%%%%%%%%%%%%%%%%%%%%%%%%%%%%%%%%%%%%%%%%%%%%%%%%%%%%%%%%%%%%%%%%%%%%%%%%%%%%%%%%%%%%%%%%%%%%%%%%%%%
\begin{theorem}\label{thm2}
Let $T\in(0,\frac{N}{2})$ and $u$ be a weak solution of the 
equation  with zero initial data
 \begin{equation}\label{eq 0.1-i0}
\left\{ \arraycolsep=1pt
\begin{array}{lll}
\partial_tu+ \lnlap u=0\quad \  &{\rm in}\ \   (0,T)\times \R^N,\\[2mm]
 \phantom{ \lnlap \ \   }
  u(0,\cdot)=0\quad \ &{\rm in}\  \   \R^N.
\end{array}
\right.
\end{equation}
 Then $u(t,x)=0$ for every $t \in(0, T)$ and for almost every $x \in \mathbb{R}^N$.
\end{theorem}
%%%%%%%%%%%%%%%%%%%%%%%%%%%%%%%%%%%%%%%%%%%%%%%%%%%%%%%%%%%%%%%%%%%%%%%%%%%%%%%%%%%%%%%%%%%%%%%%%%%%%%%%%%%%%%%%%%%%%%%%%%%%%%%%%%%%%%%%%%%%%%
\noindent\proof  {\it Step 1}. We fix   $\theta\in \cC^1_c\left(\R^N \right)$ with compact support in $\bar B_{R_0}$ for some $R_0>0$.
We aim to prove that
$$
\int_{\mathbb{R}^N} u\left(t_{0}, x\right) \theta(x) \dd x=0.
$$
For $t \in\left[0, t_{0}\right)$, we define
$$
\varphi_\theta=   \cP_{\ln} (t_{0}-t,\cdot)\ast \theta.  $$
By Corollary \ref{cr 2.1}, function $\varphi_\theta$
 is a strong solution of 
 \begin{equation}\label{eq 2.1-d}
  \left\{ \arraycolsep=1pt
\begin{array}{lll}
\partial_tv- \lnlap v=0\quad \  &{\rm in}\ \   [0,t_0)\times \R^N,\\[2mm]
 \phantom{ \lnlap \     }
  v(t_0,\cdot)=\theta\quad \ &{\rm in}\  \ \,  \R^N.
\end{array}
\right.
 \end{equation}
Then
$$
|\varphi_\theta(t,x)|=\left| \Big(\cP_{\ln}(t_{0}-t,\cdot)\ast \theta\Big) (x)\right| \leqq c\frac{\cP_0(t_0-t)}{ t_0-t }  (1+|x|)^{2(t_0-t)-N},
$$
where $\frac{\cP_0(t_0-t)}{ t_0-t }$ is bounded as $t\to t_0$ and $c>0$ depends on $R_0$, $t_{0}$    and  $ \|\theta\|_{L^1(\R^N)},\   \|\theta\|_{L^\infty(\R^N)}.$

\nind Applying (\ref{e dual}) to the derivatives of $\theta \in  \cC _c^1\left(B_{R_{0}}\right)$, we also have
\begin{eqnarray*}
|\nabla \varphi_\theta(t,x)| &=&  \left|\Big((\nabla \theta) \ast \cP_{\ln}({t_{0}-t},\cdot)\Big)(t,x) \right|
\\[2mm]&= & \left|\Big( \cP_{\ln}({t_{0}-t},\cdot)\ast (\nabla \theta)\Big)(t,x) \right|   
\\[2mm]&\leq  & c\frac{\cP_0(t_0-t)}{ t_0-t }  (1+|x|)^{2(t_0-t)-N},
\end{eqnarray*}
where   $c>0$ depends on $R_0$, $t_{0}$    and  $ \|\nabla \theta\|_{L^1(\R^N)},\   \|\nabla \theta\|_{L^\infty(\R^N)}.$

 For $R>2 R_{0}$ we define 
$$\eta_{_R}(x)=\eta_0\Big(\tfrac{|x|}{R}\Big)\qquad{\rm and}\qquad
\psi_{_R}(x, t):=\varphi_\theta(x, t) \eta_{_R}(x),
$$
where $\eta_0:[0,+\infty) \to[0,1]$ is a smooth, non-increasing  function   such that
$$
 \eta_0=1\quad  {\rm in} \ \ [0,\tfrac12]\qquad{\rm and}\qquad \eta_0=0\quad  {\rm in} \ \ [1,+\infty).
$$
As a test function in (\ref{eq 0.1-i0}) we take $\psi_{_R}$. Since $\varphi_\theta$ is a strong solution,  we have that  
\begin{eqnarray*}
\Big| \int_{\mathbb{R}^N} u\left(t_0,x\right) \theta(x) \psi_{R}(x) \mathrm{d} x \Big|  
&=&\Big|\int_{0}^{t_{0}} \int_{\mathbb{R}^N}\left[u \eta_{R}(x) \partial_t\varphi_\theta  -u \lnlap (\psi_{_R}(t,x))\right] \mathrm{d} x \mathrm{~d} t\Big| \\[2mm]
&=&\Big|\int_{0}^{t_{0}} \int_{\mathbb{R}^N} u  \lnlap \Big(\varphi_\theta(1-\eta_R)\Big)  \mathrm{d} x \mathrm{~d} t\Big| \\[2mm]
&\leqq & \int_{0}^{t_{0}} \int_{\mathbb{R}^N} | u | \,\Big|\lnlap \Big(\varphi_\theta(1-\eta_R)\Big)\Big|  \mathrm{d} x \mathrm{~d} t.
\end{eqnarray*}
The main point is to show that
 \begin{equation}\label{e 2.1-1}
\lim _{R \rightarrow \infty} \int_{0}^{t_{0}} \int_{\mathbb{R}^N} | u | \, \Big|\lnlap \Big(\varphi_\theta(1-\eta_R)\Big)\Big|  \mathrm{d} x \mathrm{~d} t=0.
\end{equation}
\smallskip

\nind {\it Step 2. We claim that there exists a positive constant $c$ independent of $t$ such that for $R>\!\!> 1$ 
 \begin{equation}\label{c 2.1-1}
 \Big|\lnlap \Big(\varphi_\theta(1-\eta_R)\Big)(t,x)\Big| \leq  c\min\Big\{R^{2(t_0-t)-N}\ln R, |x|^{2(t_0-t)-N}\ln |x|\Big\},
 \end{equation}
 for any $(t,x)\in (0,T)\times \R^N$.} 
\smallskip

 (i) 
In the sequel we take $R>\!\!> 1$.  If $x\in B_{\frac R4}$, then  $B_1(x)\subset B_{\frac{R}{4}}(x)\subset B_{\frac R2}$, thus
  $$|x-y|>\frac{1+|y|}{4}\quad {\rm for}\ \ y\in \R^N\setminus B_{\frac R2}.
  $$
Since  $1-\eta_R=0$ in $B_{\frac R2}$, we have that
 \begin{eqnarray*}
\big|\lnlap \big(\varphi_\theta(1-\eta_R)\big)(t,x)\big|  &= &c_{_N} \Big|\int_{\R^N\setminus B_{\frac R2}} \frac{-\varphi_\theta(t,y)(1-\eta_R(y))}{|x-y|^N}\dd y\Big| \\[2mm] 
&\leqq & 4^N  c_{_N}\int_{\R^N\setminus B_{\frac R2}} \frac{|\varphi_\theta(t,y)|}{(1+|y|)^N}\dd y
 \\[2mm] 
&\leqq & c\frac{\cP_0(t_0-t)}{ t_0-t } \int_{\R^N\setminus B_{\frac R2}}  (1+|y|)^{2(t_0-t)-2N} \dd y
\\[2mm] &\leqq & c' R^{2(t_0-t)-N},
\end{eqnarray*} 
the last inequality, being a consequence of the  fact that $\dsps \lim_{t\to t_0}\frac{\cP_0(t_0-t)}{ t_0-t }$ exists.   Therefore, for any $t\in[0,t_0]$ and  $x\in   B_{\frac R4}$,  
 \begin{equation}\label{e 2-p1} 
 \big|\lnlap \big(\varphi_\theta(1-\eta_R)\big)(t,x)\big|\leq c' R^{2(t_0-t)-N}.
 \end{equation}

 (ii) Next, if $x\in B_{2R}\setminus B_{\frac R4}$, we note that  $\varphi_\theta(1-\eta_R)$ is $C^1$, then
 \begin{eqnarray*}
 \big|\lnlap \big(\varphi_\theta(1-\eta_R)\big)(t,x)\big|  &\leq &c_{_N}\Big|
\int_{B_1(x)} \frac{\big(\varphi_\theta(1-\eta_R)\big)(t,x)-\big(\varphi_\theta(1-\eta_R)\big)(t,y)}{|x-y|^N} \dd y\Big|
\\[2mm]\\&&+
 c_{_N}\Big|\int_{\R^N\setminus \big(B_{\frac R2} \cup B_1(x)\big) } \frac{-\varphi_\theta(1-\eta_R)\big)(t,y) }{|x-y|^N}\dd y\Big|  +|\rho_{_N}\varphi_\theta (1-\eta_R) |. 
 \end{eqnarray*} 
From the estimate on $\varphi_\gth$ there holds if  $\frac R4\leq |x|\leq 2R$, 
$$\left|(\varphi_\theta(1-\eta_R))(t,x) \right|\leq   c R^{2(t_0-t)-N}.$$
Furthermore
\begin{eqnarray*}  
\Big| \int_{B_1(x)} \frac{\big(\varphi_\theta(1-\eta_R)\big)(t,x)-\big(\varphi_\theta(1-\eta_R)\big)(t,y)}{|x-y|^N} \dd y\Big|&\leq& \sup_{y\in B_1(x)}|\nabla\varphi_\theta(t,y)||\nabla\eta_{_R}(y)| 
\\[2mm]&\leq& \frac{c}{R}  (1+|x|)^{2(t_0-t)-N}\chi_{[\frac R2,R]}(|x|)
\\[2mm]&\leq& c R^{2(t_0-t)-N-1},
 \end{eqnarray*} 
 and
\begin{eqnarray*}  
  \Big|\int_{\R^N\setminus \big(B_{\frac R2} \cup B_1(x)\big)} \frac{-\big(\varphi_\theta(1-\eta_R)\big)(t,y) }{|x-y|^N}\dd y\big| 
 &\leq& \int_{\R^N\setminus \big(B_{\frac R2} \cup B_1(x)\big)} \frac{|\varphi_\theta (t,y)| }{|x-y|^N}\dd y
\\[2mm]&\leq&  c \int_{\R^N\setminus \big(B_{\frac R2} \cup B_1(x)\big) } \frac{ |y|^{2(t_0-t)-N}}{|x-y|^N}\dd y.
%\\[2mm]&=&c_1  \int_{B_{\varrho}(x) \setminus B_1(x) } \frac{ |y|^{2(t_0-t)-N}}{|x-y|^N}\dd y+c_1  \int_{\R^N\setminus \Big(B_{\frac R2} \cup B_\varrho(x)\Big) } \frac{ |y|^{2(t_0-t)-N}}{|x-y|^N}\dd y.
%\\[2mm]&=&c_1 |x|^{2(t_0-t)-N} \int_{\R^N\setminus B_{\frac R2}} \frac{ |z|^{2(t_0-t)-N}}{|{\bf e}_x-z|^N}\dd z
 \end{eqnarray*} 
Since $R>\!\!> 1$, we have 
\begin{eqnarray*} 
\int_{B_{R/8}(x) \setminus B_1(x) }   \frac{ |y|^{2(t_0-t)-N}}{|x-y|^N}\dd y&\leq&  (R/2)^{2(t_0-t)-N}   \int_{B_{R/8}(x) \setminus B_1(x) }   \frac{1}{|x-y|^N}\dd y
\\[2mm]& \leq &c R^{2(t_0-t)-N} \ln R.
 \end{eqnarray*} 
Clearly
\begin{eqnarray*} 
 \int_{B_{4R}\setminus \big(B_{\frac R2} \cup B_{R/8}(x)\big) } \frac{ |y|^{2(t_0-t)-N}}{|x-y|^N}\dd y
 \leq c R^{2(t_0-t)-N},
 \end{eqnarray*} 
 and
 \begin{eqnarray*} 
 \int_{\R^N \setminus B_{4R} } \frac{ |y|^{2(t_0-t)-N}}{|x-y|^N}\dd y
 \leq \int_{\R^N \setminus B_{4R} }   |y|^{2(t_0-t)-2N} \dd y =c R^{2(t_0-t)-N}.
 \end{eqnarray*} 
   Thus,  for any $t\in[0,t_0]$ and  $x\in B_{2R}\setminus B_{\frac R4}$, we have that 
 \begin{equation}\label{e 2-p2} 
 \big|\lnlap \big(\varphi_\theta(1-\eta_R)\big)(t,x)\big|\leq c R^{2(t_0-t)-N} \ln R.
 \end{equation}

 (iii) Finally, for $x\in \R^N \setminus B_{2R}$, 
 \begin{eqnarray*}
 \Big|\lnlap \big(\varphi_\theta(1-\eta_R)\big)(t,x)\Big|  &\leq &c_{_N}\Big|
\int_{B_1(x)} \frac{\varphi_\theta (t,x)-\varphi_\theta (t,y)}{|x-y|^N} \dd y\Big|
\\[2mm]\\&&+
 c_{_N}\Big|\int_{\R^N\setminus \big(B_{\frac R2} \cup B_1(x)\big) } \frac{-\big(\varphi_\theta(1-\eta_R)\big)(t,y) }{|x-y|^N}\dd y\Big|  +|\rho_{_N}\varphi(t,x) |.
  \end{eqnarray*}
We estimate the three terms on the right-hand side of this inequality as follows.
Obviously
 $$|\rho_{_N}\varphi_\theta(t,x) |\leq c|x|^{2(t_0-t)-N},$$
 next
\begin{eqnarray*}  
\Big| \int_{B_1(x)} \frac{\varphi_\theta (t,x)-\varphi_\theta (t,y)}{|x-y|^N} \dd y\Big|&\leq& \sup_{y\in B_1(x)}|\nabla\varphi_\theta(t,y)| 
 \leq  c    (1+|x|)^{2(t_0-t)-N},
 \end{eqnarray*} 
 and finally
\begin{eqnarray*}  
  \Big|\int_{\R^N\setminus \big(B_{\frac R2} \cup B_1(x)\big)} \frac{-\big(\varphi_\theta(1-\eta_R)\big)(t,y) }{|x-y|^N}\dd y\Big| 
 &\leq& \int_{\R^N\setminus \big(B_{\frac R2} \cup B_1(x)\big)} \frac{|\varphi_\theta (t,y)| }{|x-y|^N}\dd y
\\[2mm]&\leq&  c   \int_{\R^N\setminus \big(B_{\frac R2} \cup B_1(x)\big) } \frac{ |y|^{2(t_0-t)-N}}{|x-y|^N}\dd y.
 \end{eqnarray*} 
From these estimates,  we obtain the following series of inequalities (with ${\bf e}_x=\frac x{|x|}$):
\begin{eqnarray*} 
\int_{B_{\frac{|x|}{4}}(x) \setminus B_1(x) }   \frac{ |y|^{2(t_0-t)-N}}{|x-y|^N}\dd y&\leq& 
c |x|^{2(t_0-t)-N} \int_{B_{\frac{|x|}{4}}(x) \setminus B_1(x) }   \frac{1}{|x-y|^N}\dd y 
\\[2mm]&\leq& c |x|^{2(t_0-t)-N} \ln\tfrac{|x|}{4},
 \end{eqnarray*}
 \begin{eqnarray*} 
\int_{B_{\frac{|x|}{4}}(0) \setminus B_{\frac R2}(0)}   \frac{ |y|^{2(t_0-t)-N}}{|x-y|^N}\dd y&\leq& c   |x|^{-N}  \int_{B_{\frac{|x|}{4}}(0) \setminus B_{\frac R2}(0)}     |y|^{2(t_0-t)-N} \dd y
\\[2mm]& \leq& c   |x|^{2(t_0-t)-N}, 
 \end{eqnarray*}
\begin{eqnarray*} 
\int_{B_{2|x|}(x)\setminus \big(B_{\frac{|x|}4}(0) \cup B_{\frac{|x|}4}(x)\big) } \frac{ |y|^{2(t_0-t)-N}}{|x-y|^N}\dd y\phantom{-----------}
 \\[2mm] = |x|^{2(t_0-t)-N} \int_{B_{2}({\bf e}_x)\setminus \big(B_{\frac{1}4}(0) \cup B_{\frac{1}4}({\bf e}_x)\big) }\frac{ |z|^{2(t_0-t)-N}}{|{\bf e}_x-z|^N}\dd z
 \\=c |x|^{2(t_0-t)-N}, \phantom{----------------}
 \end{eqnarray*} 
 and
 \begin{eqnarray*} 
 \int_{\R^N \setminus B_{2|x|}(x) } \frac{ |y|^{2(t_0-t)-N}}{|x-y|^N}\dd y
&=& |x|^{2(t_0-t)-N}  \int_{\R^N \setminus B_{2} }  \frac{ |z|^{2(t_0-t)-N}}{|{\bf e}_x-z|^N}\dd z
\\[1mm]&=& c |x|^{2(t_0-t)-N}.
 \end{eqnarray*} 
   Therefore,  for any $t\in[0,t_0]$ and  $x\in \R^N \setminus B_{2R}$ 
 \begin{equation}\label{e 2-p3} 
 \Big|\lnlap \big(\varphi_\theta(1-\eta_R)\big)(t,x)\Big|\leq c   |x|^{2(t_0-t)-N}\ln|x|.
 \end{equation}

\nind{\it Step 3. End of the proof. }
 Assuming that (\ref{c 2.1-1}) holds, we have 
\bel{XX1-sec3}\BA{lll}\dsps
  \int_{\mathbb{R}^N} | u | \big|\lnlap \big(\varphi_\theta(1-\eta_R)\big)\big|  \mathrm{d} x  
&\leq &\dsps c \int_{\mathbb{R}^N}  | u |  \min\Big\{R^{2(t_0-t)-N}\ln R, |x|^{2(t_0-t)-N}\ln |x|\Big\}\dd x   
\\[4mm]\dsps&\leq &\dsps\int_{\mathbb{R}^N} |u(t,x)| (1+|x|)^{2T-N}  V_{_R}(x) dx,
\EA\ee
where 
\bel{XX2-sec3}\BA{lll}\dsps
V_{_R}(x)&=&(1+|x|)^{N-2T} \min\Big\{R^{2(t_0-t)-N}\ln R, |x|^{2(t_0-t)-N}\ln |x|\Big\} 
\\[2mm]&\leq & \dsps\min\Big\{R^{2(t_0-t)-2T}\ln R, |x|^{2(t_0-t)-2T}\ln |x|\Big\} 
\\[2mm]&\leq &  \dsps R^{2(t_0-t)-2T}\ln R.  
\EA\ee
Clearly the right-hand side of (\ref{XX2-sec3}) tends to $0$ uniformly in $\R^N$ as $R\to\infty$. This implies that right-hand side of (\ref{XX1-sec3}) shares this property.
It follows that 
$$\int_{\mathbb{R}^N} u\left(t_0,x\right) \theta(x) \mathrm{d} x=0
$$
by the dominated convergence theorem, thus $u(t_0,x)=0$. Because $t_0<T$ is arbitrary, $u$ is identically zero.
\hfill$\Box$

%%%%%%%%%%%%%%%%%%%%%%%%%%%%%%%%%%%%%%%%%%%%%%%%%%%%%%%%%%%%%%%%%%%%%%%%%%%%%%%%%%%%%%%%%%%%%%%%%%%%%%%%%%%%%%%%%%%%%%%%%%%%%%%%%%%%
%%%%%%%%%%%%%%%%%%%%%%%%%%%%%%%%%%%%%%%%%%%%%%%%%%%%%%%%%%%%%%%%%%%%%%%%%%%%%%%%%%%%%%%%%%%%%%%%%%%%%%%%%%%%%%%%%%%%%%%%%%%%%%%%%%%%
%%%%%%%%%%%%%%%%%%%%%%%%%%%%%%%%%%%%%%%%%%%%%%%%%%%%%%%%%%%%%%%%%%%%%%%%%%%%%%%%%%%%%%%%%%%%%%%%%%%%%%%%%%%%%%%%%%%%%%%%%%%%%%%%%%%%

\subsection{Strong positive solutions}

In this subsection we give a representation of the positive strong solutions of (\ref{eq 1.1}) as a Poisson type integral.
%and in particular we prove that strong solutions are upper bounds for the kernel convolutions:
\begin{lemma}\label{approx} Let $T\in(0,\frac{N}{2})$, $u$ be a strong solution of (\ref{eq 1.1}) and $\{\gs_\gd\}_{\gd>0}\subset C^\infty_c(\BBR^N)$ a sequence of mollifiers with support in $B_\gd$. Then $u\ast \gs_\gd$ is a strong solution of (\ref{eq 1.1}),  which converges locally uniformly to $u$ 
in $[0,T]\ti\BBR^N$.
\end{lemma}
\noindent\proof We recall that the topological dual of $\CS(\BBR^N)$ endowed with its usual Frechet topology of semi-norms 
$$\norm \varphi_{\ga,R}:=\sup_{x\in\BBR^N}|R(x)D^\ga\varphi(x)|\quad\text{for }\;\ga\in \BBN^N \text{ and }\; R(x)\in\BBR[x]$$
  is  $\CS'(\BBR^N))$, the elements of which  are called slowly increasing distributions  (see e.g. \cite{Sch}). The Fourier transform is an isomorphism of 
$\CS'(\BBR^N)$ into $\CS'(\BBR^N)$ and clearly $Dom(\CL_{\Gd})\subset \CS'(\BBR^N)$: if $\varphi\in \CS(\BBR^N)$
the following duality formulas hold
$$\BA{lll}\dsps \langle \cF(u_t+\cL_\Gd u),\varphi\rangle=\langle u_t,\cF\varphi\rangle+\langle \CF\cL_\Gd u,\varphi\rangle=\frac d{dt}\langle u,\cF\varphi\rangle+\langle \cF\cL_\Gd u,\varphi\rangle\\[2mm]\dsps 
\phantom{\langle \cF(u_t+\cL_\Gd u),\varphi\rangle}=\frac d{dt}\langle \cF u,\varphi\rangle+\langle 2\ln |.|\cF u,\varphi\rangle=\langle \partial_t\hat u+2\ln |.|\hat u,\varphi\rangle.|
\EA$$
Therefore there holds
$$\partial_t\hat u+2 \ln|\xi|\, \hat u=0.
$$
Multiplying by $\hat\gs_\gd$ and using the fact that $\hat\gs_\gd\hat u=\widehat{\gs_\gd\ast u}$, we obtain 
$$\partial_t\widehat{\gs_\gd\ast u}+2\ln|\xi|\,\widehat{\gs_\gd\ast u}=0.
$$
Since $\widehat{\gs_\gd\ast u}$ is $C^\infty$, it follows that $\gs_\gd\ast u$ is a strong solution of (\ref{L1}). As for the convergence of 
$\gs_\gd\ast u$ to $u$, 
it is a classical result since $u\in \cC([0,T)\ti\BBR^N)$.\hfill$\Box$

%%%%%%%%%%%%%%%%%%%%%%%%%%%%%%%%%%%%%%%%%%%%%%%%%%%%%%%%%%%%%%%%%%%%%%%%%%%%%%%%%%%%%%%%%%%%%%%%%%%%%%%%%%%%%%%%%%%%%%%%%%%%%%%%%%%%
\begin{lemma}\label{lm 4.1-im}
Let $T\in(0,\frac{N}{2})$ and $f$ be a nonnegative function in $ \cC(\R^N) \cap L^1(\R^N,(1+|x|)^{2T-N} dx)$.   If $u\geqq 0$ is a strong solution of (\ref{eq 1.1}), then
\begin{eqnarray*}
u(t,x)\geq  \int_{\mathbb{R}^N} \cP_{\ln}(t,x-y) f(y)\dd y   \quad {\it for}\ \, (t,x)\in \big[0, T\big)\times \R^N.
\end{eqnarray*}

\end{lemma}
\noindent\proof For $r>2$ we denote by $\varphi_{r}\in C^{0,1}_c(\R^N)$ the following cut-off function,
$$
\varphi_{r}(x)=
\left\{ \arraycolsep=1pt
\begin{array}{lll}
\displaystyle 1    &{\rm for}\ \   |x|\leq r-1,\\[2mm]
 \phantom{   }
\displaystyle   r-|x|\quad\    &{\rm{for}}\  \ r-1<|x|\leq r,\\[2mm]
\phantom{   } 0  &{\rm{for}}\  \  |x|\geq r,
\end{array}
\right.
 $$
and by $\{\gs_\gd\}_{\gd>0}\subset\C^\infty_c(\BBR^N)$ a sequence of mollifiers with support in $B_\gd$. We set
$f_{r,\gd}=\varphi_r f\ast\gs_\gd$ and define the function $v_{r,\gd}$ by the expression
$$
v_{r,\gd}(t,x)=\int_{\mathbb{R}^N} \cP_{\ln}(t, x-y) f_{r,\gd}(y) \mathrm{d} y=\cP_{\ln}(t, \cdot)\ast  \big(f_{r,\gd}  \big)  (x).
$$
By Proposition \ref{pr 2.1} $v_{r,\gd}$ is a positive strong solution of 
\begin{equation}\label{eq 4.1}
\left\{ \arraycolsep=1pt
\begin{array}{lll}
\partial_t v_{r,\gd}+\lnlap v_{r }=0  &\qquad {\rm in } \ \   (0,T)\times\mathbb{R}^N, \\[2mm]
 \phantom{  --,,}
 v_{r,\gd }(0,\cdot)=f_{r,\gd}   \qquad & \qquad  {\rm in }\ \  \mathbb{R}^N.
\end{array}
\right.
\end{equation}

\nind By (\ref{e1 2.2}) we have for $(t,x)\in (0,T)\times\mathbb{R}^N$ and independently of $\gd$ and $r$, 
\begin{eqnarray}\label{Ve1}
0 <  v_{r,\gd }(t,x)   \leqq  c (\|f\|_{L^\infty}+\|f\varphi_{_r}\|_{L^1}) \frac{ \cP_0(t) }{2t} (1+|x|)^{2t-N}, 
\end{eqnarray}
where, we recall it, $\frac{ \cP_0(t) }{2t}$ is bounded as $t\to0^+$. \smallskip

\nind Let $f_1(x)=(1+|x|)^{\tau_0}$ with $\tau_0\in(-  \frac N2,- T )$
and
$$w_0=\cP_{\ln}(t,\cdot)\ast f_1.$$
By Proposition \ref{pr 2.1}-(ii) there holds
$$w_0(t,x)\geq c(1+|x|)^{2t+\tau_0}\quad {\rm for }\ (t,x)\in (0,T)\times \R^N,$$
where $0>2t+\tau_0>2t-N$. 
From (\ref{Ve1}) for every $\varepsilon>0$ there exists $\rho=\gr(\ge)>r$ such that   
$$
0 < v_{r,\gd}(t,x) \leqq \varepsilon w_0(t,x)  \quad \  {\rm for\ any\ }|x| \geqq \rho,\ t \in[0, T),
$$
and since $u$ is nonnegative in $[0,T)\times \mathbb{R}^N $, $u_\gd:=u\ast \sigma_\gd$ is also nonnegative, thus
$$
0 < v_{r,\gd }(t,x) \leqq \varepsilon w_0(t,x) \leqq \varepsilon w_0(t,x)+u_\gd(t,x) \ \ {\rm for\ any\ }|x| \geqq \rho,\ t \in[0, T).
$$
Clearly
$$
v_{r,\gd }(0,x)=(\varphi_ru\ast\gs_\gd)(0,x) \leqq u_\gd(0,x) \leqq \varepsilon w(0,x)+u_\gd(0,x) \  {\rm for\ any\ }|x| \leqq \rho.
$$
Next we set
$$
w_\gd =u_\gd +\varepsilon w_0-v_{r,\gd}\quad {\rm in}\ \  (0,T)\times \R^N.
$$
Then
$$
w_\gd (0,\cdot)= \varepsilon w_0+((1-\varphi_{r}) f)\ast\gs_\gd\quad {\rm in}\ \   \R^N,
$$
and
$w_\gd (0,\cdot)\geq \varepsilon (1+\rho)^{\tau_0}>0\quad{\rm in}\ \bar B_\rho(0).$
By   Lemma \ref{lm comp} $w_\gd > 0$ in $ (0,T)\times \bar B_\rho(0) $, that is, 
$$
v_{r,\gd}  < \varepsilon w_0+u  \quad \ {\rm for\ }|x| \leqq \rho \ \, {\rm and\ }\,  t \in[0, T) .
$$
Because $\varepsilon>0$ could be chosen arbitrary, the previous inequality implies that,  
$$
v_{r,\gd}(t, x) \leqq u_\gd(t,x) \;\  {\rm for\ all\ } (x,t) \in \mathbb{R}^N\ti [0, T),
$$
and this is independent of $\gr$. Therefore
$$u_\gd(t,x)\geq \CP_0(t)\int_{\BBR^N}|x-y|^{2t-N}((f\varphi_r)\ast\gs_\gd)(y)dy.
$$
When $\gd\to 0$, $u_\gd\to u$ locally uniformly in $[0,T)\ti\BBR^N$ by the previous lemma and $(f\varphi_r)\ast\gs_\gd\to f\varphi_r$ uniformly in $B_r$, thus 
$$\int_{\BBR^N}|x-y|^{2t-N}((f\varphi_r)\ast\gs_\gd)(y)dy\to \int_{\BBR^N}|x-y|^{2t-N}f\varphi_r(y)dy,
$$
by the monotone convergence theorem. Therefore
$$u(t,x)\geq \CP_0(t)\int_{\BBR^N}|x-y|^{2t-N}f\varphi_r(y)dy.
$$
Because $r\mapsto v_r$ is increasing  as $\displaystyle 
\lim _{r \rightarrow \infty} \varphi_r=1$, we conclude that
$$
\lim _{r \rightarrow \infty} \CP_0(t)\int_{\BBR^N}|x-y|^{2t-N}f\varphi_r(y)dy=\CP_0(t)\int_{\BBR^N}|x-y|^{2t-N} f(y)\dd y    \leqq u(t,x)  
$$
by the monotone convergence theorem. This ends the proof.  \hfill$\Box$\medskip 

The next easy to prove result is fundamental in the sequel. 
\begin{corollary}\label{cr 4.1-1}
Let $0<\tau<T<\frac N2$ and $u$ be a strong nonnegative solution of 
\begin{equation}\label{eq 4.2}
\partial_tu+ \lnlap u=0\quad \   {\it in}\ \   (0,T)\times \R^N
\end{equation} 
such that $u(t,.)\in L^1(\R^N,(1+|x|)^{2T-N} dx)$ for all $t\in [0,T]$. Then for any $(t,x) \in (0, T-\tau) \times \mathbb{R}^N $,
\begin{equation}\label{e 4.3}
\int_{\mathbb{R}^N} \cP_{\ln}(t, x-y) u(\tau, y) \dd y \leq  u(t+\tau,x).
\end{equation} 
 As a consequence, we have for every $x \in \mathbb{R}^N$ and $t \in(0, T-\tau)$, 
\begin{equation}\label{e 4.4}
\int_{0}^{T-t} \int_{\mathbb{R}^N} \cP_{\ln}(t,x-y) u(\tau,y ) d y d \tau \leqq \int_{0}^{T-t} u(t+\tau,x) d \tau .
\end{equation} 
\end{corollary}
\proof Let $\gt\in (0,T)$. Then the function $t\mapsto u^{<\gt>}(t,.):=u(t+\gt,.)$ is a strong solution of 
\begin{equation}\label{eq 4.2v}
\partial_tu+ \lnlap u=0\quad \   {\it in}\ \   (0,T-\gt)\times \R^N
\end{equation} 
satisfying $u^{<\gt>}(0,x)=u(\gt,x)$ for all $x\in\BBR^N$. Hence it is larger than $u_{u^{<\gt>}(0,.)}$ in $(0,T-\gt)\times \R^N$. It follows that 
 for all  $(t,x)\in (0,T-\gt)\times \R^N$
\begin{equation}\label{eq 4.2vv}
u(t+\gt,x)=u^{<\gt>}(t,x)\geq u_{u^{<\gt>}(0,.)}(t,x)=\int_{\mathbb{R}^N} \cP_{\ln}(t, x-y) u(\tau, y) \dd y,
\end{equation} 
which ends the proof.\hfill$\Box$
 \medskip

  \medskip

Now we prove our main result:\medskip

\noindent{\it Proof of Theorem \ref{thm1}. }
{\it Proof of (\ref{e 1-0}).}  For $0< T\leq\frac N2$, let $u \geqq 0$ be a strong solution of equation (\ref{eq 1.1}) with initial data $f\in \cC(\BBR^N )\cap L^1(\BBR^N, (1+|x|)^{2T -N} dx)$. Our aim is prove to that  $u$ is a weak solution of (\ref{eq 1.1}).\smallskip
 
  {\it We first prove that for $T'<\frac{N}{2}$,}
 $$u  \in L^{1}\big([0, T'], L^1\big(\mathbb{R}^N, (1+|x|)^{2(T-T')-N}dx\big)  \big).$$
 Take an arbitrary $0<T^{\prime}<T<\frac N2$; then (\ref{e 4.4}) with $x=0$ and $t= T-T'$ shows that 
\begin{eqnarray*}
\int_0^{T'}\int_{\mathbb{R}^N} u(\tau, y)(1+|y|)^{2t-N}  dy d\tau&=&
\int_{0}^{T'} \cP_0(\tau)\int_{\mathbb{R}^N} u(\tau, y)(1+|y|)^{2t-N} \mathrm{d} y \mathrm{~d} \tau \\[1mm]& \leqq &\int_{0}^{T'} u( t+\tau,0) \mathrm{d} \tau 
\\[1mm]
& <& \int_{t}^{T} u(\tau,0) \mathrm{d} \tau
 <+\infty,
\end{eqnarray*}
where $\cP_0$ is bounded in $[0,T]$ because to $T<\frac N2$. \smallskip
 
\nind {\it We claim that  $u \in \cC\left((0, T),\,  L_{  {\rm loc }}^{1} (\mathbb{R}^N )\right)$}.

  For given $\epsilon>0$ and any $R>0$, if $|t_1-t_2|$ small enough,
 there holds that 
 \begin{eqnarray*}
\int_{B_R(0)}  \big| u(t_1, y)-u(t_2, y)\big| \dd y &\leq &  \epsilon\int_{B_R(0)} \dd y \,  
\end{eqnarray*}
by the continuity of $u$.
As a consequence, the strong solution of (\ref{eq 1.1})  is a weak solution.
\\
Define now the function
$$
w =u -\cP_{\ln}\ast f, 
$$
then $w \geqq 0$ by Lemma \ref{lm 4.1-im}.

It is clear that $w$ is a strong solution of (\ref{eq 0.1-i0}), so it is a weak solution by above arguments.  Now we apply  Theorem \ref{thm2} to obtain that 
$$w(x, t)=0 \quad{\rm for\ almost\ every}\ \, (t,x) \in (0,T)\times\mathbb{R}^N.$$
By the continuity of strong solutions,  we have that 
$$u=\cP_{\ln}(t,\cdot)\ast f.$$
\smallskip

\noindent{\it Proof of (\ref{e1-1}).}  Observe that 
\begin{eqnarray*}
\cP_0(t)&=&\pi^{-\frac{N}{2}}\, 4^{-t} \frac{\Gamma(\frac{N-2t}{2})}{\Gamma(t)}
\\[2mm]&=&\pi^{-\frac{N}{2}}\, 4^{-t} \frac{\Gamma(\frac{N-2t}{2}+1)}{\Gamma(t)}\, \frac{2}{N-2t},
\end{eqnarray*}
where 
$$\pi^{-\frac{N}{2}}\, 4^{-t} \frac{\Gamma(\frac{N-2t}{2}+1)}{\Gamma(t)}\to \pi^{-\frac{N}{2}}\, 4^{-\frac{N}{2}}  \Gamma\left(\frac{N}{2}\right)^{-1}\quad{\rm as}\ \ t\to\left(\frac N2\right)^-.$$
%%%%
Given an arbitrary $\epsilon\in(0,1)$ and $x\in\R^N$, there exists $\varrho\in(0,\frac14)$ such that 
$$0<\int_{ \R^N\setminus B_{\frac1\varrho(0)})  }   f(y) \dd y\leq \epsilon\quad{\rm and}\quad \int_{   B_{ \varrho(0)})  }  |x-y|^{2t-N} \dd y\leq \epsilon,$$
and for fixed $\varrho$, there is $s_0\in(\frac{N}{4}, \frac{N}{2})$ such that  for $t\in[s_0,\frac{N}{2})$
$$\big||x-y|^{2t-N}-1\big|\leq \epsilon, $$ 
uniformly in  $\bar B_{\frac1\varrho(0)} \setminus  B_\varrho(x)$. Therefore, we have that 
\begin{eqnarray*}
 && \big|\tfrac{N-2t}{2} \big(\cP_{\ln}(t,\cdot)\ast f\big) (x) - \|f\|_{L^1(\R^N)}\big| 
\\[1mm] &&\qquad\qquad \leq   \tfrac{N-2t}{2} \cP_0(t)\Big( \int_{B_{\frac1\varrho(0)} \setminus  B_\varrho(x) } \big||x-y|^{2t-N}-1 \big| f(y) \dd y 
\\[0mm]
&&\qquad\qquad\qquad  + \int_{ \R^N\setminus B_{\frac1\varrho(0)}    }\big||x-y|^{2t-N}-1 \big| f(y) \dd y+  \int_{  B_{ \varrho }(x)    }\big||x-y|^{2t-N}-1\big| f(y) \dd y \Big)
\\[0mm]
&&\qquad\qquad \leq  c\Big(  \int_{B_{\frac1\varrho(0)} \setminus  B_\varrho(x) }  f(y) \dd y\epsilon
+\int_{ \R^N\setminus B_{\frac1\varrho(0)}    }  f(y) \dd y+\|f\|_{L^\infty(\R^N)} \int_{  B_{ \varrho }(x)    }|x-y|^{2t-N} dy\Big)
\\[0mm] 
&&\qquad\qquad \leq c(\|f\|_{L^1(\R^N)}+ \|f\|_{L^\infty(\R^N)}+1)\epsilon,
\end{eqnarray*}
which implies  (\ref{e1-1}) since $\epsilon$ is arbitrary.\smallskip

The decays at infinity in $(i)$ and $(ii)$ follows  by Proposition \ref{pr 2.1}.
  \hfill$\Box$\medskip

\begin{corollary}\label{cr 4.1}
Let  $0<T\leq \frac{N}{2}$,  $  f\in \cC(\R^N) \cap L^1(\R^N,(1+|x|)^{2T-N} dx)$
and 
$u(t,x)\geqq 0$ be a strong solution of  (\ref{eq 4.2}).
Then:
    
\nind (i)  the   strong solution of (\ref{eq 1.1})  is a weak solution;\smallskip
 
\nind (ii)  for $t_1,t_2>0$, $t_1+t_2<T$, 
\begin{equation}\label{e 4.4-}
  \cP_{\ln}(t_1,\cdot)\ast \big(\cP_{\ln}(t_2,\cdot)\ast f\big)=   \cP_{\ln}(t_1+t_1,\cdot) \ast f.
\end{equation} 
\end{corollary}
\noindent\proof The main point in the proof of Theorem \ref{thm1} is to show that
  the   strong solution of (\ref{eq 1.1})  is a weak solution and 
  $\cP_{\ln}(t,\cdot)\ast f$ is the unique solution of (\ref{eq 1.1}). By the uniqueness, Corollary \ref{cr 4.1-1} could be improved in the sense that  
$$
 \cP_{\ln}(t, \cdot) \ast u(\tau, \cdot)   =  u(t+\tau,\cdot) 
$$
 and
 $$
\cP_{\ln}(t_1,\cdot)\ast \big(\cP_{\ln}(t_2,\cdot)\ast f\big)=   \cP_{\ln}(t_1+t_1,\cdot) \ast f
$$
for $t_1+t_2<T$. 
We complete the proof.  \hfill$\Box$

\subsection{Existence of weak solutions}

\nind {\it Proof of Theorem \ref{th 1.2}. } Let $f$ be nonnegative function in $L^1(\R^N)$; there exists a sequence of nonnegative  functions  $\{f_n\}_{n\in\N}$  in $\cC_c(\R^N) \cap L^1(\R^N)$ such that 
$$ 0\leq f_n\leq f_{n+1}\leq\cdots \leq f\quad{\rm in}\ \R^N $$
and 
$$\lim_{n\to+\infty}\int_{\R^N}\big(f(x)-f_n(x)\big)dx=0. $$
The Cauchy  problem 
\begin{equation}\label{eq 5.1}
\left\{ \arraycolsep=1pt
\begin{array}{lll}
\displaystyle \partial_tu+ \lnlap u=0\quad \  &{\rm in}\ \,   (0,T)\times \R^N,\\[2mm]
 \phantom{ \lnlap \ \   }
\displaystyle   u(0,\cdot)=f_n\quad \ &{\rm{in}}\  \   \R^N,
\end{array}
\right.
\end{equation}
has a unique strong solution 
$$u_n(t,x)=\cP_0(t) \int_{\R^N} |x-y|^{2t-N} f_n(y)dy,$$
and 
the mapping $\N\ni n\to u_n$ is nondecreasing and 
$$\cP_0(t) \int_{\R^N} |x-y|^{2t-N} f_n(y)dy\leq \cP_0(t) \int_{\R^N} |x-y|^{2t-N} f(y)dy:=u_f(t,x),$$
where $u_f$ is well defined since $f\in L^1(\R^N)$. From Corollary \ref{cr 4.1},
$u_n$ is a weak solution of (\ref{eq 5.1}) and then
\begin{equation}\label{def w1-n}
 \int_{0}^{T'} \int_{\mathbb{R}^N}\left[-u_{n}  \partial_t \varphi_{t} +u_n  \lnlap
\varphi \right]   \dd t\dd x  =\int_{\mathbb{R}^N} f(x) \varphi(0, x) \dd x-\int_{\mathbb{R}^N} u_n\left(T', x \right) \varphi\left(T',x \right) \dd x
\end{equation}
for every   function $\varphi \in \cC^1\left( [0, T)\times \mathbb{R}^N\right)$ 
 such that $\varphi(t,\cdot)$ has compact support in $\R^N$ and $t\in [0,T)$.

Note that   $u_f$ is the limit of $\{u_n\}$, i.e. 
$$\cP_0(t) \int_{\R^N} |x-y|^{2t-N} f (y)dy=\cP_0(t)  \lim_{n\to+\infty}\int_{\R^N} |x-y|^{2t-N} f_n(y)dy.$$\smallskip

\noindent {\it Step 1: we prove that  $u_f \in L^1\big((0, T'),\,  L^{1} (\mathbb{R}^N,(1+|x|)^{2(T-T')-N}dx\big)$ for any $T'\in (0,T)$. }
Observe that  $\cP_0(t)\leq c t$ for $t\in(0, T']$ and then
\begin{eqnarray*}
  t\int_{\R^N}u_f(t,x)  (1+|x|)^{2(T-T')-N} \dd x &\leq & 
c t \int_{\R^N}\int_{\R^N}  |x-y|^{2t-N} f(y)(1+|x|)^{2(T-T')-N} \dd y\dd x  
\\[2mm]&=& t\int_{\R^N} \Big(\int_{\R^N}  |x-y|^{2t-N}(1+|x|)^{2(T-T')-N}dx\Big)f(y) dy
\\[2mm]&\leq & c\int_{\R^N} (1+|y|)^{2t-2T'+2T-N} f(y) dy
\\[2mm]&\leq & c\int_{\R^N}  (1+|x|)^{2T-N}   f(y) dy,
\end{eqnarray*}
where we used the fact, by  Lemma \ref{lm a-1},
\begin{eqnarray*}
\int_{\R^N}  |x-y|^{2t-N}(1+|x|)^{2(T'-T)}dx\leq \frac{c}{t} (1+|y|)^{2t+2(T'-T)}. 
\end{eqnarray*}

\smallskip

\noindent {\it Step 2: we prove that $u_f \in \cC\left((0, T),\,  L_{  {\rm loc }}^{1} (\mathbb{R}^N )\right)$. }\smallskip

 Fixed $t_0\in(0,T)$ and $ \varrho_0>\!\!>1$, for an arbitrary $\epsilon\in(0,1)$,   there exist  $\varrho\in(0,\frac14)$ and $s_0\in\big(0,\frac{\min\{t_0,T-t_0\}}{4}\big)$  such that for $|t-t_0|<s_0$
$$\big|\cP_0(t)-\cP_0(t_0)\big|< \epsilon,  \qquad  \int_{   B_{ \varrho}(0)  }   f(y)  \dd y\leq \epsilon,$$
and for $x\in \bar B_{\varrho_0} (0)$ 
$$  \big||x-y|^{2t-N}-|x-y|^{2t_0-N}\big|\leq \epsilon \ \ {\rm uniformly\; for  }\ \ y\in \R^N \setminus  B_\varrho(x). $$ 
Then  
\begin{eqnarray*}
 &&\int_{B_{ \varrho_0} (0)} |u_f(t,x)-u_f(t_0,x)| dx
\\[2mm]&=&\Big|\cP_0(t)  \int_{B_{ \varrho_0} (0)}   \int_{\R^N} |x-y|^{2t-N} f(y) \mathrm{d} y\dd x- \cP_0(t_0)  \int_{B_{ \varrho_0} (0)}   \int_{\R^N}  f(y) |x-y|^{2t_0-N}  \mathrm{d} y\dd x\Big| 
\\[2mm]& \leq & \int_{B_{ \varrho_0} (0) } \int_{\R^N} \big|\cP_0(t)|x-y|^{2t-N}-\cP_0(t_0)|x-y|^{2t_0-N}\big| f(y) \mathrm{d} y\dd x  
\\[2mm]
& \leqq& \big|\cP_0(t) -\cP_0(t_0)\big| \int_{B_{ \varrho_0} (0) } \int_{\R^N} f(y)  |x-y|^{2t-N} \mathrm{d} y\dd x
\\&&\qquad +\cP_0(t_0)\int_{B_{ \varrho_0} (0) } \int_{\R^N\setminus B_\varrho(x)} \big||x-y|^{2t-N}-|x-y|^{2t_0-N}\big| f(y)  \,\mathrm{d} y\dd x
\\&&\qquad +\cP_0(t_0)\int_{B_{ \varrho_0} (0) } \int_{  B_\varrho(x)} \big(|x-y|^{2t-N}+|x-y|^{2t_0-N}\big) f(y)  \,\mathrm{d} y\dd x
\\ 
& \leq &  c \epsilon +\cP_0(t_0)\Big(  \int_{B_{ \varrho_0} (0) }  f(y) \dd y\, \epsilon+ \epsilon \varrho_0^N \int_{\R^N} f(y)\dd y
  + c(\varrho_0^{-2t}+\varrho_0^{-2t_0}) \int_{  B_{ \varrho }(x)  }f(y) dy \Big)
\\[2mm] & \leq & c\epsilon,  
\end{eqnarray*}
which implies that
 $u_f \in \cC\left((0, T),\,  L^{1}_{loc} (\mathbb{R}^N)\right)$.

Now we can pass to the limit in (\ref{def w1-n}) as $n\to+\infty$ and we obtain that 
$$
\int_{0}^{T'} \int_{\mathbb{R}^N}\left[-u_f  \partial_t \varphi_{t} +u_f  \lnlap
\varphi \right]   \dd t\dd x  =\int_{\mathbb{R}^N} f(x) \varphi(0, x) \dd x-\int_{\mathbb{R}^N} u_f (T', x  ) \varphi (T',x  ) \dd x.
$$
As a conclusion, we get $\cP_{\ln}(t,\cdot)\ast f$ is the unique weak solution by Theorem \ref{thm2}.
\hfill$\Box$ \medskip

\noindent{\it Proof of Corollary \ref{th 1.3}. } Given $T\in (0,\frac{N}{2})$, let $t_n\in (0, \frac{N-2T}{4})$ with $n\in\N$ such that
$$t_n\to0\quad {\rm as}\ n\to+\infty,$$ then 
$u_n:=\cP_{\ln}(t_n+\cdot,\cdot)$ is the strong and the weak solution of
the Cauchy  problem 
\begin{equation}\label{eq 5.1-0}
\left\{ \arraycolsep=1pt
\begin{array}{lll}
\displaystyle \partial_tu+ \lnlap u=0\quad \  &{\rm in}\ \   (0,T)\times \R^N\\[2mm]
 \phantom{ \lnlap \ \   }
\displaystyle   u(0,\cdot)=\cP_{\ln}(t_n,\cdot) \quad \ &{\rm{in}}\  \   \R^N,
\end{array}
\right.
\end{equation}
therefore 
\begin{equation}\label{eq 5.1-00}
\int_{0}^{T'} \int_{\mathbb{R}^N}\left[-u_n  \partial_t \varphi_{t} +u_n  \lnlap
\varphi \right]   \dd t\dd x  =\int_{\mathbb{R}^N}\cP_{\ln}(t_n,x)  \varphi(0, x) \dd x-\int_{\mathbb{R}^N} u_n\left(T', x \right) \varphi\left(T',x \right) \dd x.
\end{equation}
We note the the limit of $\{\cP_{\ln}(t_n+\cdot,\cdot)\}_n$ is $\cP_{\ln}$ in related space
and  $\cP_{\ln}(t_n,\cdot)\to\delta_0$, then passing to the limit of (\ref{eq 5.1-00}), we obtain
that $\cP_{\ln}$ is a weak solution of (\ref{eq 1.1-d_0}).  
\hfill$\Box$ \medskip

%%%%%%%%%%%%%%%%%%%%%%%%%%%%%%%%%%%%%%%%%%%%%%%%%%%%%%%%%%%%%%%%%%%%%%%%%%%%%%%%%%%%%%%%%%%%%%%%%%%%%%%%%%%%%%%%%%%%%%%%%%%%%%%%%%%%%%%%%%%%%%%%%%%%%%%%%%%%%%%%%%%%%%%%%%%%%%%%%%%%%%%%%%%%%%%%%%%%%%%%%%%%%%%%%%%%%%%%%%%%%%%%%%%%%%%%%%%%%%%%%%%%%%%%%%%%%%%%%%%%%%%%%%%%%%%%%%%%%%%%%%%%

 \subsection{The initial trace }
 Let $\Gw$ be a bounded Lipschitz domain in $\BBR^N$. Following \cite{CT} we define 
 \begin{equation}\label{tr*0}
 \BBH(\Gw)=\left\{u\in L^1_{loc}(\BBR^N): u{\bf 1}_{\Gw^c}=0,\int\int_{|x-y|\leq 1}\frac{(u(x)-u(y))^2}{|x-y|^N}dxdy<\infty\right\},
 \end{equation}
 and $\CC_{c,D}(\Gw)$ the space of uniformly Dini continuous functions with support in $\overline \Gw$. Then $\BBH(\Gw)$ contains $\CC_{c,D}(\Gw)$ and it is an Hilbert space for the norm
  \begin{equation}\label{tr*01}
 \norm u_{\BBH(\Gw)}=\left(\int\int_{|x-y|\leq 1}\frac{(u(x)-u(y))^2}{|x-y|^N}dxdy\right)^{\frac 12},
 \end{equation}
 and the associated inner product
   \begin{equation}\label{tr*02}
\langle u,v\rangle_{\BBH(\Gw)}=\int\int_{|x-y|\leq 1}\frac{(u(x)-u(y))(v(x)-v(y)}{|x-y|^N}dxdy.
 \end{equation}
 %%%%%%%%%%%%%%%%%%CT%%THEOREM%%%%%%%%%%%%%%%%%%%%%%%%%%%%%%%%%%%%%%%%%%%%%%%%%%%%%%%%%%%%%%%%%%%%%%%%%%%%%%%%%%%%%%%%%%%%%%%%%%%%%%%%%%%%%%%%%%%%%%%%%%%%%%%%%%%%%%%%%%%%%%%%%%%%%%%%%%%%%%%%%%%%%%%%%%%%%%%%%%%%%%%%%%%
 The following fundamental result is proved in \cite{CT}
 \begin{theorem}\label{CT1}The imbedding of $\BBH(\Gw)$ into $L^2(\Gw)$ is compact and the operator $\CL_\Gd$ admits an increasing sequence of eigenvalues $\{\gl_k^{^L}(\Omega)\}$ tending to infinity with corresponding eigenfunctions $\gf_{_{k,\Gw}}$, $k\in\BBN^*$. Furthermore $\lambda_1^{^L}(\Omega)$ is simple and $\gf_{_{1,\Gw}}$ has constant sign. Finally, $\CL_\Gd$ satisfies the maximum principle if and only if $\gl_1^{^L}(\Omega)>0$, which in turn is satisfied if one of the two conditions (i) or (ii) of Lemma \ref{sec:main-result-max-principle} is satisfied. 
\end{theorem}
%%%%%%%%%%%%%%%%TRACE%%THEOREM%%%%%%%%%%%%%%%%%%%%%%%%%%%%%%%%%%%%%%%%%%%%%%%%%%%%%%%%%%%%%%%%%%%%%%%%%%%%%%%%%%%%%%%%%%%%%%%%%%%%%%%%%%%%%%%%%%%%%%%%%%%%%%%%%%%%%%%%%%%%%%%%%%%%%%%%%%%%%%%%%%%%%%%%%%%%%%%%%%%%%%%%%%%

\begin{lemma}\label{approx} Let $\gr\in \cC(\BBR^N)$ be a positive uniformly continuous function such that for any $\gth>0$ there exists $m_\gth>0$ such that 
\begin{equation}\label{C1}
\sup_{|x-z|\leq \gth}\frac{\gr(x)}{\gr(z)}\leq m_\gth.
\end{equation}
If $\{\gs_\gd\}_{\gd\in(0,1]}$ is a sequence of smooth mollifiers with support in $B_\gd$ and $u\in L_\gr^1(\BBR^N)$ where
$$L_\gr^1(\BBR^N)=\left\{\varphi\in L^1_{loc}(\BBR^N): \int_{\BBR^N}|\varphi|\gr dx<\infty\right\},
$$
then $u\ast\gs_\gd\to u$ in $L_\gr^1(\BBR^N)$.
\end{lemma}
\proof Let $\ge>0$ and $u\in L^1_\gr(\BBR^N)$. Then there exists $R>0$ such that 
$\norm {u{\bf 1}_{B^c_R}}_{L^1_\gr}<\frac\ge4$. There exists $u_n\in C_c(B_R)$ such that $\norm {u{\bf 1}_{B_R}-u_n}_{L^1_\gr}<\frac\ge4$. Then
$$u-u\ast\gs_\gd=u-u_n+u_n-u_n\ast\gs_\gd+u_n\ast\gs_\gd-u\ast\gs_\gd.
$$
There holds
$$\BA{lll}\norm{u-u\ast\gs_\gd}_{L^1_\gr}\leq \norm{u-u_n}_{L^1_\gr}+\norm{u_n-u_n\ast\gs_\gd}_{L^1_\gr}+\norm{u_n\ast\gs_\gd-u\ast\gs_\gd}_{L^1_\gr}\\[3mm]\phantom{\norm{u-u\ast\gs_\gd}_{L^1_\gr}}\leq I+II+III.
\EA$$
Now $I\leq \frac\ge 2$. 
$$\BA{lll}\dsps II= \int_{B_{R+\gd}}\left|\int_{B_\gd}u_n(x)\gs_\gd(y)dy-\int_{B_\gd}u_n(x-y)\gs_\gd(y)dy\right| \gr(x)dx
\\[4mm]\phantom{II}\dsps
\leq \int_{B_{R+\gd}}\int_{B_\gd}\left|u_n(x)-u_n(x-y)\right|\gs_\gd(y)\gr(x)dydx
\\[4mm]\dsps\phantom{II}\leq \sup_{|x|\leq R+\gd}\gr(x)\sup_{|y|\leq \gd}|u_n(x)-u_n(x-y)|.
\EA$$
For $n$ fixed, there exists $\gd_0$ such that for $\gd\leq\gd_0$, one has $II\leq\frac\ge 4$. Finally, 
$$\BA{lll}\dsps III=\int_{\BBR^N}|(u_n-u)\ast\gs_\gd(x)|\gr(x)dx\\[4mm]
\phantom{lll}\dsps
\leq \int_{B_{R+\gd}}\int_{B_\gd}|(u_n(x-y)-u(x-y)|\gs_\gd(y)dy\gr(x)dx+\int_{B_{R+\gd}}\int_{B_\gd}|u(x-y)|\gs_\gd(y)dy\gr(x)dx
\\[4mm]
\phantom{lll}\dsps
\leq \sup_{|x|\leq R+\gd}\gr(x)\norm{u_n-u}_{L^1(B_{R+\gd})}+m_\gd\int_{B_{R+\gd}}\int_{B_\gd}|u(x-y)|\gs_\gd(y)dy\gr(x-y)dx
\\[4mm]
\phantom{lll}\dsps
\leq \left(\sup_{|x|\leq R+\gd}\gr(x)+m_\gd\right)\frac{\ge}{4}.
\EA$$
Combining the three estimates we complete the proof.\hfill$\Box$ 
%%%%%%%%%%%%%%%%TRACE%%THEOREM%%%%%%%%%%%%%%%%%%%%%%%%%%%%%%%%%%%%%%%%%%%%%%%%%%%%%%%%%%%%%%%%%%%%%%%%%%%%%%%%%%%%%%%%%%%%%%%%%%%%%%%%%%%%%%%%%%%%%%%%%%%%%%%%%%%%%%%%%%%%%%%%%%%%%%%%%%%%%%%%%%%%%%%%%%%%%%%%%%%%%%%%%%%

 \begin{theorem}\label{tr1}
Let $u$ be a nonnegative strong solution of (\ref{e 1-v}) in $(0,T)\ti\BBR^N$. Then 
\begin{equation}\label{tr*01}\dsps
u\in L^1(0,T';L^1(\BBR^N,(1+|x|)^{-N}),
\end{equation}
for any $T'<T$. There exists a nonnegative Radon measure $\gm$ in $\BBR^N$ satisfying 
(\ref{R1}) such that 
\begin{equation}\label{tr*1}\dsps
\lim_{t\to 0}\int_{\BBR^N}u(t,x)\gz(x)dx=\int_{\BBR^N}\gz(x)d\gm(x),
\end{equation}
and there holds
\begin{equation}\label{tr*2}\dsps
u(t,x)\geq \int_{\BBR^N}\CP_{\ln}(t,x-y) d\gm(y)\quad\text{for all }\;(t,x)\in (0,T)\ti\BBR^N.
\end{equation}
\end{theorem}
%%%%%%%%%%%%%%%PROOF%%%%%%%%%%%%%%%%%%%%%%%%%%%%%%%%%%%%%%%%%%%%%%%%%%
\proof {\it Step 1: Let $u\in Dom_{\CL_{\Gd}}(\BBR^N)$, then given bounded Lipschitz domain $\Omega$, the following identity holds}
\begin{equation}\label{tr*3}
\int_{\Gw}\gf_{_{1,\Gw}}\lnlap u dx=\gl_1^{^L}(\Gw)\int_{\Gw}u\gf_{_{1,\Gw}} dx -c_{_N}\int_{\Gw}\gf_{_{1,\Gw}}(x)\int_{\Gw^c}\frac{u(y)}{|x-y|^N}dydx.
\end{equation}
Indeed, we have 
$$\begin{array}{lll}\dsps
\int_{\Gw}\gf_{_{1,\Gw}}\lnlap u dx=\int_{\Gw}\gf_{_{1,\Gw}}(x)\left(c_{_N}\int_{\Gw}\frac{u(x)-u(y)}{|x-y|^N}dy-c_{_N}\int_{\Gw^c}\frac{u(y)}{|x-y|^N}dy
+\left(h_\Gw(x)+\gr_N\right)u(x)\right) dx\\[4mm]
\phantom{\dsps
\int_{\Gw}\gf_{_{1,\Gw}}\lnlap u dx}\dsps
=\int_{\Gw}u(x)\left(c_{_N}\int_{\Gw}\frac{\gf_{_{1,\Gw}}(x)-\gf_{_{1,\Gw}}(y)}{|x-y|^N}dy+\left(h_\Gw(x)+\gr_N\right)\gf_{_{1,\Gw}}(x)\right)dx\\[4mm]
\phantom{\dsps
\int_{\Gw}\gf_{_{1,\Gw}}\lnlap u dx-----------------}\dsps
-c_{_N}\int_{\Gw}\gf_{_{1,\Gw}}\int_{\Gw^c}\frac{u(y)}{|x-y|^N}dydx
\\[4mm]
\phantom{\dsps
\int_{\Gw}\gf_{_{1,\Gw}}\lnlap u dx}\dsps
=\gl_1^{^L}(\Gw)\int_{\Gw}u\gf_{_{1,\Gw}} dx -c_{_N}\int_{\Gw}\gf_{_{1,\Gw}}(x)\int_{\Gw^c}\frac{u(y)}{|x-y|^N}dydx.
\end{array}$$

\nind {\it Step 2: Construction of the initial trace.}
We multiply the equation by $\gf_{_{1,\Gw}}$ and get
\begin{equation}\label{tr*4}
\frac{d}{dt}\int_{\Gw}u(t,x)\gf_{_{1,\Gw}}(x)dx+\gl_1^{^L}(\Gw)\int_{\Gw}u(t,x)\gf_{_{1,\Gw}}(x)dx=c_{_N}\int_{\Gw}\gf_{_{1,\Gw}}(x)\int_{\Gw^c}\frac{u(y)}{|x-y|^N}dydx\geq 0.
\end{equation}
Hence 
$$t\mapsto e^{-\gl_1^{^L}(\Gw)t}\int_{\Gw}u(t,x)\gf_{_{1,\Gw}}(x)dx
$$
is an increasing nonnegative function, which therefore admits a nonnegative limit $\ell(\Gw)$ when $t\to 0$. This implies in particular that for any 
compact set $G\subset\Gw$ and $T'\in (0,T)$, 
\begin{equation}\label{tr*5}
\int_0^{T'}\int_Gu(t,x)dxdt<\infty.
\end{equation}
Let $\gz\in C^\infty_c(\Gw)$ which vanishes outside the open set $G\subset\overline G\subset\Gw$, $\gz\geq 0$, then 
$$\lnlap\gz (x)=c_{_N}\int_{G}\frac{\gz(x)-\gz(y)}{|x-y|^N}dy+\left(h_G(x)+\gr_N\right)\gz(x),
$$
thus
\begin{equation}\label{tr*51}\int_{G}\gz(x)\lnlap u(t,x) dx=\int_{G}u(t,x)\lnlap\gz (x) dx -c_{_N}\int_{G}\gz(x)\int_{G^c}\frac{u(t,y)}{|x-y|^N}dydx.
\end{equation}
It follows that
\begin{equation}\label{tr*6}
\frac{d}{dt}\left(\int_{G}u(t,x)\gz(x)dx-\norm{\lnlap\gz}_{L^{\infty}(G)}\int_t^{T'}\int_G u(s,x)dxds\right)
\geq c_{_N}\int_{G}\gz(x)\int_{G^c}\frac{u(t,y)}{|x-y|^N}dydx\geq 0.
\end{equation}
Hence $t\mapsto \int_{G}u(t,x)\gz(x)dx$ admits a nonnegative limit $\ell(\gz)$ when $t\to 0^+$, and 
 \begin{equation}\label{tr*7}
0\leq \ell(\gz)\leq \norm{\lnlap\gz}_{L^{\infty}(G)}\int_0^{T'}\int_G u(t,x)dxdt+\int_{G}u(T',x)\gz(x)dx.
\end{equation}
Now $\norm{\lnlap\gz}_{L^{\infty}(G)}\leq \norm{\gz}_{C^2(G)}$, therefore 
the mapping $\gz\mapsto \ell(\gz)$ is a positive distribution, hence a measure. 
By a partition of unity, we construct the nonnegative Radon measure $\gm$ such that 
 \begin{equation}\label{tr*8}
\lim_{t\to 0}\int_{\BBR^N}u(t,x)\gz(x) dx=\int_{\BBR^N}\gz d\gm(x)\quad\text{for all }\,\gz\in C_c(\BBR^N).
\end{equation}
This measure is uniquely defined and is called the {\it initial trace }of $u$. Finally there holds from (\ref{tr*51}),
 \begin{equation}\label{tr**8}\int_Gu(T',x)\gz(x)dx-\int_{\BBR^N}\gz d\gm(x)\\[4mm]
= -\int_0^{T'}\int_{G}u(t,x)\lnlap\gz (x) dx+c_{_N}\int_0^{T'}\int_{G}\gz(x)\int_{G^c}\frac{u(t,y)}{|x-y|^N}dydxdt.
\end{equation}
\smallskip

\nind{\it Step 3: the relation (\ref{tr*01}) holds.} We derive from (\ref{tr*6}) and (\ref{tr*8})
$$\int_{G}u(T',x)\gz(x)dx\geq \norm{\lnlap\gz}_{L^{\infty}(G)}\int_0^{T'}\!\int_G u(s,x)dxds+\int_{\BBR^N}\gz d\gm(x)+c_{_N}\int_0^{T'}\!\int_{G}\gz(x)\int_{G^c}\frac{u(t,y)}{|x-y|^N}dydxdt.
$$
Since (\ref{tr*5}) holds we obtain that 
$$% \begin{equation}\label{tr*8-1}
\int_0^{T'}\int_{G}\gz(x)\int_{G^c}\frac{u(t,y)}{|x-y|^N}dydx<+\infty.
$$%\end{equation}
We take $G=B_{2}$ and $\gz=1$ on $B_1$, then 
$$\int_0^{T'}\int_{B_1}\int_{B_2^c}\frac{u(t,y)}{|x-y|^N}dydxdt<+\infty.
$$
If $x\in B_1$ and $y\in B_2^c$, $|x-y|\geq \frac{1}{3}(|y|+1)$, hence 
$$\int_0^{T'}\int_{B_2^c}\frac{u(t,y)}{(|y|+1)^N}dydt<\infty.
$$
Combined with (\ref{tr*5}), this inequality implies (\ref{tr*01}).
\smallskip

\nind{\it Step 4: the inequality (\ref{tr*2}) is verified.} Let  $\{\gs_\gd\}_{\gd>0}$ be a sequence of smooth mollifiers with support in $B_\gd$, $\gm_{\gd}=\gm\ast\gs_\gd$ and $u_\gd=u\ast\gs_\gd$. Then 
 $$\BA{lll}\dsps 
 \int_{G}u_\gd(T',x)\gz(x)dx=\int_0^{T'}\int_G \lnlap\gz (x)u_\gd(t,x)dxdt+\int_{\BBR^N}\gz(x) u_\gd(0,x) dx\\[4mm]
 \phantom{------------------}
 \dsps+c_{_N}\int_0^{T'}\int_{G}\gz(x)\int_{G^c}\frac{u_\gd(t,y)}{|x-y|^N}dydxdt.
\EA $$
When $\gd\to 0$, $u_\gd$ converges to $u$ uniformly on compact subsets of $(0,T]\ti\BBR^N$ and in $L^1((0,T')\ti B_R)$ for any $R>0$. 
As in Step 3 we have that $u_\gd\in L^1(0,T';L^1(\BBR^N,(1+|x|)^{-N})$, and by  Lemma \ref{approx} it converges to $u$ in $L^1(0,T';L^1(\BBR^N,(1+|x|)^{-N})$.
Consequently, it follows from (\ref{tr**8}) that the following limit exists with explicit value, 
$$\lim_{\gd\to 0}\int_{\BBR^N}\gz(x) u_\gd(0,x) dx=\int_{\BBR^N}\gz(x) d\gm(x).
$$
Furthermore,
$$\BA{lll}\dsps\int_{\BBR^N}u_\gd(t,x)\gz(x)dx=\int_{\BBR^N}\left(\int_{\BBR^N}u(t,y)\gs_\gd(x-y)dy\right)\gz(x)dx\\[4mm]
\phantom{\dsps\int_{\BBR^N}u_\gd(t,x)\gz(x)dx}=\dsps\int_{\BBR^N}u(t,y)\left(\int_{\BBR^N}\gs_\gd(x-y)\gz(x)dx\right)dy
\\[4mm]
\phantom{\dsps\int_{\BBR^N}u_\gd(t,x)\gz(x)dx}=\dsps\int_{\BBR^N}u(t,y)(\gz\ast\hat\gs_\gd)(y)dy
\EA$$
where $\hat\gs_\gd(y)=\gs_\gd(-y)$. Therefore
$$\BA{lll}\dsps\int_{\BBR^N}u_\gd(0,x)\gz(x)dx=
\lim_{t\to 0}\int_{\BBR^N}u_\gd(t,x)\gz(x)dx=\int_{\BBR^N}(\gz\ast\hat\gs_\gd)(y)d\gm(y)=
\int_{\BBR^N}\gz(y)\gm_\gd(y)dy,\\[4mm]
\EA$$
where $\gm_\gd=\gm\ast\gs_\gd$. This implies that $u_\gd$ is a strong solution with initial data $\gm_\gd$.
\medskip
We apply Lemma \ref{lm 4.1-im} and get
$$% \begin{equation}\label{tr*8-2}
u_\gd\geq u_{\gm\ast\gs_\gd},
$$%\end{equation}
which yields to (\ref{tr*2}) by letting $\gd\to 0$.\hfill$\Box$ 
 
 \medskip
%%%%%%%%%%%%%%%%%%%%%%%%%%%%%%%%%%%%%%%%%%%%%%%%%%%%%%%%%%%%%%%%%%%%%%%%%%%%%%%%%%%%%%%%%%%%%%%%%%%%%%%%%%%%%%%%%%%%%%%%%%%%%%%%%%%%%%%%%%%%%%
%%%%%%%%%%%%%%%%%%%%%%%%%%%%%%%%%%%%%%%%%%%%%%%%%%%%%%%%%%%%%%%%%%%%%%%%%%%%%%%%%%%%%%%%%%%%%%%%%STATIONNARY%%%%%%%%%%%%%%%%%%%%%%%%%%%%%%%%%
%%%%%%%%%%%%%%%%%%%%%%%%%%%%%%%%%%%%%%%%%%%%%%%%%%%%%%%%%%%%%%%%%%%%%%%%%%%%%%%%%%%%%%%%%%%%%%%%%%%%%%%%%%%%%%%%%%%%%%%%%%%%%%%%%%%%%%%%%%%%%%
%%%%%%%%%%%%%%%%%%%%%%%%%%%%%%%%%%%%%%%%%%%%%%%%%%%%%%%%%%%%%%%%%%%%%%%%%%%%%%%%%%%%%%%%%%%%%%%%%%%%%%%%%%%%%%%%%%%%%%%%%%%%%%%%%%%%%%%%%%%%%%

 \section{Fundamental solutions of $\lnlap$ }
We first give the proof of the existence and asymptotic properties of the solutions of the Helmholtz equation in $\BBR^N$.
 \subsection{Existence in $\BBR^N$}
 \begin{lemma}\label{helm} If $N\geq 3$, the fundamental solution $\Phi_{1}$ to the equation 
  \begin{equation}\label {hel1}\begin{array}{lll}\dsps
-\Gd u-u=\gd_0\quad\text{in }\cD'(\BBR^N).
\end{array}\end{equation}
satisfies 
  \begin{equation}\label {hel2}\begin{array}{lll}\dsps
\Phi_{1}(x)=c_{N,1}|x|^{2-N}(1+o(1))\qquad\text{as }x\to 0,
\end{array}\end{equation}
and
 \begin{equation}\label{hel9} 
\Phi_1(x)=-\frac{\sin(|x|-\frac {(N-1)\gp}4)}{2(2\gp |x|)^{\frac{N-1}{2}}}(1+O(|x|^{-1}))\quad\text{as }|x|\to \infty
 \end{equation}

 \end{lemma}
 \noindent{\it Proof.} The equation reduces to the classical Helmoltz equation
   \begin{equation}\label {hel3}\begin{array}{lll}\dsps
-\Gd u-u=\gd_0\quad\text{in }\, \BBR^N.
\end{array}\end{equation}
The method for obtaining $\Phi_1$ is standard but for the sake of completeness, we recall some elements. 
The function $\Phi_1$ is radial, thus, if we define $w$ by  $\Phi_1(r)=r^{1-\frac N2}w(r)$, it satisfies the Bessel equation of order $\frac N2-1$
 \begin{equation}\label{eq H} 
w''(r)+\frac{w'(r)}{r}+\left(1-\frac{(N-2)^2}{4r^2}\right)w(r)=0\quad\text{on }\,(0,\infty).
 \end{equation}
The solution of (\ref{eq H}) is a linear combination of two Hankel functions respectively of the first kind $H^{(1)}_{\frac N2-1}$ and of the second kind  $H^{(2)}_{\frac N2-1}$, that is
 \begin{equation}\label{hel4} 
w(r)=AH^{(1)}_{\frac N2-1}(r)+BH^{(2)}_{\frac N2-1}(r).
 \end{equation}
 We find $A$ and $B$ if we use the relation 
 $$\dsps|S^{N-1}|\lim_{r\to 0}r^{N-1}\Phi_1'r)=-|S^{N-1}|,$$ 
 which yields with $B=0$, 
 $$\lim_{r\to 0}\left(\left(1-\frac N2\right)e^{-\frac N2}H^{(1)}_{\frac N2-1}(r)+r^{1-\frac N2}\left(H^{(1)}_{\frac N2-1}\right)'(r)\right)=
 -|S^{N-1}|.$$
 The Hankel functions have the following behaviours when $r\to 0$ (see e.g. \cite{WW}),
  \begin{equation}\label{hel5} 
\BA{lll}
\phantom{--}\dsps
 H^{(1)}_{\frac N2-1}(r)=-{\rm i}\frac{\Gg(\frac N2-1)}{\gp}\left(\frac{2}{r}\right)^{\frac N2-1}(1+o(1)),\\[4mm]
\!\!\dsps\left(H^{(1)}_{\frac N2-1}\right)'(r)={\rm i} \frac{(N-2)\Gg(\frac N2-1)}{4\gp}\left(\frac{2}{r}\right)^{\frac N2}(1+(1)).
\EA
 \end{equation}
 Putting altogether these estimates we find 
$$
A=\frac{{\rm i}}{4}\left(\frac 1{2\gp}\right)^{\frac N2-1},
$$
 and 
   \begin{equation}\label{hel6} 
\Phi_1(r)=\frac{{\rm i}}{4}\left(\frac 1{2\gp r}\right)^{\frac N2-1} H^{(1)}_{\frac N2-1}(r),
 \end{equation}
 which implies
    \begin{equation}\label{hel8} 
\Phi_1(r)=\frac{\Gg(\tfrac{N}{2}-1)}{4\gp^{\tfrac N2}}r^{2-N}(1+o(1))\quad\text{as }r\to 0.
 \end{equation}
 The Hankel function is oscillatory at infinity with asymptotic frequency $2\gp$. Expression (\ref{hel9}) can be found in \cite[see (5.16.1)-(5.16.3]{Lebed}.
Note in particular the following rough a priori estimate,
 \begin{equation}
  \label{hel10}
  |\Phi_1(x)| \le C_0 \max \{|x|^{2-N},|x|^{\frac{1-N}{2}}\} \quad
\end{equation}   
for $x \in \R^N \setminus \{0\}$.
%%%%%%%%%%%%%%%%%%%%%%%%%%%%%%%%%%%%%%%%%%%%%%%%%%%%%%%%%%%%%%%%%%%%%%%%%%%%%%%%%%%%%%%%%%%%%%%%%%%%%%%%%%%%%%%%%%%%%%%%%%%%%%%%%%%%%%%%%%%%%%%%%%%%%%%%%%%%%%%\\
      This ends the proof.\hfill$\Box$ \medskip
    
    %%%%%%%%%%%%%%%%%%%%%%%%%%%%%%%%%%%%%%%%%%%%%%%%%%%%%%%%%%%%%%%%%%%%%%%%%%%%%%%%%%%%%%%%%%%%%%%%%%%%%%%%%%%%%LEMMA%%%%%%%%%%%%%%%%%%%%%%%%%%%%%%%%%%%%%%%%%%%%%%%%%%%%%%%%%%%%%%%%%%%%%%%%%%%%%%%%%%%%%%%%%%%%%%%%%%%%%%%%%%%%%%%%%%%%%%%%%%%
 
The next estimate involving the Riesz kernel $I_{2\ga}(x):=|x|^{2\ga-N}$ will be used several times in the sequel.
 \begin{lemma}\label{est-l} If $N\geq 3$  and $\ga\in (0,1)$, then there holds for $x\in B_{\frac1e}\setminus\{0\}$ and some $c>0$,
 \begin{equation}\label{AS6}
 |I_{2\ga}\ast\Phi_1(x)|\leq \left\{ \BA{lll}c|x|^{2\ga+2-N}\qquad&\text{if }N\geq 4\text{ and }0<\ga<1,\text{ or }N=3\text{ and }0<\ga< \frac 12\\[2.5mm]
 c\ln\frac 1{|x|}\qquad&\text{if }N=3\text{ and }\ga=\frac 12\\[2.5mm]
 c\qquad&\text{if }N=3\text{ and }\frac 12<\ga<1.
 \EA\right.
 \end{equation}
  There also holds if $x\in B_1^c$,
     \begin{equation}\label{hel 2}
  |I_{2\ga}\ast \Phi_1|(x)\leq \left\{\BA{lll}c|x|^{2\ga-\frac{N+1}{2}}\qquad&\text{ if }\,\frac 12<\ga<1, \\[2mm]
  C_\ge|x|^{2\ga-\frac{N+1}{2}+\ge}\qquad&\text{ if }\, 0<\ga\leq \frac 12,
  \EA\right.
   \end{equation}
  \end{lemma}
  %%%%%%%%%%%%%%%%%%%%%%%%%%%%%%%%%%%%%%%%%%%%%%%%%%%%%%%%%%%%%%%%%%%%%%%%%%%%%%%%%%%%%%%%%%%%%%%%%%%%%%%%%%PROOF OF LEMMA%%%%%%%%%%%%%%%%%%%%%%%%%%%%%%%%%%%%%%%%%%%%%%%%%%%%%%%%%%%%%%%%%%%%%%%%%%%%%%%%%%%%%%%%%%%%%%%%%%%%%%%%%%%%%%%%%%%%%%%%%%%%%
  \noindent{\it Proof.} We have
     \begin{equation}\label{eq H6-1} \BA{lll}
 \dsps |I_{2\ga}\ast \Phi_1|(x)=c_{_N}\left|\int_{\R^N} |x-y|^{2\ga-N}\Phi_1(y)dy\right|\\[4mm]
\phantom{I_{2\ga}\ast |\Phi_1|(x)}\dsps\leq c\int_{B_1} |x-y|^{2\ga-N}|y|^{2-N}dy+c_1\left|\int_{B_1^c} |x-y|^{2\ga-N}\Phi_1(y)dy\right|
\\[4mm]
\phantom{I_{2\ga}\ast |\Phi_1|(x)}\dsps \leq  C_\ga(x)+D_\ga(x).
 \EA\end{equation}
%%%
{\it Step 1}. We have near $x=0$,
 $$C_\ga(x)=c|x|^{2\ga+2-N}\int_{B_{\frac 1{|x|}}}|z|^{2-N}|{\bf e}_N-z|^{2\ga-N}dz,
 $$
 %and
% $$I\!I_\ga(x)=c_2|x|^{2\ga+\frac{1-N}{2}}\int_{B^c_{\frac 1{|x|}}}|z|^{\frac{1-N}{2}}|{\bf e}_N-z|^{2\ga-N}dz,$$
 where ${\bf e}_N=(0,0,...1)$. If $N\geq 4$ and for any $\ga\in (0,1)$, or if $N=3$ and $0<\ga<\frac 12$ , we have, 
\begin{equation}\label{AS1}C_\ga(x)= c|x|^{2\ga+2-N}\left(\int_{\BBR^N}|z|^{2-N}|{\bf e}_N-z|^{2\ga-N}dz\right)(1+o(1))\,\text{ when }x\to 0.
\end{equation}
If $N=3$ and $\ga=\frac 12$, we have
\begin{equation}\label{AS2}C_\ga(x)= c\left(\ln\tfrac{1}{|x|}\right)(1+o(1))\,\text{ when }x\to 0.
\end{equation}
Finally, if $N=3$ and $\frac 12\leq \ga<1$, we have
 \begin{equation}\label{AS3}C_{\ga}(x)=c'(1+o(1))\quad\text{when }x\to 0.
\end{equation}
Concerning the second term $D_\ga$,  we use the oscillations of $\Phi_1(x)$ when $|x|\to\infty$. There holds
$$|x-y|^{2\ga-N}=|y|^{2\ga-N}+O(|y|^{2\ga-N-1})\quad\text{when }|y|\to \infty,
$$
and
$$\Phi_1(y)=c^*_N|y|^{\frac {1-N}{2}}\sin\left(|y|-\frac{(N-1)\gp}{4}\right)+O\left(|y|^{-\frac {1+N}{2}}\right)\quad\text{when }|y|\to\infty.
$$
Furthermore these two estimates hold uniformly with  respect to $x\in B_1$. Therefore, for  $x\in B_1$ small enough and $A>1$ large enough, we have
$$\BA{lll}\dsps
\left|\int_{B_{1}^c} |x-y|^{2\ga-N}\Phi_1(y)dy\right|\leq \left|\int_{B_{A}\setminus B_1} |x-y|^{2\ga-N}\Phi_1(y)dy\right|+\left|\int_{B_{A}\setminus B_1} |x-y|^{2\ga-N}\Phi_1(y)dy\right|\\[4mm]
\phantom{\dsps
\left|\int_{B_{1}^c} |x-y|^{2\ga-N}\Phi_1(y)dy\right|}\dsps
=c'+c\int_A^\infty r^{2\ga-\frac {N+1}{2}}\sin\left(r-\frac{(N-1)\gp}{4}\right)dr+O\left(\int_A^\infty r^{2\ga-\frac {N+3}{2}}dr\right).
\EA$$
By integration by parts, 
$$\BA{lll}\dsps\int_A^\infty r^{2\ga-\frac {N+1}{2}}\sin\left(r-\frac{(N-1)\gp}{4}\right)dr=-\left[r^{2\ga-\frac {N+1}{2}}\cos\left(r-\frac{(N-1)\gp}{4}\right)\right]_{r=A}^\infty\\[4mm]\dsps
\phantom{-----------------}+\left(2\ga-\frac {N+1}{2}\right)\int_A^\infty r^{2\ga-\frac {N+3}{2}}\cos\left(r-\frac{(N-1)\gp}{4}\right)dr.
\EA$$
This implies 
 \begin{equation}\label{AS5}D_\ga(x)=
 \left|\int_{B_{1}^c} |x-y|^{2\ga-N}\Phi_1(y)dy\right|\leq c'+c''A^{2\ga-\frac {N+1}{2}},
\end{equation}
where $c'=c'(A)>0$ and this term is bounded whenever $A$ is fixed. Combining (\ref{eq H6-1}) and (\ref{AS5}) we obtain estimate (\ref{AS6}) when $0<|x|\leq 1$
for some $c>0$. \smallskip
 
  \noindent{\it Step 2}. When $|x|$ is large, we write
  \begin{equation}\label{AS7} \BA{lll}
 \dsps |I_{2\ga}\ast \Phi_1|(x)\leq c_1\int_{B_1} |x-y|^{2\ga-N}|y|^{\frac{1-N}{2}}dy+c_1\left|\int_{B_1^c} |x-y|^{2\ga-N}\Phi_1(y)dy\right|
\\[4mm]
\phantom{I_{2\ga}\ast |\Phi_1|(x)}\dsps \leq  \widetilde C_\ga(x)+\widetilde D_\ga(x).
 \EA\end{equation}
 By a change of variable 
 $$\int_{B_1} |x-y|^{2\ga-N}|y|^{\frac{1-N}{2}}dy=|x|^{2\ga+\frac{1-N}2}\int_{B_\frac {1}{|x|}} |{\bf e}_N-z|^{2\ga-N}|z|^{\frac{1-N}{2}}dz.
 $$
If $|x|\geq 2$, then $|z|\leq \frac12$ thus $ |{\bf e}_N-z|\geq \frac 12$, hence 
$$\int_{B_\frac {1}{|x|}} |{\bf e}_N-z|^{2\ga-N}|z|^{\frac{1-N}{2}}dz\leq 2^{N-2\ga}\gw_{_{N}}\int_0^{\frac 1{|x|}}r^{\frac{N-1}{2}}dr=\frac{2^{N+1-2\ga}\gw_N}{N+1}
|x|^{-\frac{N+1}{2}}.$$
This yields
  \begin{equation}\label{hel17}
\widetilde C_\ga(x)\leq \frac{2^{N+1-2\ga}\gw_N}{N+1}|x|^{2\ga-N}.
 \end{equation}
Concerning $\widetilde D_\ga(x)$, we  have for $|y|\geq 1$, 
 $$\Phi_1(y)=c|y|^{\frac {1-N}{2}}\sin\left(|y|-\frac{(N-1)\gp}{4}\right)+\gb(y)|y|^{-\frac {1+N}{2}},
 $$
 where $\gb$ is a bounded function. Then 
     \begin{equation}\label{hel18}\BA{lll}
\dsps \widetilde D_\ga(x)\leq  c\left|\int_{B^c_1} |x-y|^{2\ga-N}|y|^{\frac {1-N}{2}}\sin\left(|y|-\frac{(N-1)\gp}{4}\right)\right| +\norm\gb_{L^\infty}\int_{B^c_1} |x-y|^{2\ga-N}|y|^{-\frac {1+N}{2}}dy.
\EA \end{equation}
We first notice that 
 \begin{equation}\label{hel18+}\BA{lll}\dsps\int_{B^c_1} |x-y|^{2\ga-N}|y|^{-\frac {1+N}{2}}dy=|x|^{2\ga-\frac{N+1}{2}}\int_{B^c_{\frac 1{|x|}}}|{\bf e}_N-z|^{2\ga-N}|z|^{-\frac {1+N}{2}}dz
\\[5mm]
\phantom{\dsps\int_{B^c_1} |x-y|^{2\ga-N}|y|^{-\frac {1+N}{2}}dy}\dsps
=|x|^{2\ga-\frac{N+1}{2}}\left(\int_{\BBR^N}|{\bf e}_N-z|^{2\ga-N}|z|^{-\frac {1+N}{2}}dz\right)(1+o(1)).
\EA\end{equation}
For the second term, its treatment depends on the value of $\ga$ with respect to $\frac 12$. \\
 - If $\frac 12<\ga<1$ we use spherical coordinates $y=(r,\gs)$ and $x=(\gr,\gth)$, then
$$\BA{lll}\dsps
\int_{B^c_{1}} |x-y|^{2\ga-N}|y|^{\frac {1-N}{2}}\sin\left(|y|-\frac{(N-1)\gp}{4}\right)dy \\[4mm]
\phantom{---------}\dsps=
\int_{1}^\infty\int_{S^{N-1}}r^{\frac{N-1}{2}}(r^2+\gr^2-2r\gr\langle\gs,\gth\rangle)^{\ga-\frac N2}\sin\left(r-\frac{(N-1)\gp}{4}\right)dSdr\\[4mm]
\phantom{---------}\dsps=
-\int_{S^{N-1}}\left[r^{\frac{N-1}{2}}(r^2+\gr^2-2r\gr\langle\gs,\gth\rangle)^{\ga-\frac N2}\cos\left(r-\frac{(N-1)\gp}{4}\right)\right]_{r=1}^\infty dS
\\[4mm]
\phantom{----------}\dsps+\int_{1}^\infty\int_{S^{N-1}}\frac{d}{dr}\left(r^{\frac{N-1}{2}}(r^2+\gr^2-2r\gr\langle\gs,\gth\rangle)^{\ga-\frac N2}\right)\cos\left(r-\frac{(N-1)\gp}{4}\right) dSdr\\[4mm]
\phantom{---------}\dsps
:=-J_1+J_2.
\EA$$
We have that 
$$|J_1|\leq c(1+|x|^2)^{\ga-\frac{N}{2}}.
$$
Furthermore, since
\begin{equation}\label{der1}\BA{lll}\dsps \frac{d}{dr}\left(r^{\frac{N-1}{2}}(r^2+\gr^2-2r\gr\langle\gs,\gth\rangle)^{\ga-\frac N2}\right)
=\frac{N-1}{2}r^{\frac{N-3}{2}}\left(r^2+\gr^2-2r\gr\langle\gs,\gth\rangle\right)^{\ga-\frac N2}\\[4mm]\phantom{--------------}
\dsps
+\left(2\ga-N\right)r^{\frac{N-1}{2}}(r-\gr\langle\gs,\gth\rangle)(r^2+\gr^2-2r\gr\langle\gs,\gth\rangle)^{\ga-\frac {N+2}2},
\EA\end{equation}
uniformly with respect to $\gr\leq 2r$ and $(\gs,\gth)\in S^{N-1}\ti S^{N-1}$, we obtain that $J_2$ satisfies , 
$$|J_2|\leq c|x|^{2\ga-\frac{N+1}{2}}.
$$
Combining these estimates on $J_1$ and $J_2$ with (\ref{hel17}), we obtain
          \begin{equation}\label{hel19}
|I_{2\ga}\ast \Phi_1|(x)\leq c'|x|^{2\ga-\frac{N+1}{2}}\qquad\text{if }|x|\geq 1.
     \end{equation}
     %%%%%%%%%%%%%%%%%%%%%%%%%%%%%%%%%%%%%%%%%%%%%%%%%%%%%%%%%%%%%%%%%%%%%%%%%%%%%%%%%%%%%%%%%%%%%FIRST CASE TREATED%%%%%%%%%%%%%%%%%%%%%%%%%%%%%%%%%%%%%%%%%%%%%%%%%%%%%%%%
- If $0<\ga\leq \frac 12$, we cannot use directly the integration by parts method since the function 
$(r,\gs)\mapsto (r-\gr\langle\gs,\gth\rangle)(r^2+\gr^2-2r\gr\langle\gs,\gth\rangle)^{\ga-\frac {N+2}2}$ is not integrable near $(\gr,\gth)$. 
We have for $|x|>2$
$$\BA{lll}\dsps
\int_{B^c_1} \!\!|x-y|^{2\ga-N}|y|^{\frac {1-N}{2}}\sin\left(|y|-\frac{(N-1)\gp}{4}\right)dy=\int_{1<|y|<2|x|} \!\!\!\!\!\!\!\! |x-y|^{2\ga-N}|y|^{\frac {1-N}{2}}\sin\left(|y|-\frac{(N-1)\gp}{4}\right)dy\\[4mm]
\phantom{\dsps=
\int_{B^c_1} \!\!|x-y|^{2\ga-N}|y|^{\frac {1-N}{2}}\sin\left(|y|-\frac{(N-1)\gp}{4}\right)dy}
\dsps +\int_{|y|>2|x|} \!\!\!\!\!\!|x-y|^{2\ga-N}|y|^{\frac {1-N}{2}}\sin\left(|y|-\frac{(N-1)\gp}{4}\right)dy
\\[4mm]
\phantom{\dsps
\int_{B^c_1} \!\!|x-y|^{2\ga-N}|y|^{\frac {1-N}{2}}\sin\left(|y|-\frac{(N-1)\gp}{4}\right)dy}
=K_1+K_2.
\EA$$
We estimate the terms $K_j$ by two different methods. 
Using again the spherical coordinates introduced above,
$$\BA{lll}
\dsps |K_2|\leq \left|\int_{2|x|}^\infty\int_{S^{N-1}}r^{\frac{N-1}{2}}(r^2+\gr^2-2r\gr\langle\gs,\gth\rangle)^{\ga-\frac N2}\sin\left(r-\frac{(N-1)\gp}{4}\right)dSdr\right|
\\[4mm]
\phantom{\dsps |K2|}\dsps\leq
\left|\int_{S^{N-1}}\left[r^{\frac{N-1}{2}}(r^2+\gr^2-2r\gr\langle\gs,\gth\rangle)^{\ga-\frac N2}\cos\left(r-\frac{(N-1)\gp}{4}\right)\right]_{r=2|x|}^\infty dS\right|
\\[4mm]
\phantom{=\dsps |K2|---}\dsps+\left|\int_{2|x|}^\infty\int_{S^{N-1}}\frac{d}{dr}\left(r^{\frac{N-1}{2}}(r^2+\gr^2-2r\gr\langle\gs,\gth\rangle)^{\ga-\frac N2}\right)\cos\left(r-\frac{(N-1)\gp}{4}\right) dSdr\right|\\[4mm]
\phantom{\dsps |K2|}
=K_{2,1}+K_{2,2}.
\EA$$
Clearly $K_{2,1}\leq c|x|^{2\ga-\frac{N+1}{2}}$. Using (\ref{der1}) we have
$$\left|\frac{d}{dr}\left(r^{\frac{N-1}{2}}(r^2+\gr^2-2r\gr\langle\gs,\gth\rangle)^{\ga-\frac N2}\right)\right|
\leq cr^{2\ga-\frac{N+3}{2}},
$$
which yields $K_{2,2}\leq c|x|^{2\ga-\frac{N+1}{2}}$ by integration. \\
For estimating $K_1$ we make the change of variable already but with a totally different point of view,
\begin{equation}\label{hel33}\BA{lll}\dsps K_1= |x|^{2\ga-\frac{N-1}{2}}\int_{\frac1{|x|}<|z|<2} \!\!\!\!\!\!\!\! |{\bf e}_N-z|^{2\ga-N}|z|^{\frac {1-N}{2}}\sin\left(|x||z|-\frac{(N-1)\gp}{4}\right)dz\\
[4mm]\phantom{\dsps K_1}\dsps
=|x|^{2\ga-\frac{N-1}{2}}\int_{B_2} |{\bf e}_N-z|^{2\ga-N}|z|^{\frac {1-N}{2}}\sin\left(|x||z|-\frac{(N-1)\gp}{4}\right)dz +O(|x|^{2\ga-N}).
\EA\end{equation}
We estimate the main term by real interpolation theory. Let $T_\gr$ be the mapping from $L^1(B_2)$ to $\BBR$ defined by
\begin{equation}\label{int1}
f\mapsto T_\gr(f):=\int_{B_2}f(z)\sin\left(\gr|z|-\frac{(N-1)\gp}{4}\right)dz,
\end{equation}
clearly
\begin{equation}\label{int3}
|T_\gr(f)|\leq \norm f_{L^1(B_2)}.
\end{equation}
Also, if $f\in W^{1,1}(B_2)$, then $r\mapsto f(r,.)$ is continuous from $W^{1,1}(B_2)$ into $L^1(S^{N-1})$, therefore
$$\BA{lll}\dsps \int_{B_2}f(z)\sin\left(\gr|z|-\frac{(N-1)\gp}{4}\right)dz=-\frac {\cos\left(2\gr-\frac{(N-1)\gp}{4}\right)}{\gr}\int_{|z|=2}f(2,.))dS\\[4mm]
\phantom{-----------------}
\dsps + \frac{1}{\gr}\int_{B_2}\langle \nabla f,\frac z2\rangle \cos\left(\gr|z|-\frac{(N-1)\gp}{4}\right)dz.
\EA$$
This implies 
\begin{equation}\label{int3}
|T_\gr(f)|\leq \frac{c}{\gr}\norm f_{W^{1,1}(B_2)},
\end{equation}
where $c>0$ is independent of $\gr$. By Lions-Petree's real interpolation method theorem  \cite{Li-Pe} and because 
$$W^{s,1}(B_2)=\left[W^{1,1}(B_2),L^1(B_2)\right]_{s,1},
$$
for $0<s<1$, with the notations therein, we derive
\begin{equation}\label{int4}
|T_\gr(f)|\leq \frac{c}{\gr^s}\norm f_{W^{s,1}(B_2)}.
\end{equation}
The function $z\mapsto F(z):= |{\bf e}_N-z|^{2\ga-N}|z|^{\frac {1-N}{2}}$ belongs to $W^{1,1}(B_2)$ for any $\ga>\frac 12$, to $W^{2\ga-\ge,1}(B_2)$ for any $\ge>0$ if $\ga\leq \frac 12$. This implies that for any $\ge>0$ there exists $C_\ge>0$ such that 
\begin{equation}\label{int5}
|K_1|\leq C_\ge|x|^{2\ga+\ge-\frac{N+1}{2}}\quad\text{for all }|x|\geq 2,
\end{equation}
and finally we prove
 \begin{equation}\label{hel34}
|I_{2\ga}\ast \Phi_1|(x)\leq c'|x|^{2\ga+\ge-\frac{N+1}{2}}\qquad\text{if }|x|\geq 1.
     \end{equation}
    %%%%%%%%%%%%%%%%%%%%%%%%%%%%%%%%%%%%%%%%%%%%%%%%%%%%%%%%%%%%%%%%%%%%%%%%%%%%%%%%%%%%%%%%%%%%%FIRST CASE TREATED%%%%%%%%%%%%%%%%%%%%%%%%%%%%%%%%%%%%%%%%%%%%%%%%%%%%%%%%
   Combining (\ref{hel17}), (\ref{hel19})  and (\ref{hel34}), we obtain the estimate (\ref{hel 2}), valid for any $\ge>0$ if $\ga\leq\frac 12$ and $|x|\geq 1$ ,
which is the claim. \hfill$\Box$ \medskip 
%%%%%%%%%%%%%%%%%%%%%%%%%%%%%%%%%%%%%%%%%%%%%%%%%%%%%%%%%%%%%%%%%%%%%%%%%%%%%%%%%%%%%%%%%%%%%%%%%%%%%%%%%%%%%%%%%%%%%%%%%%%%%%%%%%%%%%%%%%%%%%
%%%%%%%%%%%%%%%%%%%%%%%%%%%%%%%%%%%%%%%%%%%%%%%%%%%%%%%%%%%%%%%%%%%%%%%%%%%%%%%%%%%%%%%%%%%%%%%%%STATIONNARY%%%%%%%%%%%%%%%%%%%%%%%%%%%%%%%%%
%%%%%%%%%%%%%%%%%%%%%%%%%%%%%%%%%%%%%%%%%%%%%%%%%%%%%%%%%%%%%%%%%%%%%%%%%%%%%%%%%%%%%%%%%%%%%%%%%%%%%%%%%%%%%%%%%%%%%%%%%%%%%%%%%%%%%%%%%%%%%%
%%%%%%%%%%%%%%%%%%%%%%%%%%%%%%%%%%%%%%%%%%%%%%%%%%%%%%%%%%%%%%%%%%%%%%%%%%%%%%%%%%%%%%%%%%%%%%%%%%%%%%%%%%%%%%%%%%%%%%%%%%%%%%%%%%%%%%%%%%%%%%
 \medskip

 \noindent{\it Proof of Theorem \ref{th 8.1}. } 
 Let 
\begin{equation}\label{star1}u_*(x)=\int_0^1\cP_{\ln}(t,x) \dd t=\int_0^1\cP_{0}(t)|x|^{2t-N} \dd t\quad {\rm for}\ x\not=0. \end{equation}
Clearly $u_*$ is radially symmetric and decreasing with respect to $|x|$. Since $\cP_{\ln}$ satisfies (\ref{e 1-v}), we integrate in $t$ over $(0,1)$. Using the the Fourier isomorphism between $\CS(\BBR^N)$ and $\CS'(\BBR^N)$, we have for any $\varphi\in\cS(\BBR^N)$,
$$\BA{lll}\dsps\langle \cL_\Gd u_*,\hat\varphi\rangle=\langle \cL_\Gd\left(\int_0^1\cP_{\ln}(t,\cdot)dt\right) ,\hat\varphi\rangle=\langle2\ln |\cdot|\cF\left[\int_0^1\cP_{\ln}(t,\cdot)dt\right],\varphi\rangle\\[4mm]\phantom{\langle \cL_\Gd u_*,\hat\varphi\rangle}\dsps=\int_0^1\langle2\ln |.|\cF\left[\cP_{\ln}(t,\cdot)\right],\varphi\rangle dt
=\int_0^1\langle\partial_t\cF[\cP_{\ln}(t,\cdot)],\varphi\rangle dt\\[4mm]\phantom{\langle \cL_\Gd u_*,\hat\varphi\rangle}\dsps
=\langle \cF[\cP_{\ln}(0,\cdot)],\varphi\rangle -\langle \cF[\cP_{\ln}(1,\cdot)],\varphi\rangle=\langle \cP_{\ln}(0,\cdot)-\cP_{\ln}(1,\cdot),\hat \varphi\rangle
\EA$$
 and  we get
 \begin{equation}\label {XX1}\begin{array}{lll}\dsps
\lnlap u_*(x)  =\cP_{\ln}(0,x)- \cP_{\ln}(1,x)=-\frac{\Gg(\frac{N}{2}-1)}{\Gg(1)\sqrt{2}\gp^{N/2}}|x|^{2-N}\quad{\rm for}\  |x|>0.
\end{array}\end{equation}
 It follows from \cite[Theorem 5]{St} that $\cP_{\ln}(1,\cdot)$  is the fundamental solution of $-\Gd$ in $\BBR^N$ expressed by $I_2(x)=c_{N,1}|x|^{2-N}$. \smallskip

\noindent{\it Step 1: Properties of the function $u_*$. } \\
We have that for  $0<|x|<1$,  
 \begin{eqnarray*}
u_*(x)&=& \int_0^1 \cP_0(t) e^{(2t-N)\ln |x|} \, \dd t 
 \\[2mm] &=&  \left[\frac1{2\ln |x|}\cP_0(t) |x|^{2t-N}\right]_{t=0}^{t=1} -   \frac1{2\ln |x|}  \int_0^1 \cP_0'(t) e^{(2t-N)\ln |x|} \dd t
 \\[2mm] &=&  \left[\frac1{2\ln |x|}\cP_0(t) |x|^{2t-N}\right]_{t=0}^{t=1} - \left[\frac1{(2\ln |x|)^2}\cP_0'(t) |x|^{2t-N}\right]_{t=0}^{t=1} 
  \\[2mm] && \phantom{-------------}\  +\frac1{(2\ln |x|)^2}  \int_0^1 \cP_0''(t) e^{(2t-N)\ln |x|} \dd t
  \\[2mm] &=&   \frac{\cP_0(1)}{2\ln |x|}  |x|^{2-N}  +  \frac{\cP_0'(0)}{(2\ln |x|)^2} |x|^{-N}- 
  \frac{\cP_0'(1)}{(2\ln |x|)^2} |x|^{2-N} 
  \\[2mm] && \phantom{-------------} \  +\frac1{(2\ln |x|)^2}  \int_0^1 \cP_0''(t) e^{(2t-N)\ln |x|} \dd t.
\end{eqnarray*}
From Lemma \ref{re 2.1}, $\cP_0'(0)=\pi^{-N/2}\Gamma(\tfrac N2)$ and $\cP_0''$ is bounded on $[0,1]$, hence
 \begin{eqnarray*}
 \Big|\frac1{(2\ln |x|)^2}  \int_0^1 \cP_0''(t) e^{(2t-N)\ln |x|} \dd t\Big| \leq 
 \frac c{(\ln |x|)^2} \int_0^1   e^{(2t-N)\ln |x|} \dd t
 \leq   c \frac {|x|^{-N}}{(\ln |x|)^3}.  
 \end{eqnarray*}
Therefore, we obtain that for $0<|x|<\frac1e$,
   \begin{equation}\label{eq 8.4}
\Big| u_*(x)- \frac{\Gamma(\tfrac N2)}{4\pi^{N/2}} \frac{|x|^{-N}}{(\ln |x|)^2}\Big|\leq c \frac {|x|^{-N}}{(-\ln |x|)^3}.
  \end{equation} 
   It is worth noting that for $|x|=1$,  we have 
  $$u_*(x)=\int_0^1\cP_{0}(t)\dd t.$$
 For the asymptotic behavior of $u_*$ when $|x|\to\infty$, we have that for  $ |x|>1$,  
 \begin{equation}\label {eqX}\BA{lll}\dsps
0<u_*(x)\leq c\int_0^1 t e^{(2t-N)\ln |x|} \, \dd t 
 \\[4mm]\phantom{0<u_*(x)}\dsps = c\left(\left[ \frac1{2\ln |x|} t |x|^{2t-N}\right]_{t=0}^{t=1} -  \left[\frac1{(2\ln |x|)^2}  |x|^{2t-N}\right]_{t=0}^{t=1} \right)
  \\[5mm]\phantom{0<u_*(x)}\dsps \leq   \frac{c}{ \ln |x|}   |x|^{2-N}+ \frac{c}{ (\ln |x|)^2}   |x|^{-N}.   
\EA\end{equation}
Since $\cP_{\ln}$ is a weak solution of (\ref{eq 1.1-d_0}), then  taking $T'=1$ and $\varphi(t,\cdot)=\varphi\in C^1_c(\R^N) $ in (\ref{def w1-d_0})  we infer  that 
\begin{eqnarray*}
0 &=& \int_0^{1} \int_{\mathbb{R}^N} \big(-\partial_t \varphi+\lnlap \varphi\big)\cP_{\ln} \dd x\dd t-\varphi(0,0) +\int_{\mathbb{R}^N} \cP_{\ln}\left(1, x \right) \varphi(1,x)dx
\\[2mm]&=&  \int_{\mathbb{R}^N} u_* \lnlap \varphi  \dd x -\varphi(0) +\int_{\mathbb{R}^N} \cP_{\ln}\left(1, x \right)  \varphi\dd x.
\end{eqnarray*}
As a consequence,  $u_*$ is a weak solution of 
 \begin{equation}\label{e 7.1} 
\lnlap  u_* = \delta_0-\cP_{\ln}(1,\cdot)\quad {\rm in}\ \cD'(\R^N), 
\end{equation} 
and a classical solution of 
 \begin{equation}\label{e 7.2} 
\lnlap  u_* =  -\cP_{\ln}(1,\cdot)\quad {\rm in}\ \R^N\setminus\{0\}.
\end{equation}

\noindent {\it Step 2: There exists $v_1$ such that
 \begin{equation}\label{eq 8.2}
 \lnlap v_1 =\cP_{\ln}(1,\cdot)\quad {\rm in}\ \cD'(\R^N).
  \end{equation}}
  Given $f\in L^1(\BBR^N, (1+|x|)^{2T-N}dx)$, denote  
 $$   u_f(x)=\int_0^1 \cP_{\ln}(t,\cdot)\ast f(t,x) \dd t. $$ 
Then 
 \begin{eqnarray*}
\lnlap u_f(x)  &=&  \int_0^1 \int_{\R^N}\lnlap\cP_{\ln}(t,x-y)f(y) \dd y  \dd t 
\\[2mm] &=&- \int_0^1 \partial_t\int_{\R^N} \cP_{\ln}(t,x-y)f(y) \dd y  \dd t 
 =   f(x)- \cP_{\ln}(1,\cdot)\ast f (x).
\end{eqnarray*}
 Now we aim to choose $f$ such that
 $$ f(x)- \cP_{\ln}(1,\cdot)\ast f (x)=\cP_{\ln}(1,x).$$
 Since  $\cP_{\ln}(1,\cdot)$  is the fundamental solution of $-\Gd$,  $f$ satisfies 
 \begin{equation}\label{eq 8.3} 
 -\Delta f-f=\delta_0\quad {\rm in}\ \ \cD'(\R^N).
 \end{equation}
 Therefore,  $\Phi_{1}$ is a solution of (\ref{eq 8.3}) and  
 \begin{equation}  \label{eq H6}
v_1(x)=\int_0^1(\cP_{\ln}(t,\cdot)\ast\Phi_{1})(x,t)dt=\int_0^1\cP_{0}(t)(I_{2t}\ast\Phi_{1})(x)dt,
\end{equation}
satisfies (\ref{eq 8.2}).

 \medskip

\smallskip

\noindent {\it Step 3: Asymptotic estimates of  $v_1$}\\
{\it  3-(i): Estimate near $x=0$.} Since for $0\leq t\leq 1$ we have that $\frac {1}{c}t\leq \cP_{0}(t)\leq ct,$ we deduce from Lemma \ref{est-l} that if  $|x|<\frac 1e$ we have 
when $N\geq 4$,
 \begin{equation}  \label{hel55}|v_1(x)|\leq c|x|^{2-N}\int_0^1t|x|^{2t}dt= c|x|^{2-N}\left(\frac{|x|^2}{2\ln |x|}-\frac{|x|^2}{4\ln^2|x|}+\frac{1}{4\ln^2|x|}\right)\leq c\frac{|x|^{2-N}}{4\ln^2|x|},
\end{equation}
and when $N=3$ 
 \begin{equation}  \label{hel56}|v_1(x)|\leq c|x|^{-1}\int_0^{\frac12}t|x|^{2t}dt+\frac c2\leq  \frac{c}{4|x|\ln^2|x|}+\frac c2.
\end{equation}
{\it  3-(i): Estimate as $|x|\to\infty$.} 
If  $|x|>2$ we have for any $\ge>0$
 \begin{equation}  \label{hel57}|v_1(x)|\leq C_\ge|x|^{-\frac{N+1}{2}+\ge}\int_0^{\frac 12}t|x|^{2t}dt+|x|^{-\frac{N+1}{2}}\int_{\frac 12}^1t|x|^{2t}dt\leq \frac{c'}{|x|^{\frac{N-3}{2}}\ln |x|}.
\end{equation}

\noindent {\it Step 4: Asymptotics of the fundamental solution of $\CL_\Gd$} \smallskip

\noindent From (\ref{eq 7.2})-(\ref{eq 8.2}) the function $\Phi_{\ln}=u_*+v_1$ satisfies 
$$\CL_{\Gd}\Phi_{\ln}=\gd_0\qquad\text{in }\CD'(\BBR^N).
$$
Using (\ref{eq 8.4}), it follows from (\ref{hel55})-(\ref{hel56}) that for $0<|x|\leq e^{-1}$ there holds for some $c>0$
   \begin{equation}\label{AS15}
\left| \Phi_{\ln}(x)- \frac{\Gamma(\tfrac N2)}{4\pi^{N/2}} \frac{|x|^{-N}}{(\ln |x|)^2}\right|\leq c\frac{|x|^{-N}}{\ln^3|x|}.
  \end{equation} 
Using (\ref{hel57}) we obtain, for $|x|$ large enough and some positive constant $c>0$,
   \begin{equation}\label{AS16}
|\Phi_{\ln}(x)|\leq c\frac{|x|^{\frac{3-N}{2}}}{\ln |x|}.
  \end{equation} 
This ends the proof. \hfill$\Box$
%\medskip
 
 \begin{remark}\label{int} It is important to notice from (\ref{AS15}) that $\Phi_{\ln}\in L^1_{loc}(\R^N)$.
 \end{remark}

% \subsection{Poisson problem in $\R^N$}

 \subsection{Existence in bounded domains}

 Let $\Omega$ be a bounded Lipchitz domain and $\cH(\Omega)$ denote the space of all measurable functions $u:\R^N\to \R$ with $u \equiv 0$ in $\R^N \setminus \Omega$ and
$$
\int \!\!\!\! \int_{\stackrel{x,y \in \R^N}{\text{\tiny $|x-y|\!\le\! 1$}}} \frac{(u(x)-u(y))^2}{|x-y|^N} dx dy <+\infty.
$$
Note that  $\cH(\Omega)$ is a Hilbert space under the inner product
$$
\mathcal{E}(u,w)=\frac{c_{_N}}2 \int \!\!\int_{\stackrel{x,y \in \R^N}{\text{\tiny $|x-y|\!\le\! 1$}}} \frac{(u(x)-u(y))(w(x)-w(y))}{|x-y|^N}
dx dy
$$
and with the associated norm $\norm{u}_{\cH(\Omega)}=\sqrt{\mathcal{E}(u,u)}$. By \cite[Theorem 2.1]{CP}, {\it the embedding: $\cH(\Omega) \hookrightarrow L^2(\Omega)$ is compact. }  
Let
$$
\cE_L: \cH(\Omega) \times \cH(\Omega) \to \R, \qquad \cE_L(u,w)= \mathcal{E}(u,w) - \int_{\R^N} \Bigl( \j * u  -\rho_N u \Bigr)w\,dx.
$$
where 
$$
\j: \R^N \to \R,\qquad \j(z)= c_{_N} 1_{\R^N \setminus B_1}(z)|z|^{-N}.
$$
For additional properties of $\cE_L$ see Sections 3, 4 in \cite{CT}. 
If $\lambda_1^{^L}(\Omega)>0$, there is some $c>0$ such that
$$\cE_L(u,u)\geq c\cE(u,u).  $$

We recall that $\eta_0\in C_c^\infty(\R^N)$ is the cut-off  function defined in the proof of Theorem \ref{thm2} and $\eta(x)=\eta_0(|x|)$

\begin{lemma}\label{trunc} The function $\CL_\Gd(\eta\Phi_{\ln})$ is uniformly bounded in $\BBR^N\setminus\{0\}$.
\end{lemma}
\proof Set $v=\eta\Phi_{\ln}$ in $\BBR^N\setminus\{0\}$, then 
\begin{equation}\label{aprrox1}\BA{lll}\dsps
\CL_\Gd v(x)=c_{_N}\int_{B_1(x)}\frac{v(x)-v(y)}{|x-y|^N}dy-c_{_N}\int_{B^c_1(x)}\frac{v(y)}{|x-y|^N}dy+\gr_Nv(x)\\[4mm]
\phantom{\CL_\Gd v(x)}\dsps
=c_{_N}\eta(x)\int_{B_1(x)}\frac{\Phi_{\ln}(x)-\Phi_{\ln}(y)}{|x-y|^N}dy+c_{_N}\int_{B_1(x)}(\eta(x)-\eta(y))\frac{\Phi_{\ln}(y)}{|x-y|^N}dy\\[4mm]
\phantom{\CL_\Gd v(x)----------------}\dsps-c_{_N}\int_{B^c_1(x)}\frac{\eta(y)\Phi_{\ln}(y)}{|x-y|^N}dy+\gr_Nv(x)\\[4mm]
\phantom{\CL_\Gd v(x)}\dsps
=c_{_N}\int_{\BBR^N}\frac{(\eta(x)-\eta(y))\Phi_{\ln}(y)}{|x-y|^N}dy.
\EA\end{equation}
Then
$$\BA{lll}\dsps\CL_\Gd v(x)=c_{_N}\int_{B_{\frac 12}}\!\frac{(\eta(x)-1)\Phi_{\ln}(y)}{|x-y|^N}dy+c_{_N}\int_{B^c_{ 1}}\!\frac{\eta(x)\Phi_{\ln}(y)}{|x-y|^N}dy
+c_{_N}\int_{B_1\setminus B_{\frac 12}}\!\!\frac{(\eta(x)-\eta(y))\Phi_{\ln}(y)}{|x-y|^N}dy\\[2mm]
\phantom{\CL_\Gd v(x)}=A_1(x)+A_2(x)+A_3(x).
\EA$$
(i) If $y\in B_1\setminus B_{\frac 12}$, the function $\Phi_{\ln}$ is bounded. If $x\in B_2$, then $B_1\setminus B_{\frac 12}\subset B_3(x)$, hence
$$
|A_3(x)|\leq c_{_N}\norm{\Phi_{\ln}}_{L^\infty(B_1\setminus B_{\frac 12})}\norm{\nabla\eta}_{L^\infty}\int_{B_3(x)}|x-y|^{1-N}dy
\leq 3c_{_N}\gw_{_{N}}\norm{\Phi_{\ln}}_{L^\infty(B_1\setminus B_{\frac 12})}\norm{\nabla\eta}_{L^\infty}.
$$
We recall that $\gw_{_{N}}$ is the (N-1)-volume of the unit sphere $S^{N-1}$. If $x\in B^c_2$, then 
$$\left|\frac{(\eta(x)-\eta(y))\Phi_{\ln}(y)}{|x-y|^N}\right|\leq 2\norm{\Phi_{\ln}}_{L^\infty(B_1\setminus B_{\frac 12})}.
$$
Combining these two estimates we obtain
$$|A_3(x)|\leq c_1\quad\text{for all }x\in\BBR^N\setminus\{0\}
$$
for some $c_1>0$.\smallskip

\noindent (ii) If $y\in B^c_1$, we write
$$\BA{lll}\dsps A_2(x)=c_{_N}\int_{B^c_{ 1}}\!\frac{\eta(x)\Phi_{\ln}(y)}{|x-y|^N}dy=c_{_N}\int_{B_2\setminus B_{ 1}}\!\frac{\eta(x)\Phi_{\ln}(y)}{|x-y|^N}dy+c_{_N}\int_{B^c_2}\!\frac{\eta(x)\Phi_{\ln}(y)}{|x-y|^N}dy\\[3mm]
\phantom{B(x)}
\dsps=A_{21}(x)+A_{22}(x).
\EA$$
Then $A_{21}(x)=0$ if $|x|>1$ and
$$A_{21}(x)\leq c_{_N}\norm{\nabla\eta}_{L^\infty}\norm{\Phi_{\ln}}_{L^\infty(B_2\setminus B_{1})}\int_{B_2\setminus B_{ 1}}\frac{dy}{|x-y|^{N-1}}\quad\text{for all }x\in B_1.
$$
This implies that $B_1(x)$ is uniformly bounded. Concerning $B_2(x)$ which vanishes also in $B_1^c$, we have
$$A_{22}(x)\leq c_{_N}\int_{B^c_2}\frac{|\Phi_{\ln}(y)|}{(|y|-1)^N}dy\quad\text{for all }x\in B_1.
$$
By (\ref{AS16}) the integral is convergent and we obtain 
Combining these two estimates we obtain that
$$|A_2(x)|\leq c\quad\text{for all }x\in\BBR^N\setminus\{0\},
$$
 for some $c>0$.
\smallskip

\noindent (iii) If $x\in B_{\frac12}$, $A_1(x)=0$. If $|x|\geq 1$, 
$$|A_1(x)|\leq \frac{c_{_N}}{(|x|-\frac12)^N}\int_{B_\frac12}|\Phi_{\ln}(y)|dy,
$$
a quantity which is uniformly bounded. If $\frac12<|x|<1$ we write since $\eta(y)=1$ therein
$$\BA{lll}\dsps
|A_1(x)|\leq c_{_N}\int_{B_{\frac12}\setminus B_{\frac13}}\frac{|\eta(y)-\eta(x)|\varphi_{\ln}(y)}{|x-y|^N}dy+ c_{_N}\int_{B_{\frac13}}\frac{|\eta(y)-\eta(x)|\varphi_{\ln}(y)|}{|x-y|^N}dy\\[4mm]
\phantom{|A(x)|}\dsps
\leq c_{_N}\norm{\Phi_{\ln}}_{L^\infty(B_{\frac12}\setminus B_{\frac13})}\norm{\nabla\eta}_{L^\infty}\int_{B_{\frac12}\setminus B_{\frac13}}\frac{1}{|x-y|^{N-1}}dy
+c_{_N}6^N\int_{B_{\frac13}}|\varphi_{\ln}(y)|dy\\[4mm]
\phantom{|A(x)|}\dsps
\leq c_{_N}\norm{\Phi_{\ln}}_{L^\infty(B_{\frac12}\setminus B_{\frac13})}\norm{\nabla\eta}_{L^\infty}\int_{B_{\frac32}}\frac{dz}{|z|^{N-1}}+
c_{_N}6^N\norm{\varphi_{\ln}}_{L^1(B_{\frac13})}.
\EA$$
Hence it follows that for some $c>0$
$$|A_1(x)|\leq c\quad\text{for all }x\in\BBR^N\setminus\{0\}.
$$
Combining the three estimates we obtain the claim.  \hfill$\Box$
\begin{lemma}\label{trunc} Set $h=\CL_{\Gd}(\eta\Phi_{\ln}){\bf1}_{\Gw}$. If $0\notin\gs(\CL_\Gd,\Gw)$ there exists a unique function $w\in \BBH\cap \cC(\overline\Gw)$ such that 
\begin{equation}\label{approx2}
\left\{ \arraycolsep=1pt
\begin{array}{lll}
  \CL_\Gd w=h\ \   & {\rm in\ }\, \Gw,\\[2mm]
 \phantom{ \lnlap     }
    w=0 \ \ &{\rm in\ }\, \Gw^c.
\end{array}
\right.
\end{equation}
\end{lemma}
\noindent\proof The function $h$ is bounded hence it belongs to $L^2(\Gw)$. By \cite[Theorem 2.1]{CP} the embedding of $\BBH(\Gw)$ into 
$L^2(\Gw)$ is compact and by  \cite[Theorem 3.1]{CT} $C^\infty_c(\Gw)$ is dense in $\BBH(\Gw)$. Then the operator $\CL_\Gd$ is a Fredholm operator and since $0$ is not an eigenvalue, there exists a unique $w\in\BBH(\Gw)$ such that (\ref{approx2}) holds in the weak sense, and by density such that 
$$\CE_L(w,\gz)=\int_{\Gw}h\gz dx\quad\text{for all }\gz\in\BBH(\Gw).
$$
By \cite[Theorem 1.11]{CT} the function $w$ is continuous and it satisfies
\begin{equation}\label{approx3}
|w(x)|=O\left(\frac{1}{|\ln\gr(x)|^\gt}\right)\quad\text{for all }\gt\in \left(0,\frac 12\right),
\end{equation}
near $\partial\Gw$ where we recall that $\gr(x)=$dist$(x,\partial\Gw)$.  This ends the proof. \hfill$\Box$

%%%%%%%%%%%%%%%%%%%%%%%%%%%%%%%%%%%%%%%%%%%%%%%%%%%%%%%%%%%%%%%%%%%%%%%%%%%%%%%%%%%%%%%%%%%%%%%%%%%%%%%%%%%%%%%%%%%%%%%%%%%%%%%%%%%%%%%%%%%%%%%%%%%%%%%%%%%%%%%%%%%%%%%%%%%%%%%%%%%%%%%%%%%%%%%%%%%%%%%%%%%%%%%%%%%%%%%%%%%%%%%%%%%%%%%%%%%%%%%%%%%%%%%%%%%%%%%%%%%%%%%%%%%%%%%%%%%%%%%%%%%%
 \medskip
 
\noindent{\it Proof of Theorem \ref{th 7.1}. }   Set $\Phi_{\ln}^\Gw=\eta\Phi_{\ln}-w$. Then $\Phi_{\ln}^\Gw$  vanishes in $\Gw^c$ and satisfies
$$\CL_\Gd\Phi_{\ln}^\Gw=0\qquad\text{in }\ \Gw\setminus\{0\}.
$$
Estimates (\ref{E1})  follows from (\ref{AS15}) and (\ref{approx3}). The same holds with (\ref{7xx}) and actually this estimate is already obtained in \cite[Theorem 1.11]{CT}). Now let $\gz\in C^\infty_c(\Gw)$, then
$$\int_{\R^N}\Phi_{\ln}^\Gw\CL_\Gd\gz dx=\int_{\R^N}\eta\Phi_{\ln}\CL_\Gd\gz dx-\CE_L(w,\gz) =A_4-A_5.
$$
We have that
$$\BA{lll}\dsps A_4=c_{_N}\int_{\R^N}\eta(x)\Phi_{\ln}(x)\int_{B_1(x)}\frac{\gz(x)-\gz(y)}{|x-y|^N}dy dx-c_{_N}\int_{\R^N}\eta(x)\Phi_{\ln}(x)\int_{B^c_1(x)}\frac{\gz(y)}{|x-y|^N}dy dx\\[4mm]\phantom{--------------------------}\dsps+\gr_N\int_{\R^N}\eta(x)\Phi_{\ln}(x)\gz(x) dx\\[4mm]\phantom{A}\dsps
= c_{_N}\int_{\R^N}\Phi_{\ln}(x)\int_{B_1(x)}\frac{(\eta\gz)(x)-(\eta\gz)(y)}{|x-y|^N}dy dx-c_{_N}\int_{\R^N}\Phi_{\ln}(x)\int_{B^c_1(x)}\frac{(\eta\gz)(y)}{|x-y|^N}dydx\\[4mm]\phantom{A----}\dsps-c_{_N}\int_{\R^N}\Phi_{\ln}(x)\int_{\BBR^N}\frac{\eta(x)-\eta(y)}{|x-y|^N}\gz(y)dydx+\gr_N\int_{\R^N}\eta(x)\Phi_{\ln}(x)\gz(x) dx
\\[4mm]\phantom{A}\dsps=\gz(0)-c_{_N}\int_{\R^N}\Phi_{\ln}(x)\int_{\BBR^N}\frac{\eta(x)-\eta(y)}{|x-y|^N}\gz(y)dydx,
\EA,$$
and, by (\ref{aprrox1}),
$$\BA{lll}\dsps A_5
=\int_{\R^N} \gz\CL_\Gd(\eta\Phi_{\ln})  dx =c_{_N}\int_{\R^N} \gz(x)\int_{\BBR^N}\frac{\eta(x)-\eta(y)}{|x-y|^N}\Phi_{\ln}(y)dydx\\[4mm]
\phantom{\dsps B
=\int_{\R^N} \gz\CL_\Gd(\eta\Phi_{\ln})  dx}\dsps
=-c_{_N}\int_{\R^N}\Phi_{\ln}(x)\int_{\R^N}\frac{\eta(x)-\eta(y)}{|x-y|^N}\gz(y)dydx.
\EA$$
 This implies that $\CL_\Gd\Phi_{\ln}^\Gw=\gd_0$ in $\CD'(\Gw)$.\hfill$\Box$
%%%%%%%%%%%%%%%%%%%%%%%%%%%%%%%%%%%%%%%%%%%%%%%%%%%%%%%%%%%%%%%%%%%%%%%%%%%%%%%%%%%%%%%%%%%%%%%%%%%%%%%%%%%%%%%%%%%%%%%%%%%%%%%%%%%%%%%%%%%%%%%%%%%%%%%%%%%%%%%%%%%%%%%%%%%%%%%%%%%%%%%%%%%%%%%%%%%%%%%%%%%%%%%%%%%%%%%%%%%%%%%%%%%%%%%%%%%%%%%%%%%%%%%%%%%%%%%%%%%%%%%%%%%%%%%%%%%%%%%%%%%%%%%%%%%%%%%
%%%%%%%%%%%%%%%%%%%%%%%%%%%%%%%%%%%%%%%%%%%%%%%%%%%%%%%%%%%%%%%%%%%%%%%%%%%%%%%%%%%%%%%%%%%%%%%%%%%%%%%%%%%%%%%%%%%%%%%%%%%%%%%%%%%%%%%%%%%%%%%%%%%%%%%%%%%%%%%%%%%%%%%%%%%%%%%%%%%%%%%%%%%%%%%%%%%%%%%%%%%%%%%%%%%%%%%%%%%%%%%%%%%%%%%%%%%%%%%%%%%%%%%%%%%%%%%%%%%%%%%%%%%%%%%%%%%%%%%%%%%%
%%%%%%%%%%%%%%%%%%%%%%%%%%%%%%%%%%%%%%%%%%%%%%%%%%%%%%%%%%%%%%%%%%%%%%%%%%%%%%%%%%%%%%%%%%%%%%%%%%%%%%%%%%%%%%%%%%%%%%%%%%%%%%%%%%%%%%%%%%%%%%%%%%%%%%%%%%%%%%%%%%%%%%%%%%%%%%%%%%%%%%%%%%%%%%%%%%%%%%%%%%%%%%%%%%%%%%%%%%%%%%%%%%%%%%%%%%%%%%%%%%%%%%%%%%%%%%%%%%%%%%%%%%%%%%%%%
\section{Regularity estimates}

This section is devoted to the proof of some regularity  estimates. We decompose the fundamental solution $\Phi_{\ln}$ of $\CL_\Gd$ in $\R^N$ under the form
\begin{equation}\label{E1}
\Phi_{\ln}=\Phi_{\ln}{\bf 1}_{B_{1}}+\Phi_{\ln}{\bf 1}_{B^c_{1}}=\Phi_{\ln,1}+\Phi_{\ln,2}.
 \end{equation}
 %%%%%%%%%%%%%%%%%%%%%%%%%%%%%%%%%%%%%%%%%%%%
%%%%%%%%%%%%%%%%%%%%%%%%%%%%%%%%%%%%%%%%%%%%%%%%%%%%%%%%%%%%%%%%%%%%%%%%%%%%%%%%%%%%%%%%%%%%%%%%%%%%%%%%%%%%%%%%%%%%%%%%%%%%%%%%%%%%%%%%%%%%%%%%%%%%%%%%%%%%%%%%
\begin{lemma}\label{P1}
Let $N\geq 3$.
If
$f\in L^{\infty}(\R^N)$
then  $\Phi_{\ln,1}\ast f$ is uniformly continuous and  there exists $c>0$ such that 
 \begin{equation}\label{EI-2}
 \|\Phi_{\ln,1}\ast f\|_{\infty}\leq c\|  f\|_{\infty}.
 \end{equation}  
\end{lemma} 
\noindent{\it Proof.}  By estimate (\ref{7xy})  $\Phi_{\ln,1}\in L^1(\R^N)$. Estimate (\ref{EI-2}) follows and classicaly $\Phi_{\ln,1}\ast f$ is uniformly continuous. \hfill$\Box$\medskip
%%%%%%%%%%%%%%%%%%%%%%%%%%%%%%%%%%%%%%%%%%%%
%%%%%%%%%%%%%%%%%%%%%%%%%%%%%%%%%%%%%%%%%%%%%%%%%%%%%%%%%%%%%%%%%%%%%%%%%%%%%%%%%%%%%%%%%%%%%%%%%%%%%%%%%%%%%%%%%%%%%%%%%%%%%%%%%%%%%%%%%%%%%%%%%%%%%%%%%%%%%%%%

Similarly we derive from (\ref{best inf}),

\begin{lemma}\label{P2}
Let $N\geq 3$ , then $\Phi_{\ln,2}\in L^\frac{2N}{N-2}(\R^N)$. 
If
$f\in L^{q'}(\R^N)$ for all $q'\leq \frac{2N}{N+2}$, then  $ \Phi_{\ln,2}\ast f$ is uniformly continuous and 
\begin{equation}\label{EI-3}
 \|\Phi_{\ln,2}\ast f\|_{\infty}\leq c\|  f\|_{\frac{2N}{N+2}}.
 \end{equation}  
\end{lemma} 
%%%%%%%%%%%%%%%%%%%%%%%%%%%%%%%%%%%%%%%%%%%%
%%%%%%%%%%%%%%%%%%%%%%%%%%%%%%%%%%%%%%%%%%%%%%%%%%%%%%%%%%%%%%%%%%%%%%%%%%%%%%%%%%%%%%%%%%%%%%%%%%%%%%%%%%%%%%%%%%%%%%%%%%%%%%%%%%%%%%%%%%%%%%%%%%%%%%%%%%%%%%%%
In the next result we prove some regularity estimate of extended  Marcinkiewicz type for $\Phi_{\ln}$.
 \begin{proposition}\label{Ma1}
 Assume $N\geq 3$, then for some $c=c(N)>0$ there holds
 \begin{equation}\label{Ma-1}
E_{\ln,1}(\gl):=\left\{x\in B_1:|\Phi_{\ln,1}(x)|>\gl\right\}\subset \left\{x\in B_1:|x|\leq c\left(\gl\ln^2\gl\right)^{-\frac1N}\right\}\qquad\text{for all }\gl>1,
 \end{equation}
 hence
  \begin{equation}\label{Ma-1'}
meas\left\{E_{\ln,1}(\gl)\right\}\leq \frac{c}{\gl\ln^2\gl}.
 \end{equation}
Furthermore, if $N>3$, 
  \begin{equation}\label{Ma-2}
E_{\ln,2}(\gl):=\left\{x\in B^c_1:|\Phi_{\ln,2}(x)|>\gl\right\}\subset \left\{x\in B^c_1:|x|\leq \frac{c}{\gl^{\frac 2{N-3}}\ln\left(\frac{1}{\gl}\right)}\right\}\qquad\text{for all }0<\gl<1
 \end{equation}
 which implies
   \begin{equation}\label{Ma-2'}
meas\left\{E_{\ln,2}(\gl)\right\}\leq \frac{c}{\gl^{\frac{2N}{N-3}}|\ln\gl|^N},
 \end{equation}
 and if $N=3$,
   \begin{equation}\label{Ma-3}
E_{\ln,2}(\gl):=\left\{x\in B^c_1:|\Phi_{\ln,2}(x)|>\gl\right\}\subset \left\{x\in B^c_1:|x|\leq e^{\frac c\gl-1}\right\}\qquad\text{for all }0<\gl<1,
 \end{equation}
which yields
    \begin{equation}\label{Ma-3'}
meas\left\{E_{\ln,2}(\gl)\right\}\leq c_{_N}e^{\frac {Nc}\gl}.
 \end{equation}
 \end{proposition}
 %%%%%%%%%%%%%%%%%%%%%%%%%%%%%%%%%%%%%%%%%%%%%%%%%%%%%%%%%%%%%%%%%%%%%%%%%%%%%%%%%%%%%%%%%%%%%%%%%%%%%%%%%%%%%%%%%%%%%%%%%%%%%%%%%%%%%%%%%%%%%%%%%%%%%%%%%%%%%%%
\noindent{\it Proof.}  If $x\in E_{\ln,1}(\gl)$, then from (\ref{7xy}),
$$(1+\ln^2|x|)|x|^N\leq \frac c\gl,
$$
which implies in particular $|x|^N\leq \frac c\gl$, thus $1+\ln^2|x|\geq 1+\frac1{N^2}\ln^2\left(\frac c\gl\right)$, which yields
$$|x|^N\leq \frac {cN^2}{\gl\left(N^2+\ln^2\left(\frac c\gl\right)\right)},
$$
from what estimates inequalities (\ref{Ma-1}) and (\ref{Ma-1'}) follow.\smallskip

\noindent Next, if $x\in E_{\ln,2}(\gl)$, then by (\ref{best inf}), for some constant $a_1>0$
$$\gl\leq a_1\frac{|x|^{\frac {3-N}{2}}}{1+\ln |x|},
$$
which implies $(1+\ln |x|)|x|^{\frac {N-3}{2}}\leq \frac {a_1}\gl$. If $N=3$ we deduce that (\ref{Ma-3}) and (\ref{Ma-3'}) hold.
If $N\geq 4$, we have in particular $1\leq |x|\leq \left(\frac{a_1}{\gl}\right)^{\frac 2{N-3}}$. From the monotonicity properties of the function 
$$r\mapsto (1+\ln r)r^{\frac {N-3}{2}}:=\gm(r) \qquad\text{for }r>1,
$$
we can find some $\gth>0$ such that 
$$\gm\left(\gth\frac{\left(\frac{a_1}{\gl}\right)^{\frac 2{N-3}}}{\ln\left(\frac{c_1}{\gl}\right)} \right)\geq \frac{a_1}{\gl}.
$$
This implies 
$$ E_{\ln,2}(\gl)\subset \left\{x:1\leq |x|\leq \gth\frac{\left(\frac{a_1}{\gl}\right)^{\frac 2{N-3}}}{\ln\left(\frac{a_1}{\gl}\right)} \right\}.
$$
Hence we obtain (\ref{Ma-2}) and (\ref{Ma-2'}).\hfill$\Box$\medskip

In the next result we derive estimate of $\Phi_{\ln}$ in Orlicz classes $L_M$ for some suitable functions $M$. 

%%%%%%%%%%%r%%%%%%%%%%%%%%%%%%%%%%%%%%%%%%%%%
%%%%%%%%%%%%%%%%%%%%%%%%%%f%%%%%%%%%%%%%%%%%%%%%%%%%%%%%%%%%%%%%%%%%%%%%%%%%%%%%%
\begin{theorem}\label{orl}
Let $M:\BBR_+\to\BBR_+$ be a $C^1$ increasing function such that $M(0)=0$.\smallskip

\noindent 1- If $M$ satisfies
    \begin{equation}\label{orl-1}
\int_2^\infty \frac{M(t)}{t^2\ln^2t}dt<\infty,
 \end{equation}
 then $M(|\varphi_{\ln,1}|)\in L^1(B_1)$.\smallskip

\noindent 2- If $N>3$ and $M$ satisfies
     \begin{equation}\label{orl-2}
\int_0^\frac12 \frac{M(t)}{t^\frac{3(N-1)}{N-3}|\ln t|^N}dt<\infty,
 \end{equation}
  then $M(|\varphi_{\ln,2}|)\in L^1(B^c_1)$.\smallskip
  
  \noindent 3- If $N=3$ and $M$ satisfies
       \begin{equation}\label{orl-3}
\int_0^\frac12 \frac{M(t)e^{\frac{Nc}{\gl}}}{\gl^2}<\infty,
 \end{equation}
   then $M(|\varphi_{\ln,2}|)\in L^1(B^c_1)$.
\end{theorem}
\noindent{\it Proof.}  1- Set $e_{\ln,1}(\gl)=meas\left(E_{\ln,1}(\gl)\right)$ and let $\ga\in (0,1)$ such that $\varphi_{\ln,1}(x)\geq 2$ if $0<|x|\leq\ga$. Then
$$\BA{lll}\dsps\int_{B_1}M(|\varphi_{\ln,1}(x)|)dx=\int_{B_\ga}M(|\varphi_{\ln,1}(x)|)dx+\int_{B_1\setminus B_\ga}M(|\varphi_{\ln,1}(x)|)dx\\[4mm]
\phantom{\dsps\int_{B_1}M(|\varphi_{\ln,1}(x)|)dx}\dsps \leq-\int_2^\infty M(t)de_{\ln,1}(t)+\int_{B_1\setminus B_\ga}M(|\varphi_{\ln,1}(x)|)dx.
\EA$$
By standard results on Stieljes integral
$$\BA{lll}\dsps-\int_2^{t_n} M(t)de_{\ln,1}(t)=-\left[M^{\phantom{\frac pp}}\!\!\!\!\!(t)e_{\ln,1}(t)\right]_2^{t_n}+\int_2^{t_n} M'(t)e_{\ln,1}(t)dt\\[4mm]
\phantom{\dsps-\int_2^{t_n} M(t)de_{\ln,1}(t)}\dsps\leq -\left[M^{\phantom{\frac pp}}\!\!\!\!\!(t)e_{\ln,1}(t)\right]_2^{t_n}
+\int_2^{t_n} \frac{M'(t)}{t\ln^2t}dt\\[4mm]
\phantom{\dsps-\int_2^{t_n} M(t)de_{\ln,1}(t)}\dsps\leq -\left[M^{\phantom{\frac pp}}\!\!\!\!\!(t)e_{\ln,1}(t)\right]_2^{t_n}
+\left[\tfrac{M(t)}{t\ln^2t}\right]_2^{t_n}+\int_2^{t_n} \frac{M(t)}{t^2\ln^2t}\left(1+\frac2{\ln t}\right)dt,
\EA$$
where $\{t_n\}$ is a sequence tending to infinity and such that 
$$\dsps\lim_{n\to\infty}\int_2^{t_n} \frac{M(t_n)}{t_n\ln^2t_n}dt=0.
$$
Hence 
$$\lim_{n\to\infty}\int_2^{t_n}M(t)de_{\ln,1}(t)\leq \left(1+\frac2{\ln 2}\right)\int_2^{\infty} \frac{M(t)}{t^2\ln^2t}dt<\infty.
$$
This implies 1. The proof of 2 and 3 is similar.  \hfill$\Box$\medskip
%%%%%%%%%%%r%%%%%%%%%%%%%%%%%%%%%%%%%%%%%%%%%
%%%%%%%%%%%%%%%%%%%%%%%%%%f%%%%%%%%%%%%%%%%%%%%%%%%%%%%%%%%%%%%%%%%%%%%%%%%%%%%%
%%%%%%%%%%%r%%%%%%%%%%%%%%%%%%%%%%%%%%%%%%%%%
%%%%%%%%%%%%%%%%%%%%%%%%%%f%%%%%%%%%%%%%%%%%%%%%%%%%%%%%%%%%%%%%%%%%%%%%%%%%%%%%%%%%%%%%%%%%%%%%%%%%%%%%%%%%%%%%%%%%%%%%%%%%%%%%%%%%%%%%%%%%%%%%%%%%%%%%%%%%%%%%
\begin{lemma}\label{gradPhi}
Let $N\geq 3$, then there holds
\begin{equation}\label{EI-5}
|\nabla \Phi_{1}(x)|\leq c\max\left\{|x|^{1-N},\, |x|^{-\frac{N-1}{2}}\right\}\quad\text{for all  }x\in\BBR^N\setminus\{0\}.
 \end{equation}  
 More precisely, since $\Phi_1$ is radial, 
 \begin{equation}\label{EI-6}
\Phi'_{1}(r)= \left\{\BA {lll}cr^{1-N}(1+o(1))\qquad&\text{as }\, r\to 0,\\[2mm]
c\frac{\sin\left(r-\frac{(N+1)\gp}{4}\right)}{r^{\frac {N-1}2}}(1+O\left(\tfrac 1r)\right)\qquad&\text{as }\, r\to +\infty.
\EA\right.
 \end{equation}  
\end{lemma}
\noindent\proof If $\Phi_1(x)=\Phi_1(|x|)$ is the solution of the Helmholtz equation in $\BBR^N$, then $Z=\Phi'_1$ satisfies
$$Z''+\frac{N-1}{r}Z'+(1-\frac {N-1}{r^2})Z=0\quad\text{on }(0,\infty).
$$
Set $z(r)=r^{\frac N2-1}Z(r)$, we have that 
$$z''+\frac{1}{r}z'+\left(1-\frac {N^2}{4r^2}\right)z=0\quad\text{on }(0,\infty),
$$
which is a classical Bessel equation. Using the asymptotic expressions of Hankel functions we obtain 
$$Z(r)=cr^{1-N}(1+o(1))\quad\text{as }r\to 0
$$
and
$$Z(r)=c\frac{\sin\left(r-\frac{(N+1)\gp}{4}\right)}{r^{\frac {N-1}2}}\left(1+O\left(\frac 1r\right)\right)\quad\text{as }r\to +\infty.
$$
Thus 
\begin{equation}\label{EI-6}
\Phi'_{1}(r)= \left\{\BA {lll}cr^{1-N}(1+o(1))\qquad&\text{as }r\to 0,\\[4mm]\dsps
c\frac{\sin\left(r-\frac{(N+1)\gp}{4}\right)}{r^{\frac {N-1}2}}\left(1+O\left(\frac 1r\right)\right)\qquad&\text{as }r\to +\infty,
\EA\right.
 \end{equation}  
which ends the proof.
\hfill$\Box$\medskip

%%%%%%%%%%%%%%%%%%%%%%%%%%%%%%%%%%%%%%%%%%%%%%%%%%%%%%%%%%%%%%%%%%%%%%%%%%%%%%%%%%%%%%%%%%%%%%%%%%%%%%%%%%%%%%%%%%%%%%%%%%%%%%%%%%%%%%%%%%%%%%%%%%%%%%%%%%%%%%%%
 \begin{theorem}\label{P3}
 Let $N\geq 4$, then there exists $c>0$ such that 
 \begin{equation}\label{EI-7}
\big|\nabla\Phi_{\ln}(x)\big|\leq c\max\left\{\frac{|x|^{-N-1}}{1+\ln^2|x|},\frac{1}{\sqrt{1+\ln^2|x|}}\right\} \quad\text{for all }x\in\BBR^N\setminus\{0\}.
 \end{equation} 
\end{theorem} 
\noindent\proof We recall that 
$$\Phi_{\ln}(x)=u_*(x)+v_1(x)=\int_0^1\CP_0(t)|x|^{2t-N}dt+\int_0^1\CP_0(t)(|.|^{2t-N}\ast\Phi_1(.))(x)dt.$$
{\it Step 1: Estimate of $\nabla u_*$}. There holds 
$$\BA{lll}\dsps
|\nabla u_*(x)|= \int_0^1\CP_0(t)(N-2t)|x|^{2t-N-1}dt\\[4mm]
\phantom{|\nabla u_*(x)|}
\dsps \leq c|x|^{-1-N}\int_0^1t|x|^{2t}dt 
\\[4mm]
\phantom{|\nabla u_*(x)|}
\dsps =\frac{c|x|^{1-N}}{2\ln |x|}-\frac{c|x|^{1-N}}{4\ln^2 |x|}+\frac{cx|^{-1-N}}{4\ln^2 |x|}.
\EA$$
if $|x|\notin\{0, 1\}$. One has
 \begin{equation}\label{EI-71}\BA{lll}\dsps
|\nabla u_*(x)|\leq \frac{cx|^{-1-N}}{4\ln^2 |x|}(1+o(1))\qquad\text{as }x\to 0
\EA \end{equation}
and
 \begin{equation}\label{EI-72}\BA{lll}\dsps
|\nabla u_*(x)|\leq \frac{cx|^{1-N}}{2\ln |x|}(1+o(1))\qquad\text{as }|x|\to \infty.
\EA \end{equation}
By standard asymptotic expansion near $|x|=1$,
 \begin{equation}\label{EI-73}\BA{lll}\dsps
\lim_{|x|\to 1}|\nabla u_*(x)|\leq c\lim_{|x|\to 1}\left(\frac{c|x|^{1-N}}{2\ln |x|}-\frac{c|x|^{1-N}}{4\ln^2 |x|}+\frac{cx|^{-1-N}}{4\ln^2 |x|}\right)=c.
\EA \end{equation}
when $|x|=1$. Therefore $x\mapsto |\nabla u_*(x)|$ is uniformly bounded on $B_\ge^c$ for any $\ge>0$. \smallskip

\noindent {\it Step 2: Estimate of $\nabla v_1$}.
We write 
$$\BA{lll}\dsps v_1(x)=\int_0^1\CP_0(t)(|.|^{2t-N}\ast\Phi_1(.))(x)dt\\[4mm]
\phantom{\dsps v_1(x)}\dsps =\int_0^1\CP_0(t)\int_{B_1}|x-y|^{2t-N}\Phi_1(y))dydt+\int_0^1\CP_0(t)\int_{B^c_1}|x-y|^{2t-N}\Phi_1(y))dydt,
\EA$$
thus
$$\BA{lll}\dsps\nabla v_1(x)=\int_0^1\CP_0(t)(|.|^{2t-N}\ast\nabla \Phi_1(.))(x)dt.
\EA
$$
Setting ${\bf e}_y=\frac y{|y|}$, we have
$$\BA{lll}\dsps
|\nabla v_1(x)|\leq c\int_0^1\int_{B_1}\CP_0(t)|x-y|^{2t-N}|y|^{1-N}dydt+\left|\int_0^1\CP_0(t)\int_{B^c_1}|x-y|^{2t-N}\langle\nabla\Phi_1(y),{\bf e}_y\rangle dydt\right|\\[4mm]
\phantom{|\nabla v_1(x)|}=A_6(x)+A_7(x).
\EA
$$
\noindent {\it Step 2-1: Estimate of $A_6(x)$.} 
We have
$$\BA{lll}\dsps A_6(x)\leq c'|x|^{1-N}\int_0^1t|x|^{2t}\int_{B_\frac1{|x|}}|{\bf e}_N-z|^{2t-N}|z|^{1-N}dzdt
\\[4mm]
\phantom{A_6(x)}\dsps
=c'|x|^{1-N}\int_{B_\frac1{|x|}}|{\bf e}_N-z|^{-N}|z|^{1-N}\left(\int_0^1t|x|^{2t}|{\bf e}_N-z|^{2t}dt\right)dz\\[4mm]
\phantom{A_6(x)}\dsps
=c'|x|^{1-N}\int_{B_\frac1{|x|}}|{\bf e}_N-z|^{-N}|z|^{1-N}\left(\frac{|x|^2|{\bf e}_N-z|^2}{2\ln(|x||{\bf e}_N-z|)}+\frac{1-|x|^2|{\bf e}_N-z|^2}{4\ln^2(|x||{\bf e}_N-z|)}\right)dz.
\EA$$
Set $\gz=|x|(z-{\bf e}_N)$, then 
\begin{equation}\label{A32}\BA{lll}\dsps 
A_6(x)\leq c'|x|^{1-N}\int_{B_1(-|x|{\bf e}_N)}\!\!\!|\gz|^{-N}\left|{\bf e}_N+|x|^{-1}\gz\right|^{1-N}\left(\frac{|\gz|^2}{2\ln(|\gz|)}+\frac{1-|\gz|^2}{4\ln^2(|\gz|)}\right)d\gz.
\EA\end{equation}
We set $|\gz|=1+m$, then by standard expansion we have
$$\frac{|\gz|^2}{2\ln(|\gz|)}+\frac{1-|\gz|^2}{4\ln^2(|\gz|)}=2+O(||\gz|-1|)\quad\text{as }|\gz|\to 1.
$$
{\it Step 2-1-(i): Estimate of $A_6(x)$ when $x\to 0$}. If $|x|\leq \frac12$, then
$$\BA{lll}\dsps 
A_6(x)\leq 2^Nc'|x|^{1-N}\int_{B_1(-|x|{\bf e}_N)}\left|{\bf e}_N+|x|^{-1}\gz\right|^{1-N}\left(\frac{|\gz|^2}{2\ln(|\gz|)}+\frac{1-|\gz|^2}{4\ln^2(|\gz|)}\right)d\gz\\[4mm]
\phantom{A(x)}\dsps\leq 2^Nc'\int_{B_1(-|x|{\bf e}_N)}||x|{\bf e}_N+\gz|^{1-N}(2+O(||\gz|-1|))d\gz,
\EA$$
therefore
\begin{equation}\label{A33}\dsps
\limsup_{x\to 0}A_6(x)\leq 2^Nc'\int_{B_1}|\gz|^{1-N}(2+O(||\gz|-1|))d\gz.
\end{equation}
{\it Step 2-1-(ii): Estimate of $A_6(x)$ when $|x|\to \infty$}. We have directly from (\ref{A32})
\begin{equation}\label{A34}\BA{lll}\dsps 
A_6(x)\leq c'\frac{|x|^{2-N}}{2\ln|x|}\left(\int_{B_1(-|x|{\bf e}_N)}||x|{\bf e}_N+\gz|^{1-N}d\gz\right)(1+o(1))\\[4mm]
\phantom{A(x)}\dsps=c'\frac{|x|^{2-N}}{2\ln|x|}\left(\int_{B_1}|z|^{1-N}dz\right)(1+o(1)).
\EA\end{equation}
Combining (\ref{A33}) and (\ref{A34}) we deduce that 
\begin{equation}\label{A35}\BA{lll}\dsps 
A_6(x)\leq c\min\left\{1,\frac{|x|^{2-N}}{\ln(1+|x|)}\right\}\quad\text{for all }x\in\BBR^N\setminus\{0\}.
\EA\end{equation}

\noindent \noindent {\it Step 2-2: Estimate of $A_7(x)$.} Using (\ref{EI-6}) we have as in (\ref{hel18})
$$\BA{lll}\dsps 
A_7(x)\leq\left|\int_0^1\CP_0(t)\int_{B^c_1}|x-y|^{2t-N}|y|^{\frac{1-N}{2}}\sin\left(|y|-\frac{(N-1)\gp}{4}\right)dydt\right|\\[5mm]
\phantom{-----}\dsps
+\int_0^1\CP_0(t)\int_{B^c_1}|x-y|^{2t-N}|y|^{-\frac{N+1}{2}}dydt\\[5mm]
\phantom{B(x)}\dsps
\leq K_1(x)+K_2(x).
\EA$$
Clearly
$$\BA{lll}\dsps K_2(x)
\leq c\int_{B^c_1}|x-y|^{-N}|y|^{-\frac{N+1}{2}}\int_0^1t|x-y|^{2t}dtdy\\[4mm]
\phantom{K_2(x)}
\dsps \leq c\int_{B^c_1}|x-y|^{-N}|y|^{-\frac{N+1}{2}}\left(\frac{|x-y|^2}{2\ln|x-y|}+\frac{1-|x-y|^2}{4\ln^2|x-y|}\right)dy.
\EA$$
We have
$$K_2(x)=c|x|^{-\frac{N+1}{2}}\int_{B^c_{\frac1{|x|}}}|{\bf e}_N-z|^{-N}|z|^{-\frac{N+1}{2}}\left(\frac{|x|^2|{\bf e}_N-z|^2}{2(\ln|x|+\ln|{\bf e}_N-z|)}
+\frac{1-|x|^2|{\bf e}_N-z|^2}{4(\ln|x|+\ln|{\bf e}_N-z|)^2}\right)dz.$$
{\it When $|x|$ is small}. The integral defining $K_2(x)$ is convergent since $N\geq 4$ and we have by integration by part
$$%\begin{equation}\label{A36}
K_2(x)=\frac{2cc_{_N}|x|^2}{(N+1)|\ln|x||}(1+o(1))\qquad\text{as }x\to 0.
$$%\end{equation} 
The term $K_1(x)$ has already been partially estimated in Lemma \ref{est-l}. By (\ref{AS5}) we have
$$\left|\int_{B_1^c}|x-y|^{2t-N}|y|^{\frac{1-N}{2}}\sin\left(|y|-\frac{(N-1)\gp}{4}\right)dy\right|\leq c\quad\text{if }|x|\leq 1,
$$
uniformly with respect to $\gt$. This implies that 
$$\left|\int_0^1\CP_0(t)\int_{B^c_1}|x-y|^{2t-N}|y|^{\frac{1-N}{2}}\sin\left(|y|-\frac{(N-1)\gp}{4}\right)dydt\right|\leq c,
$$
and finally 
$$%\begin{equation}\label{A36}
\BA{lll}\dsps 
A_7(x)\leq c\quad\text{for all }x\in B_1\setminus\{0\}.
\EA
$$%\end{equation}
{\it When $|x|$ is large}. we use the results of Lemma \ref{est-l}, in{\it Step II-Estimate of $\widetilde I\!I_\ga$}.
By (\ref{hel18+}) we have
$$\BA{lll}\dsps\int_0^1\CP_0(t)\int_{B^c_1}|x-y|^{2t-N}|y|^{-\frac{N+1}{2}}dydt\leq \int_0^1t|x|^{2t-\frac{N+1}{2}}\int_{B^c_1}\left(|{\bf e}_N-z|^{2t-N}|z|^{-\frac{N+1}{2}}\right)dz\\[4mm]
\phantom{\dsps\int_0^1\int_{B^c_1}|x-y|^{2t-N}|y|^{-\frac{N+1}{2}}dydt}\dsps
\leq |x|^{-\frac{N+1}{2}}\int_{B^c_{\frac 1{|x|}}}|{\bf e}_N-z|^{-N}|z|^{-\frac{N+1}{2}}\int_0^1t|x|^{2t}|{\bf e}_N-z|^{2t}dtdz\\[4mm]\phantom{-----------}\leq 
\dsps |x|^{-\frac{N+1}{2}}\int_{B^c_1(-|x|{\bf e}_N)}|\gz|^{-N}\left|{\bf e}_N+|x|^{-1}\gz\right|^{-\frac{N+1}{2}}\left(\frac{|\gz|^2}{2\ln|\gz|}+\frac{1-|\gz|^2}{4\ln^2|\gz|}\right) d\gz,
\EA$$
since the last term has already been estimated in (\ref{A32}). Therefore
$$\BA{lll}\dsps\int_0^1\CP_0(t)\int_{B^c_1}|x-y|^{2t-N}|y|^{-\frac{N+1}{2}}dydt\\[4mm]\phantom{------}\leq 
\dsps |x|^{-\frac{N+1}{2}}\int_{B^c_1(-|x|{\bf e}_N)}|\gz|^{-N}\left|{\bf e}_N+|x|^{-1}\gz\right|^{-\frac{N+1}{2}}\left(\frac{|\gz|^2}{2\ln|\gz|}+\frac{1-|\gz|^2}{4\ln^2|\gz|}\right) d\gz:=M(x).
\EA$$
The above integral is convergent, but the limit when $|x|\to\infty$ is a little more tricky to obtain:
Let $a>2$, for $|x|>2a$, we write
$$\BA{lll}\dsps|x|^{-\frac{N+1}{2}}\int_{B^c_1(-|x|{\bf e}_N)}|\gz|^{-N}\left|{\bf e}_N+|x|^{-1}\gz\right|^{-\frac{N+1}{2}}\left(\frac{|\gz|^2}{2\ln|\gz|}+\frac{1-|\gz|^2}{4\ln^2|\gz|}\right) d\gz\\[4mm]
\phantom{------}\dsps=
|x|^{-\frac{N+1}{2}}\left(\int_{B_a}+\int_{B^c_1(-|x|{\bf e}_N)\cap B_a^c}\right)|\gz|^{-N}\left|{\bf e}_N+|x|^{-1}\gz\right|^{-\frac{N+1}{2}}\left(\frac{|\gz|^2}{2\ln|\gz|}+\frac{1-|\gz|^2}{4\ln^2|\gz|}\right) d\gz\\[4mm]
\phantom{------}\dsps=M_1(x)+M_2(x).
\EA$$
Then
$$\BA{lll}\dsps
M_1(x)
\leq |x|^{-\frac{N+1}{2}}\int_{B_a}|\gz|^{-N}\left|{\bf e}_N+|x|^{-1}\gz\right|^{-\frac{N+1}{2}}\left(\frac{|\gz|^2}{2\ln|\gz|}+\frac{1-|\gz|^2}{4\ln^2|\gz|}\right) d\gz
\\[4mm]\phantom{M_1(x)}\dsps\leq |x|^{-\frac{N+1}{2}}2^{\frac {N+1}{2}}\int_{B_a}|\gz|^{-N}\left(\frac{|\gz|^2}{2\ln|\gz|}+\frac{1}{4\ln^2|\gz|}\right) d\gz
\\[4mm]\phantom{M_1(x)}\dsps\leq c(a)|x|^{-\frac{N+1}{2}}.
\EA$$
Next, since $a>2$, 
$$\BA{lll}\dsps
M_2(x) \leq |x|^{-\frac{N+1}{2}}\int_{B_a^c}\left|{\bf e}_N+|x|^{-1}\gz\right|^{-\frac{N+1}{2}}\frac{|\gz|^{2-N}}{\ln|\gz|} d\gz
\\[4mm]\phantom{M_2(x)}\dsps
=\int_{B_\frac{a}{|x|}^c}\left|{\bf e}_N+t\right|^{-\frac{N+1}{2}}\frac{|t|^{2-N}}{\ln|t|+\ln |x|} d
  \leq \frac{c'(a)}{\ln|x|}.
\EA$$
Therefore, 
\begin{equation}\label{A37}\BA{lll}\dsps 
M(x)\leq \frac{c}{1+\ln|x|}\qquad\text{for all }x\in B^c_1.
\EA\end{equation}
Combining (\ref{EI-71}), (\ref{EI-72}), (\ref{EI-73}), (\ref{A35}) and (\ref{A37}), we obtain (\ref{EI-7}),
which ends the proof.
\hfill$\Box$\medskip

\noindent {\bf Remark} {The restriction $N\geq 4$ appears to be purely technical, but our technique which is based upon the delicate use of a series of semi-convergent integrals appears very sensitive to the dimension.  However we conjecture that the  estimates holds whenever $N\geq 3$}.\medskip

     %%%%%%%%%%%%%%%%%%%%%%%%%%%%%%%%%%%%%%%%%%%%%%%%%%%%%%%%%%%%%%%%%%%%%%%%%%%%%%%%%%%%%%%%%%%%%%%%%%%%%%%%%%%%%%%%%%%%%%%%%%%%%%%%%%%%%%%%%%%%%%%%%%%%%%%%%%%%%%%%%%%%%%%%%%%%%%%%%%%%%%%%%%%%%%%%%%%%%%%%%%%%%%%%%%%%%%%%%%%%%%%%%%%%%%%%%%%%%%%%%%%%%%%%%%%%%%%%%%%%%%%%%%%%%%%%%%%%%%%%%%%
%\color{black}{
 %%%%%%%%%%%%%%%%%%%%%%%%%%%%%%%%%%%%%%%%%%%%%%%%%%%%%%%%%%%%%%%%%%%%%%%%%%%%%%%%%%%%%%%%%%%%%%%%%%%%%%%%%%%%%%%%%%%%%%%%%%%%%%%%%%%%%%%%%%%%%%%%%%%%%%%%%%%%%%%%%%%%%%%%%%%%%%%%%%%%%%%%%%%%%%%%%%%%%%%%%%%%%%%%%%%%
 We recall that $\Phi_{\ln}=u_*+v_1$ with the notations of Theorem \ref{th 8.1}

\begin{lemma}\label{reg1} 
Let $N\geq 4$ and  $f\in L^1(\R^N)$, then $v_1\ast f$ is Lipschitz continuous and 
 \begin{equation}\label{R1-sec5-0}
 \big|v_1\ast f(x)-v_1\ast f(y)\big|\leq c|\norm f_{L^1}|x-y|.
  \end{equation}  
\end{lemma} 

\noindent\proof By Theorem \ref{P3}-Step 2 we have
$$|\nabla v_1(x)|\leq c'{\bf 1}_{B_1}(x)+\frac{c}{1+\ln |x|}{\bf 1}_{B^c_1}(x).
$$
The proof follows by Young's inequality for convolution. \hfill$\Box$\medskip

%%%%%%%%%%%%%%%%%%%%%
The regularity concerning $u_*\ast f$ is more involved. 
\begin{lemma}\label{reg2} Let $\gth>1$, 
and  $f\in L^\infty(\R^N;(1+|y|)^\gth)$. Then $u_*\ast f$ is Dini continuous and for any $R>0$ there exists 
for some constant $c=c(N,R)>0$,
 \begin{equation}\label{R1-sec5}
\left|u_*\ast f(x)-u_*\ast f(x')\right|\leq \frac{c\norm f_{\infty,\gth}}{1+\ln^2(|x-x'|) }\quad\text{for all }x,x'\in B_R,
  \end{equation} 
  where 
$\norm f_{\infty,\gth}$ is defined in (\ref{AAD}).
\end{lemma}
%%%%%%%%%%%%%%%%%%%%%
\noindent\proof Let $R>0$, we have
$$\BA{lll}\dsps
\left|u_*\ast f(x)-u_*\ast f(x')\right|=\left|\int_0^1\CP_0(t)\int_{\BBR^N}(|x-y|^{2t-N}-|x'-y|^{2t-N})f(y)dydt\right|
\\[4mm]
\phantom{\left|u_*\ast f(x)-u_*\ast f(x')\right|}\dsps\leq 
 \norm f_{\infty,\gth}\int_0^1\CP_0(t)\int_{\BBR^N}\left||x-y|^{2t-N}-|x'-y|^{2t-N}\right|(1+|y|)^{-\gth}dydt
\\[4mm]
\phantom{\left|u_*\ast f(x)-u_*\ast f(x')\right|}\dsps
\leq \norm f_{\infty,\gth}\int_0^1\CP_0(t)\int_{\BBR^N}\left||x-x'+y|^{2t-N}-|y|^{2t-N}\right|(1+|y-x'|)^{-\gth}dydt\\[4mm]
\phantom{\left|u_*\ast f(x)-u_*\ast f(x')\right|}\dsps
\leq c \norm f_{\infty,\gth}\int_0^1t\int_{B^c_R}\left||x-x'+y|^{2t-N}-|y|^{2t-N}\right|(1+|y-x'|)^{-\gth}dydt\\[4mm]
\phantom{\left|u_*\ast f(x)-u_*\ast f(x')\right|---}\dsps
+ c \norm f_{\infty,\gth}\int_0^1t\int_{B_R}\left||x-x'+y|^{2t-N}-|y|^{2t-N}\right|(1+|y-x'|)^{-\gth}dydt\\[4mm]
\phantom{\left|u_*\ast f(x)-u_*\ast f(x')\right|}\dsps
\leq c \norm f_{\infty,\gth}\left(A(|x-x'|,x')+B(|x-x'|,x')^{\phantom{g^r}}\!\!\!\!\!\right),
\EA$$
where we have denoted 
$$A(|x-x'|,x')=\int_0^1t\int_{B^c_R}\left||x-x'+y|^{2t-N}-|y|^{2t-N}\right|(1+|y-x'|)^{-\gth}dydt,
$$
and
$$B(|x-x'|,x')=\int_0^1t\int_{B_R}\left||x-x'+y|^{2t-N}-|y|^{2t-N}\right|(1+|y-x'|)^{-\gth}dydt.
$$
We suppose that $\max\{|x|,|x'|\}\leq \frac R8$.
$$\BA{lll}\dsps A(|x-x'|,x')\leq c\int_0^1t\int_{B^c_R}\left||x-x'+y|^{2t-N}-|y|^{2t-N}\right|(1+|y|)^{-\gth}dydt\\[4mm]
\phantom{A(|x-x'|,x')}\dsps\leq c\int_0^1t|x-x'|^{2t-\gth}\int_{B^c_{\frac R{|x-x'|}}}\left||{\bf e}_{x-x'}-z|^{2t-N}-|z|^{2t-N}\right|(|x-x'|^{-1}+|z|)^{-\gth}dydt\\[4mm]
\phantom{A(|x-x'|,x')}\dsps\leq c\int_0^1|tx-x'|^{2t-\gth}\int_{B^c_{\frac R{|x-x'|}}}\left||{\bf e}_{x-x'}-z|^{2t-N}-|z|^{2t-N}\right|(|x-x'|^{-1}+|z|)^{-\gth}dydt\\[4mm]
\phantom{A(,x')}\dsps\leq c\int_0^1t|x-x'|^{2t-\gth}\int_{B^c_{\frac R{|x-x'|}}}|z|^{2t-N}\left|||z|^{-1}{\bf e}_{x-x'}-{\bf e}_{z}|^{2t-N}-1\right|(|x-x'|^{-1}+|z|)^{-\gth}dydt.
\EA$$
Since by standard Taylor expansion
$$\left|||z|^{-1}{\bf e}_{x-x'}-{\bf e}_{z}|^{2t-N}-1\right|\leq \frac{N-2t}{|z|}(1+o(1)) \quad\text{as }|z|\to\infty,
$$
we deduce that
 \begin{equation}\label{R1-1}\BA{lll}\dsps A(|x-x'|,x')\leq c'\int_0^1t|x-x'|^{2t-\gth}\int_{\frac R{|x-x'|}}^\infty r^{2t-2-\gth}drdt\\[4mm]
\phantom{A(|x-x'|,x')}\dsps \leq c''|x-x'|R^{-1-\gth}\int_0^1tR^{2t}dt\\[4mm]
\phantom{A(|x-x'|,x')}\dsps \leq c''\left(\frac{R^2}{2\ln R}+\frac{1-R^2}{4\ln^2 R}\right)\frac{|x-x'|}{R^{1+\gth}}.
\EA\end{equation}
Since $|x-x'|\leq \frac R4$, we have for $a<4$, 
$$\BA{lll}\dsps B(|x-x'|,x')\leq \int_0^1t\int_{B_R}\left||x-x'+y|^{2t-N}-|y|^{2t-N}\right|(1+|y-x'|)^{-\gth}dydt\\[4mm]
\phantom{B(|x-x'|,x')}\dsps\leq \int_0^1t|x-x'|^{2t}\int_{B_{\frac R{|x-x'|}}}\left||{\bf e}_{x-x'}-z|^{2t-N}-|z|^{2t-N}\right|dzdt.
\\[4mm]
\phantom{-}\dsps\leq 2\int_0^1t|x-x'|^{2t}\int_{B_{a}}|z|^{2t-N}dzdt+\int_0^1t|x-x'|^{2t}\int_{B_{\frac R{|x-x'|}}\setminus B_a}\left||{\bf e}_{x-x'}-z|^{2t-N}-|z|^{2t-N}\right|dzdt
\\[4mm]
\phantom{}\dsps \leq A_8+A_9.
\EA$$
Then
 \begin{equation}\label{R1-2}\BA{lll}\dsps A_8\leq c_{_N}\int_0^1t|x-x'|^{2t}a^{2t}dt=c_{_N}\left(\frac{|x-x'|^{2}a^2}{2\ln(|x-x'|a)}+\frac{1-|x-x'|^{2}a^2}{4\ln^2(|x-x'|a)}\right)\\[4mm]
 \phantom{B_1}\dsps\leq \frac{2c_{_N}}{\ln^2(|x-x'|a)}(1+o(|x-x'|),
\EA\end{equation}
and 
 \begin{equation}\label{R1-3}\BA{lll}\dsps 
 A_9\leq \int_0^1t|x-x'|^{2t}\int_{B_{\frac R{|x-x'|}}\setminus B_a}\left|||z|^{-1}{\bf e}_{x-x'}-{\bf e}_{z}|^{2t-N}-1\right||z|^{2t-N}dzdt\\
 \phantom{B_2}\dsps\leq c\int_0^1t|x-x'|^{2t}\int_a^{\frac R{|x-x'|}}r^{2t-2}drdt\\[4mm]
 \phantom{B_2}\dsps
\leq c\int_0^1|x-x'|^{2t}\left(\frac{R^{2t-1}}{|x-x'|^{2t-1}}-a^{2t-1}\right)\frac{tdt}{2t-1}.
\EA\end{equation}
We have that
$$\int_0^1|x-x'|^{2t}\left(\frac{R^{2t-1}}{|x-x'|^{2t-1}}-a^{2t-1}\right)\frac{tdt}{2t-1}=|x-x'|\int_0^1\left(R^{2t-1}-(a|x-x'|)^{2t-1}\right)\frac{tdt}{2t-1}
$$
and
$$\BA{lll}\dsps R^{2t-1}-(a|x-x'|)^{2t-1}=e^{(2t-1)\ln R}-e^{(2t-1)(\ln a|x-x'|)}\\[2mm]
\phantom{R^{2t-1}-(a|x-x'|)^{2t-1}}\dsps =(2t-1) \ln \left(\frac{R}{a|x-x'|}\right)+O\left((2t-1)^2(\ln^2R+\ln^2a|x-x'|)\right).
\EA$$
By integration,
$$\left|\int_0^1\left(R^{2t-1}-(a|x-x'|)^{2t-1}\right)\frac{tdt}{2t-1}\right|\leq \ln \left(\frac{R}{a|x-x'|}\right)+O\left(\ln^2R+\ln^2a|x-x'|)\right).
$$
This finally yields to 
 \begin{equation}\label{R1-4}
 \BA{lll}\dsps A_9\leq c|x-x'|\left[\ln \left(\frac{R}{a|x-x'|}\right)+O\left(\ln^2R+\ln^2a|x-x'|)\right)\right].
\EA\end{equation}
We notice that the leading term when $x-x'\to 0$ is coming from the estimate of $B_1$ and it gives the Dini modulus of continuity  expressed by (\ref{R1-sec5}).
This ends the proof.\hfill$\Box$\medskip

\noindent{\it Proof of Theorem \ref{cr 5.1}. } It follows directly from Lemmas \ref{reg1} and \ref{reg2}.
\hfill$\Box$

 \bigskip\bigskip

   \noindent{\bf \small Acknowledgements:}    {\footnotesize 
H. Chen is supported by  NNSF of China (Nos. 12001252, 12361043) 
and Jiangxi Natural Science Foundation (No. 20232ACB211001).  The authors are much grateful to the anonymous referee for the very comprehensive review of their work. }

%\newpage

%\subsection{Regularity}

%{\color{blue} {\it Regularity of weak solution: } Since $h \in \cC^1(\bar\Omega\setminus\{0\})$, given any 
 %point $x_0\not=0$, it is standard method to show that 
 %$u_h$ is $C^1_{loc}(\Omega\setminus\{0\})$. Therefore, $u_h$ is a classical solution of (\ref{eq 7.3}) and $u_h$ is continuous up to the boundary see the proof of Theorem 5.4 in \cite{CT}. }

\end{document}